\documentclass[letterpaper,11pt]{article}

\usepackage{ucs}
\usepackage[utf8x]{inputenc}
\usepackage{graphicx}
\usepackage{amsfonts}
\usepackage{dsfont}
\usepackage{amssymb}
\usepackage{amsmath}
\usepackage{amsthm}
\usepackage{enumitem}
\usepackage{stmaryrd}
\usepackage{fullpage}
\usepackage{ifthen}
\usepackage{subfigure}
\usepackage{epic}
\usepackage{authblk}
\usepackage{textcomp}
\usepackage{mathrsfs}
\usepackage{bm}
\usepackage[small]{caption}
\usepackage{mathtools}
\usepackage{multirow}

\graphicspath{{Figures_manuscript/}}

\usepackage{color}
\usepackage{titlesec}

\graphicspath{{eps/}{pdf/}}
\newcommand{\bqq}{\begin{equation}}
\newcommand{\eqq}{\end{equation}}
\newcommand{\bqs}{\begin{equation*}}
\newcommand{\eqs}{\end{equation*}}

\newcommand{\calL}{\mathcal{L}}

\newcommand{\C}{\mathbb{C}}

\newcommand{\R}{\mathbb{R}} 

\newcommand{\W}{\mathbb{W}}

\renewcommand{\O}{\mathcal{O}}

\newcommand{\RPBP}{(\lambda-\lambda^\textnormal{rp})\sqrt{\zeta(\lambda-\lambda_v^\textnormal{bp})}}
\newcommand{\BP}{\sqrt{\zeta(\lambda-\lambda_v^\textnormal{bp})}}

\newcommand{\bI}{\mathbf{I}}
\newcommand{\boldJ}{\mathbf{J}}

\newcommand{\bG}{\mathbf{G}}

\newcommand{\boldH}{\mathbf{H}}
\newcommand{\bJ}{\mathbf{J}}

\newcommand{\cL}{\mathcal{L}}

\newcommand{\cG}{\mathcal{G}}
\newcommand{\X}{\mathcal{X}}
\newcommand{\Y}{\mathcal{Y}}
\newcommand{\vp}{\varphi}
\newcommand{\md}{\,\mathrm{d}}
\newcommand{\me}{\mathrm{e}}
\newcommand{\mbi}{\mathbf{i}}

\renewcommand{\Re}{\mathrm{Re}\,}
\renewcommand{\Im}{\mathrm{Im}\,}

\def\calL{\mathcal{L}}
\def\Sess{\Sigma_{\textnormal{ess}}}
\def\Sabs{\Sigma_{\textnormal{abs}}}

\newcommand{\rme}{\textnormal{e}}

\newtheorem{lem}{Lemma}[section]
\newtheorem{thm}{Theorem}
\newtheorem{prop}[lem]{Proposition}
\newtheorem{conj}{Conjecture}

\newtheorem{rmk}[lem]{Remark}
\newtheorem{defi}[lem]{Definition}

\newenvironment{Proof}[1][.]%
 {\begin{trivlist}\item[]\textbf{Proof#1 }}%
 {\hspace*{\fill}$\rule{0.3\baselineskip}{0.35\baselineskip}$\end{trivlist}}
 
  {\begin{trivlist}\item[]{\bf Hypothesis #1 }\em}{\end{trivlist}}

\numberwithin{equation}{section}

\title{Invasion into remnant instability: a case study of front dynamics}

\author[1]{Gr\'egory Faye\footnote{Corresponding author email: \texttt{gregory.faye@math.univ-toulouse.fr}}}
\author[2]{Matt Holzer}%\footnote{email: \texttt{mholzer@gmu.edu}}}
\author[3]{Arnd Scheel}%\footnote{email: \texttt{scheel@umn.edu}}}
\author[4]{Lars Siemer}%\footnote{email: \texttt{lars.siemer@uni-bremen.de}}}
\affil[1]{\small CNRS, UMR 5219, Institut de Math\'ematiques de Toulouse, 31062 Toulouse Cedex, France}
\affil[2]{\small Department of Mathematical Sciences, George Mason University, Fairfax, VA, USA}
\affil[3]{\small School of Mathematics, University of Minnesota, Minneapolis, MN, USA }
\affil[4]{\small  Department of Mathematics, University of Bremen, 28359 Bremen, Germany}

\begin{document}

\renewcommand{\theenumi}{(\roman{enumi})}%
\maketitle

\begin{abstract} We study the invasion of an unstable state by a propagating front in a peculiar but generic situation where the invasion process exhibits a remnant instability. Here, remnant instability refers to the fact that the spatially constant invaded state is linearly unstable in any exponentially weighted space in a frame  moving with the linear invasion speed. Our main result is the nonlinear asymptotic stability of the selected invasion front for a prototypical model coupling spatio-temporal oscillations and monotone dynamics. We establish stability through a decomposition of the perturbation into two pieces: one that is bounded in the weighted space and a second that is unbounded in the
weighted space but which converges uniformly to zero in the unweighted space at an exponential
rate.  Interestingly, long-time numerical simulations reveal an apparent instability in some cases. We exhibit how this instability is caused by round-off errors that introduce linear resonant coupling of otherwise non-resonant linear modes, and we determine the accelerated invasion speed.

\end{abstract}


% {\noindent \bf MSC numbers:}\\

\section{Introduction}
The stability analysis of fronts invading unstable states is complicated by the presence of unstable essential spectrum. This unstable spectrum is a consequence of the instability of the homogeneous state that the front is propagating into. Perturbations of the front that are placed sufficiently far ahead of the front interface will grow in norm at an exponential temporal rate from the homogeneous state. Stability typically is recovered by compensating temporal growth with exponential spatial decay of allowable perturbations away from the front interface. Technically, this amounts to viewing the stability problem in an exponentially weighted function space. The hope is that some weight can be chosen so that the essential spectrum is stabilized in the comoving frame of the front and the stability analysis can then proceed in a fashion similar to that of fronts connecting stable states; see for example~\cite{sattinger}. In this paper we study a model system where no such stabilization is possible and investigate the consequences for both the stability analysis of the front and the associated dynamics of the full partial differential equation.

As a specific example throughout, we consider the system of equations 
\bqq
\begin{cases}
\partial_t u&= d \partial_{xx} u+s\partial_x u+f(u)+\beta v, \\
\partial_t v&= -(\partial_{xx} +1)^2v+s\partial_x v+\mu v,
\end{cases} \quad t>0, \quad x\in\R,
\label{eq:main}
\eqq
where $f(u)\coloneqq\alpha u(1-u^2)$ for some $\alpha>0$. Parameters are set such that $d>0, \mu<0, \beta\neq 0$, and $s>0$. The $u$ component is governed by a Fisher-KPP type equation while the $v$ component is the Swift-Hohenberg equation typically used as a model for pattern formation, although we consider the case of $\mu<0$ where patterns are suppressed and the zero state is stable. The two equations are coupled through the $\beta v$ term in the first equation. In terms of traveling fronts, there exist monotone front solutions $(Q(x),0)$ connecting $(u_-,v_-)=(1,0)$ at $-\infty$ to $(u_+,v_+)=(0,0)$ at $+\infty$ for any wavespeed $s\geq s_*\coloneqq2\sqrt{df'(0)}=2\sqrt{d\alpha}$; see for example~\cite{AW78}. Similar systems have been studied in~\cite{gallay04, ghazaryan07}, focusing however on fronts between stable states $u=\pm 1$. 

Investigating stability of these traveling front solutions, we linearize \eqref{eq:main} near the traveling front solution $(Q_*(x),0)$. The essential spectrum of the resulting linear operator can be characterized in terms of the asymptotic systems as $x\to \pm \infty$, that is, linearizing at $(u_\pm,v_\pm)$, respectively.  The boundaries of the essential spectrum are given by spectra of these asymptotic systems, determined by solutions of the form $\me^{\mbi kx}$ with temporal growth rate $\lambda(k)$. Focusing on the asymptotic system near $x=+ \infty$ for \eqref{eq:main}, it is easy to see that the essential spectrum is stable for the $v$ component since $\mu<0$. On the other hand, the essential spectrum for the $u$ component is unstable since $f'(0)=\alpha>0$. In an exponentially weighted space, the essential spectrum instead records the temporal growth rate associated to modes of the form $\me^{(\mbi k+\eta)x}$ for some $\eta\in\R$. It is known that the essential spectrum for the $u$-component is real, $\{\lambda\leq \alpha-\frac{s^2}{4d}\}$,  for the specific choice $\eta=-s/(2d)$; see~\cite{sattinger}. It is possible however that the $v$-equation will possess unstable essential spectrum and therefore exhibit exponential growth for this specific choice of exponential weight $\eta=-s/(2d)$. This paper focuses precisely on this situation, in the case of the critical front with speed $s_*=2\sqrt{\alpha d}$ when the essential spectrum of the $u$-component is marginally stable, $\{\lambda\leq 0\}$. We shall refer to such instabilities, where the linearization at the solution $(Q_*(x),0)$ possesses essential spectrum in $\{\Re\lambda>0\}$ for any choice of exponential weights as \emph{remnant instabilities}, a term introduced in~\cite{sandstede00}.

In order to understand the possibility of stabilizing a system using exponential weights more systematically, we analyze the linearization at $(u_+,v_+)$ in more detail, referring in particular  to the \emph{absolute spectrum}, $\Sigma_\textnormal{abs}$, introduced in~\cite{sandstede00}. Roughly speaking, stability of the absolute spectrum implies that there exists a continuous choice of exponential weights $\eta(\lambda)$ in $\{\Re\lambda>0\}$  so that $\lambda$ does not belong to the essential spectrum in the space with exponential weight $\eta(\lambda)$. More algebraically, one computes all possible ``spatial eigenvalues'' $\nu(\lambda)$ corresponding to solutions $\rme^{\lambda t + \nu(\lambda) x}$ for a given $\lambda$ and determines if this collection can be consistently separated by real part along a curve with $\Re\lambda\to+\infty$. This consistent splitting breaks whenever $\lambda\in\Sigma_\textnormal{abs}$. In this way, the existence of unstable absolute spectrum implies what we shall refer to as a remnant instability~\cite{sandstede00}, an instability that cannot be stabilized using exponential weights. 

The absolute spectrum forms curves in the complex plane which terminate on double roots (sometimes referred to as branch points), where two values $\nu_{1/2}(\lambda)$ collide; see~\cite{rademacher07}. These double roots typically lead to \emph{resonance poles}, which are singularities of the resolvent kernel (also referred to as the pointwise Green's function),  and their location in the complex plane prescribes the pointwise exponential growth or decay rate of solutions to the homogeneous equation; see for example~\cite{bers84,brevdo96,briggs,holzerscheel14,huerre90,vansaarloos03}. Typically, but not always, those resonance poles form the rightmost part of the absolute spectrum. 

\begin{figure}
    \centering
     \subfigure[$\mu=-1/3$.]{\includegraphics[width=0.375\textwidth]{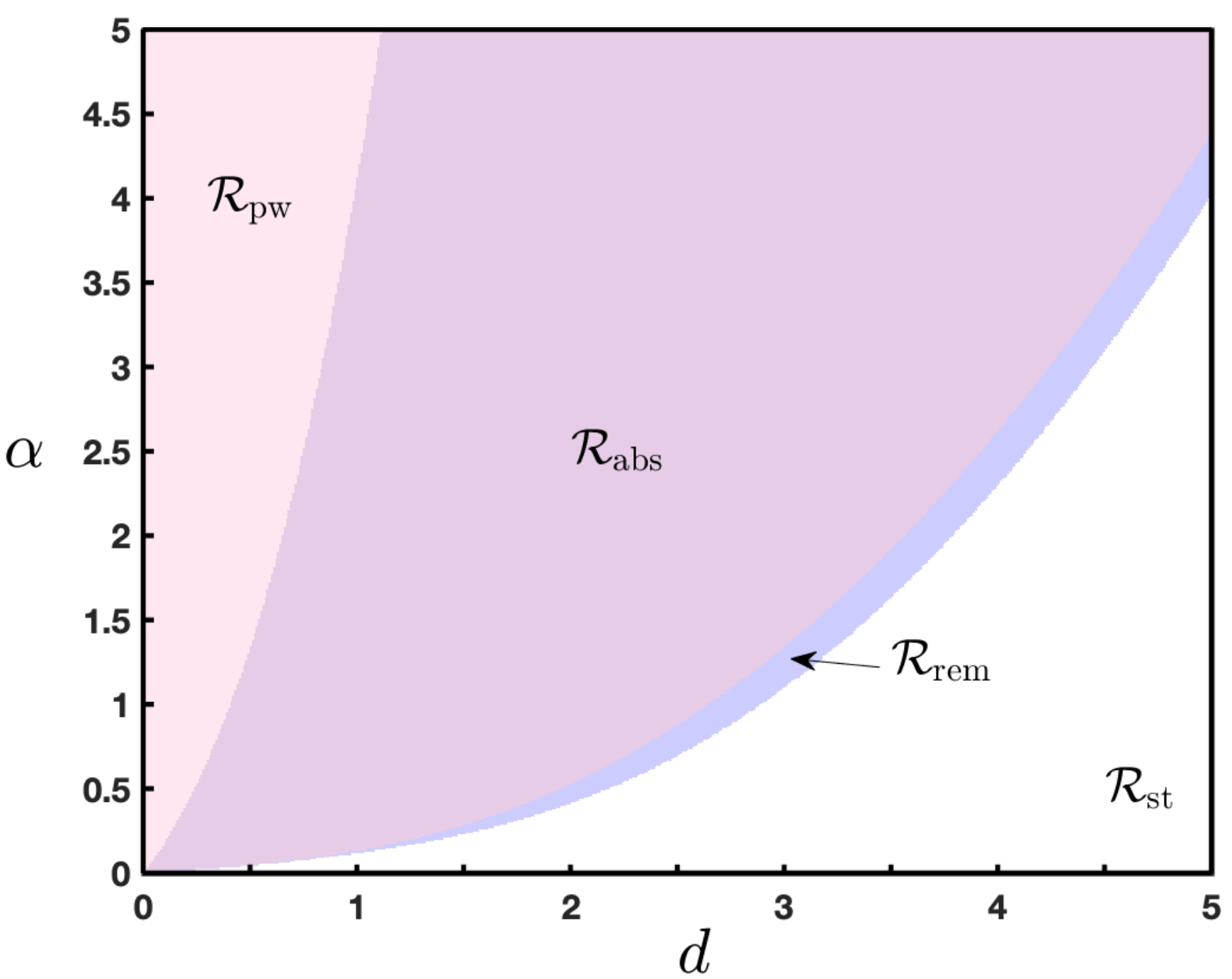}}
 \subfigure[$\alpha=3$.]{\includegraphics[width=0.375\textwidth]{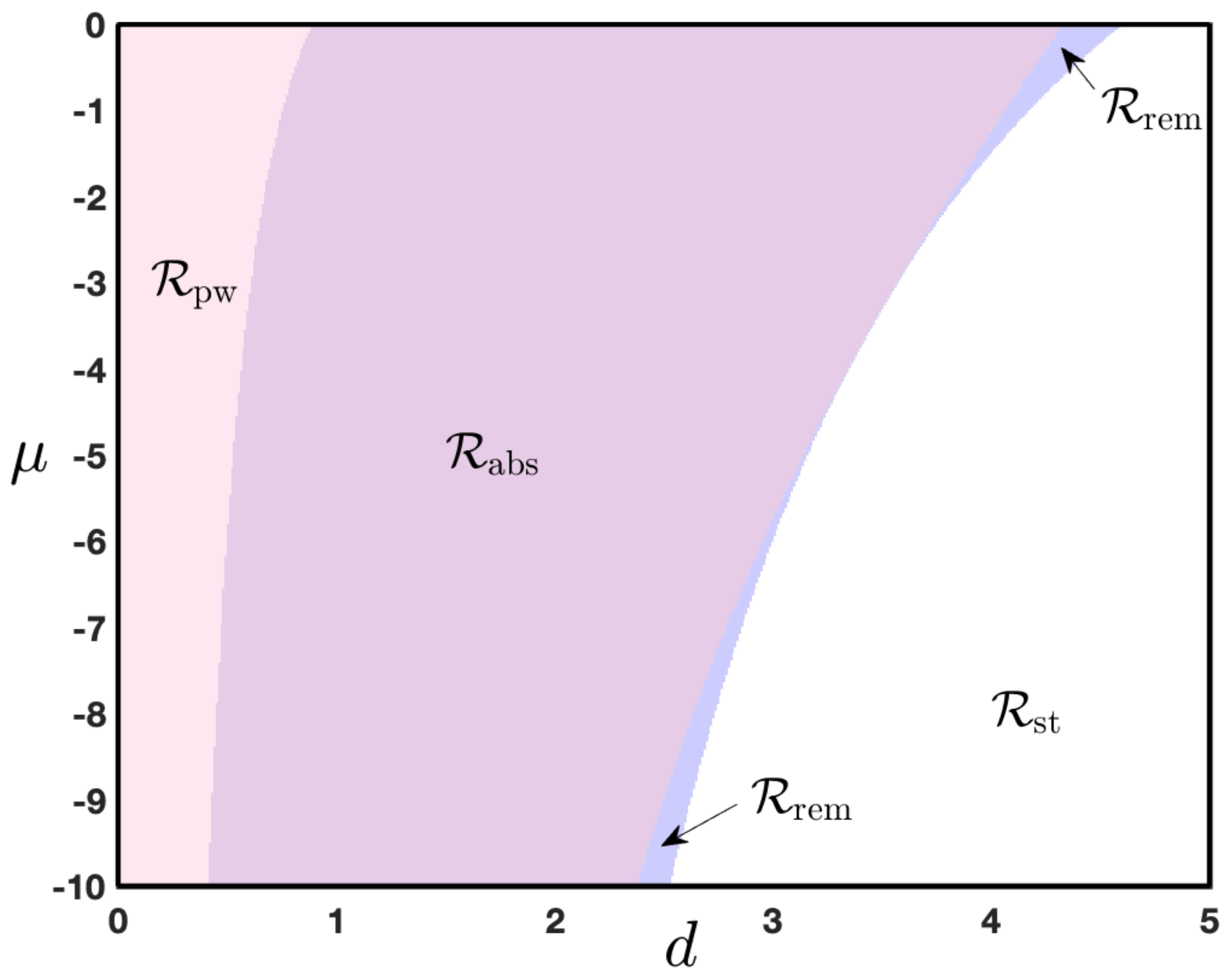}}
%  \subfigure[$d=3/4$.]{\includegraphics[width=0.32\textwidth]{Figures_manuscript/AllParmSpaceRegionsd075.pdf}}
   \caption{Numerically computed regions $\mathcal{R}_{\mathrm{st}}$, $\mathcal{R}_{\mathrm{rem}}$, $\mathcal{R}_{\mathrm{abs}}$, and $\mathcal{R}_{\mathrm{pw}}$ for $\mu=-1/3$ and $\alpha=3$ in \eqref{eq:main}, respectively; see  Table~\ref{table} for a region labels and  Figure~\ref{fig:representative} for space-time plots of solutions in different regions.}
    \label{fig:RegionsRj}
\end{figure}

In summary, unstable resonance poles imply unstable absolute spectrum, which in turn implies a remnant instability. In most cases in the literature, these three concepts of instability are in fact equivalent, that is, for given parameter values, a system exhibits all or none of those three instabilities. The main interest in our system \eqref{eq:main} is the fact that it provides a simple example where this equivalency fails, that is, there is a parameter region $\mathcal{R}_\mathrm{rem}$ with remnant instability but no unstable absolute spectrum and a parameter region $\mathcal{R}_\mathrm{abs}$ with unstable absolute spectrum but no unstable resonance poles. Table \ref{table} summarizes the characterization of these regions in comparison to Region $\mathcal{R}_{\mathrm{st}}$, where the essential spectrum of the unstable state can be (marginally) stabilized using an exponential weight, and Region $\mathcal{R}_{\mathrm{pw}}$ with unstable resonance poles. One expects stability in $\mathcal{R}_{\mathrm{st}}$ using a generalization of the argument for the scalar case, see~\cite{gallay94, faye19}, and instability in Region $\mathcal{R}_{\mathrm{pw}}$, where the unstable state possesses resonance poles in the right half of the complex plane and pointwise instability of the front is expected\footnote{The resonance poles here are indeed \emph{relevant} in the sense of \cite{holzer14, holzerscheel14}.}. Figure~\ref{fig:RegionsRj} depicts these regions in parameter space. 
We next outline the main contributions of this paper which provide stability results for fronts in the intriguing intermediate regions $\mathcal{R}_\mathrm{rem}$ and $\mathcal{R}_\mathrm{abs}$.

\begin{table}
\begin{center}
\begin{tabular}{c c c c l}
\hline
\hline
Region & \parbox[t]{1.1in}{\centering Remnantly Unstable} & \parbox[t]{1.3in}{\centering Absolute Spectrum Unstable} & \parbox[t]{1.2in}{\centering Resonance Poles Unstable}  & Front Stability \\[0.25in]

\hline
\hline
$\mathcal{R}_{\mathrm{st}}$ & no & no & no & stable, adaptation of~\cite{faye19} \\
 $\mathcal{R}_{\mathrm{rem}}$ & yes & no & no & stable, Theorem~\ref{thm:maininformal} \\
 $\mathcal{R}_{\mathrm{abs}}$ & yes & yes & no & stable, Theorem~\ref{thm:maininformal} \\
  $\mathcal{R}_{\mathrm{pw}}$ & yes & yes & yes & pointwise unstable~\cite{holzerscheel14}  \\
 \hline
 \hline
%  \label{table}
\end{tabular}
\end{center}
\vspace{-5mm}
\caption{Overview of different types of linear instability and results on stability of the critical front of~\eqref{eq:main}. Regions outlined in parameter space in Figure~\ref{fig:RegionsRj}, below.}
\label{table}
\end{table}

\paragraph{Nonlinear Asymptotic Stability.}

The main analytical result of this paper establishes nonlinear stability of the traveling front $(Q_*(x),0)$ propagating at the linear spreading speed $s_*=2\sqrt{d\alpha}$ for parameters in  $\mathcal{R}_{\mathrm{rem}}\cup\mathcal{R}_{\mathrm{abs}}$.   We delay a precise statement of this main result until Section~\ref{sec:Proof} and present an informal statement, here. 

We write perturbations of the front as $u(t,x)=Q_*(x)+P(t,x)$,  such that the pair $(P(t,x),v(t,x))$ solves
\bqq
\begin{cases}
\partial_t P&= d \partial_{xx} P+s_*\partial_x P+f(Q_*+P)-f(Q_*)+\beta v, \\
\partial_t v&= -(\partial_{xx} +1)^2v+s_*\partial_x v+\mu v,
\end{cases}\quad \quad t>0, \quad x\in\R.
\label{eq:mainPv}
\eqq
For $\delta>0$ sufficiently small, let define a positive, bounded, and smooth weight function with
\begin{equation}
\omega(x)=\left\{ 
\begin{array}{ll}
\me^{-\frac{s_*}{2d}x}, & x \geq 1,\\
1,& x=0,\\
\me^{\delta x}, & x \leq -1,
\end{array}
\right.\label{eqomega}
\end{equation}
where we recall that this weight is the one required for the stability analysis of the scalar Fisher-KPP equation; see once again \cite{faye19,gallay94}.

\begin{thm}[Nonlinear stability -- informal statement]\label{thm:maininformal}
For parameters in a subregion $\Pi\subset \mathcal{R}_{\mathrm{rem}}\cup\mathcal{R}_{\mathrm{abs}}$, solutions to~\eqref{eq:mainPv} with sufficiently small initial conditions $(P_0,v_0)\in \X\times\Y$, with weighted function spaces defined in Theorem \ref{thmMain}, exist for all $t>0$ and converge to zero in the following sense. There is a pointwise decomposition
\[
P(t,x)=E_u(t,x)+E_v(t,x), \quad x\in\R,
\]
such that 
\begin{align*}
|E_v(t,x)| &\leq C \me^{-\theta t} \me^{-\gamma |x|} \|v_0\|_\X, \quad t\geq 1,\ x\in\R,\\
|E_u(t,x)| & \leq C t^{-3/2} (1+|x|)\omega(x) \left(\|v_0\|_\X+\|P_0\|_\Y\right),\\
|v(t,x)|&\leq C\me^{-\theta t}\me^{-\delta|x|} \|v_0\|_\X,
\quad t \geq 1, \ x\in\R,
\end{align*}
for some constants $C>0$, $\theta>0$, $0<\gamma<\frac{s_*}{2d}$,

On the other hand, there exist $\varepsilon>0$ and initial conditions $v_0\in \X$ arbitrarily small such that for some $t>0$, 
\[
  \sup_{x\in\R}\frac{|P(t,x)|}{\omega(x)} \geq \varepsilon,
\]
that is, the front is not asymptotically stable in the fixed weight. 
\end{thm}
We refer the reader to Theorem~\ref{thmMain} for a precise statement of the stability result. The region $\Pi$ is described in Definition~\ref{defiPi} and excludes potential nonlinear 3:1-resonances in the sense of \cite{faye17}. The proof  relies on  pointwise semigroup methods, developed in~\cite{zumbrun98} for stability problems in the presence of marginally stable essential spectrum; see  \cite{beck14,howard02,zumbrun11} among others.  The method was recently applied in~\cite{faye19} to re-derive the nonlinear stability of the critical Fisher-KPP front~\cite{gallay94,kirchgassner92, bricmont92, eckmann94, faye20}. The pointwise representation allows us to separate the weak algebraic decay in exponentially weighted spaces intrinsic to the $u$-equation and the exponential unweighted decay in $u$ induced by coupling to the $v$-equation.

\begin{figure}[!t]
\centering
\subfigure[$(d,\alpha,\mu) \in \mathcal{R}_{\mathrm{st}}\ \text{or} \ \mathcal{R}_{\mathrm{rem}}$.]{\includegraphics[width=0.32\textwidth]{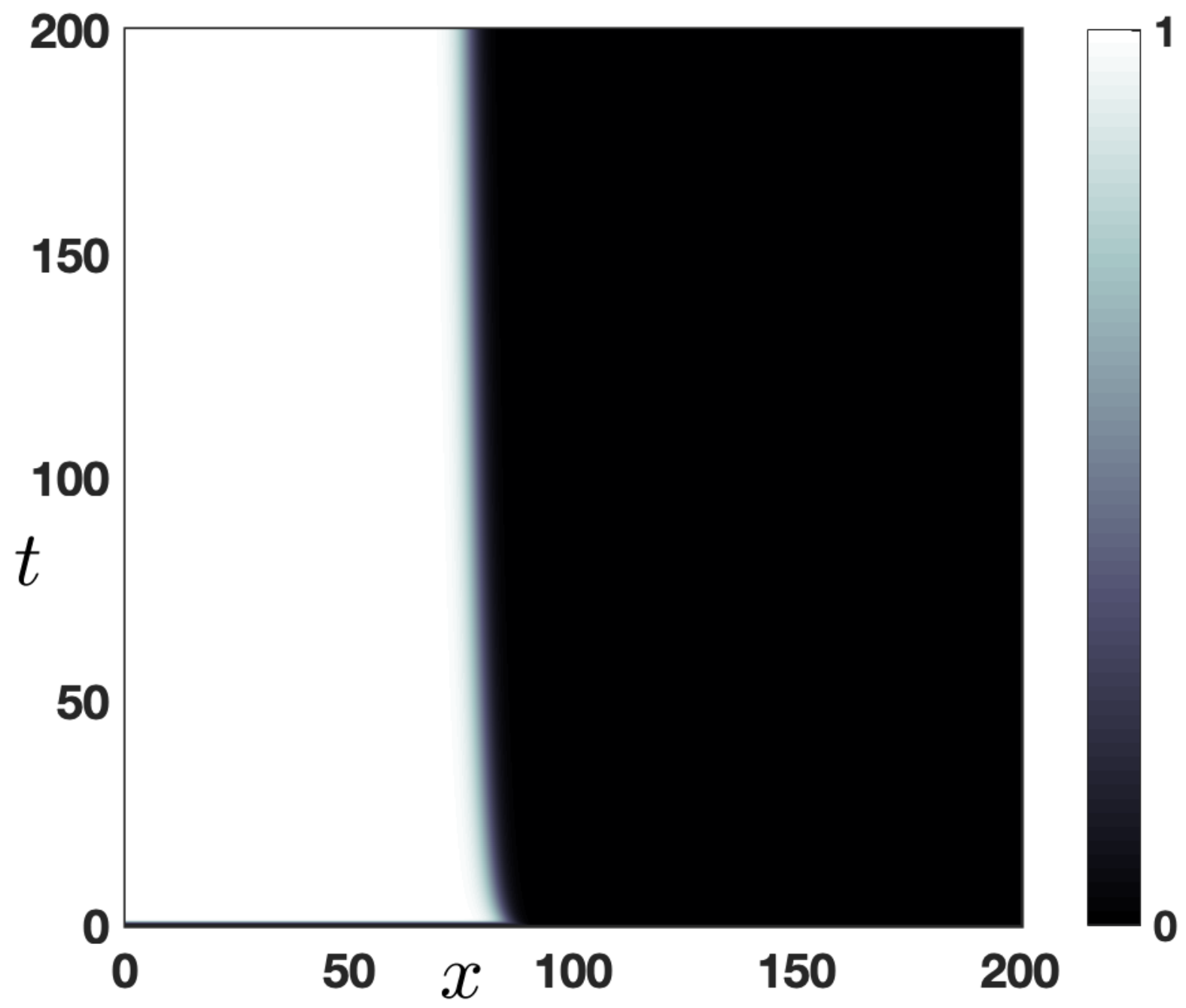}}
\subfigure[$(d,\alpha,\mu)\in \mathcal{R}_{\mathrm{abs}}$.]{\includegraphics[width=0.32\textwidth]{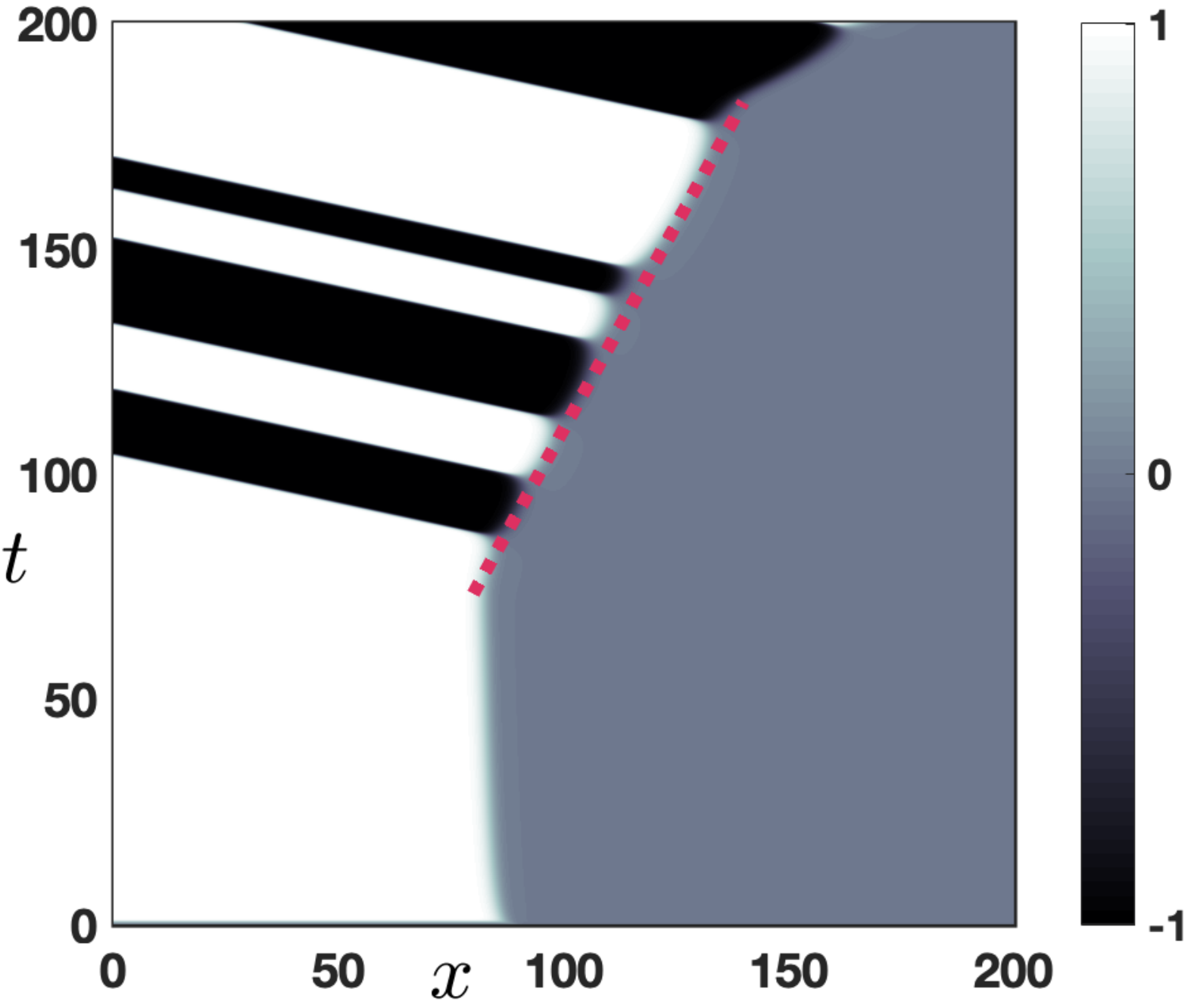}}
\subfigure[$(d,\alpha,\mu)\in \mathcal{R}_{\mathrm{pw}}$.]{\includegraphics[width=0.32\textwidth]{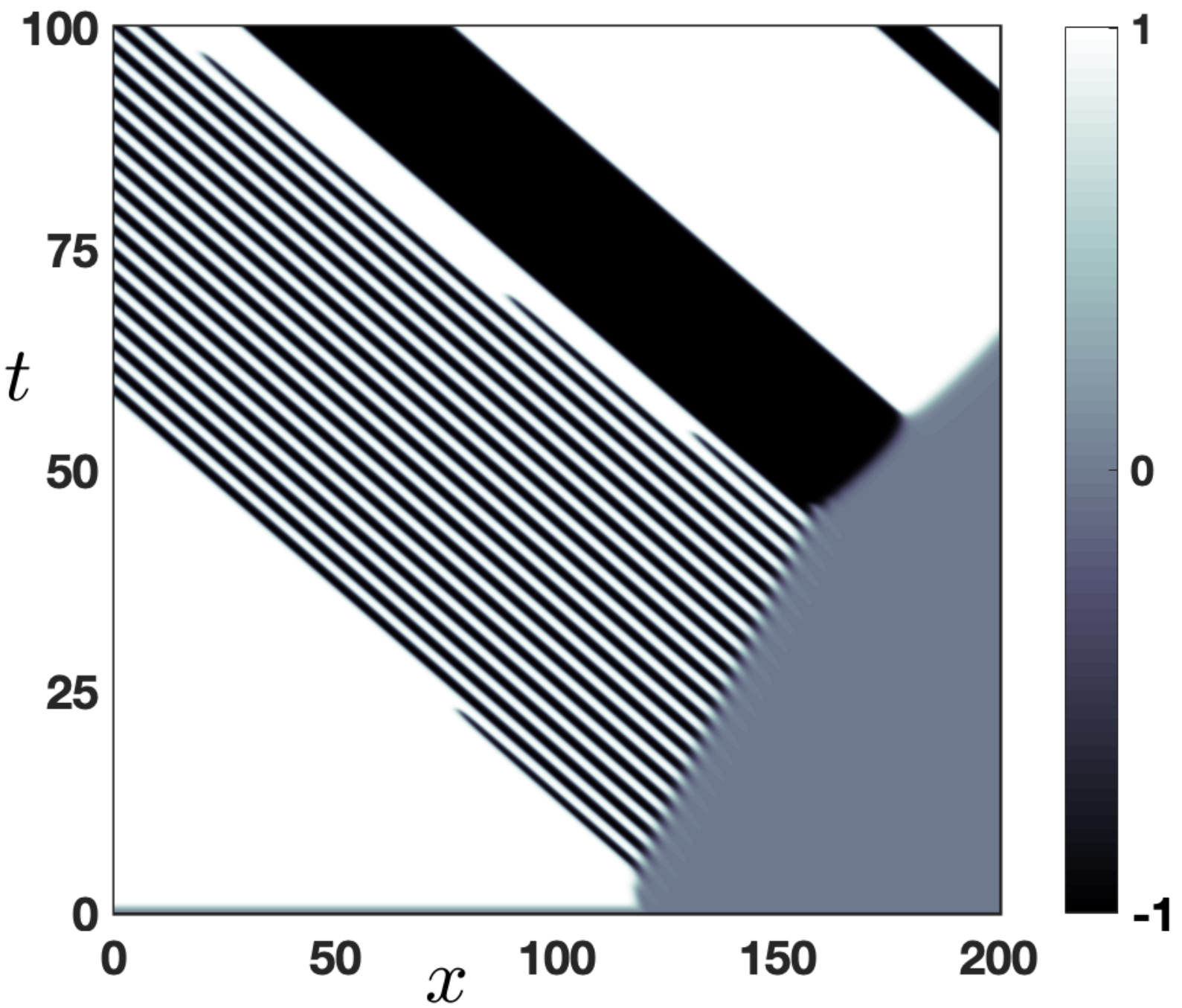}}
\caption{Space-time plots for~\eqref{eq:main}  in the comoving frame $s=s_*$, $\beta=1$, with parameters in $\mathcal{R}_{\mathrm{rem}}$~(a), $\mathcal{R}_{\mathrm{abs}}$~(b) , and $\mathcal{R}_{\mathrm{pw}}$~(c), respectively. (a) Stability of critical front, which remains stationary despite remnant instability (plots in region $\mathcal{R}_{\mathrm{st}}$ would be identical); see Theorem \ref{thm:maininformal} and \cite{faye19}, resp.. (b) Stability despite unstable absolute spectrum for long times, Theorem \ref{thm:maininformal}, and instability and acceleration to the absolute spreading speed for very long times due to round-off, Conjecture \ref{prop:sabs}. (c) Instability against pattern-forming front due to unstable resonance pole, then absolute spreading speed from round-off error. }
\label{fig:representative}
\end{figure}

\paragraph{Wavespeed selection and the role of the absolute spectrum.}

Numerical simulations illustrating stability are shown in Figure~\ref{fig:representative}. 
For parameters in $\mathcal{R}_{\mathrm{abs}}$, the front is stable, stationary in the comoving frame for a long time period as predicted but overtaken by a faster invasion mode after very long times. We will demonstrate that the speed of this faster invasion mode is approximately the speed $s_{\textnormal{abs}}$ at which the absolute spectrum becomes (marginally) stable. We will present numerical evidence that this faster invasion speed is induced by numerical round-off errors. Corroborating this numerical evidence is the following result that we formulate as a conjecture here.

\begin{conj}\label{prop:sabs} Consider the following modification of (\ref{eq:main}) where the coupling term $\beta v$ is replaced with $\beta \sigma(x) v$ for some inhomogeneity $\sigma(x)$, 
\bqq
\begin{cases}
\partial_t u&= d \partial_{xx} u+s\partial_x u+f(u)+\beta \sigma(x) v, \\
\partial_t v&= -(\partial_{xx} +1)^2v+s\partial_x v+\mu v,
\end{cases} \quad t>0, \quad x\in\R.
\label{eq:mainmod}
\eqq 
Suppose that the Fourier transform of $\sigma(x)$ has full support. Then for parameters in $\mathcal{R}_{\mathrm{abs}}$, the traveling wave $(Q_*(x),0)$ is pointwise unstable due to resonance poles accumulating on the unstable absolute spectrum $\Sigma_{abs}$.  
\end{conj}
We motivate this conjecture using a pointwise analysis and present numerical evidence in Section~\ref{sec:Num}.

Taking a step back, an underlying motivation of ours is to understand wavespeed selection mechanisms and determine invasion speeds based upon properties of the linearization at the unstable state whenever possible; see~\cite{vansaarloos03} for a review of wavespeed selection principles with many applications. In other words, one seeks to identify a ``linear'' spreading speed based on linear information that gives the nonlinear asymptotic spreading speed associated to compactly supported perturbations of the unstable state.
In the most common case, the linear spreading speed is induced by a simple branch point located  on the imaginary axis. Fronts invading at this linear spreading speed are called {\em pulled}, distinguishing them from a case of \emph{pushed fronts} where the spreading speed is faster, determined and driven by nonlinearity; see agin~\cite{vansaarloos03}. Going beyond these classical predictions, we showed in~\cite{faye17} that generic nonlinear mixing of linear modes can lead to faster invasions whose speeds can be determined based solely on information from the linearized unstable state. Conjecture~\ref{prop:sabs}, motivated in part by~\cite{goh11} suggests that, in the presence of generic multiplicative fluctuations, the absolute spreading speed is most relevant for invasion problems.

\paragraph{Estimates on the boundaries of instability regions.} Complementing the previous nonlinear analysis, we now demonstrate that our simple model actually exhibits all phenomena discussed above, that is, Theorem~\ref{thm:maininformal} holds in the regions $\mathcal{R}_{\mathrm{rem}}$ and $\mathcal{R}_{\mathrm{abs}}$ which are both nonempty.

\begin{thm}\label{thm:paraminformal}
For $\alpha,d>0,\ s=s_*$ fixed, there exist 
\begin{equation}\label{e:muremabs}
\mu^{\textnormal{rem}}(\alpha,d)\coloneqq 2\alpha -\frac{4\alpha}{d} - \frac{8\alpha^2}{d^2}\leq \mu^\textnormal{abs}_0(\alpha,d)\coloneqq \frac{d^2}{4} -\frac{4\alpha}{d} -\frac{4\alpha^2}{d^2}, 
\end{equation}
such that the origin has
\begin{itemize}[itemsep=-0pt,topsep=-0pt]
\item[(i)]  unstable absolute spectrum if $\mu^\textnormal{abs}_0<\mu<0$, that is,  $\{(\mu,\alpha,d)|\,\mu\in (\mu^\textnormal{abs}_0,0)\}\subset \mathcal{R}_\mathrm{abs} \cup \mathcal{R}_\mathrm{pw}$;
\item[(ii)] a  remnant instability if $\mu^\textnormal{rem}<\mu<0$,  that is, $\{(\mu,\alpha,d)|\,\mu\in (\mu^\textnormal{rem},0)\}\subset \mathcal{R}_\mathrm{abs}\cup \mathcal{R}_\mathrm{rem} \cup \mathcal{R}_\mathrm{pw}$.
\end{itemize}
Moreover, 
\begin{itemize}[itemsep=-0pt,topsep=-0pt]
 \item [(iii)] $\mathcal{R}_\mathrm{abs}\cup \mathcal{R}_\mathrm{rem}\neq \emptyset$ when $\mu^{\textnormal{rem}}<0$;
 \item[(iv)]$\mathcal{R}_\mathrm{abs}\neq \emptyset$ when  $\mu^\textnormal{abs}_0<0$ and $\alpha>d+\frac{d^2}{2}$;
 \item[(v)] $\mathcal{R}_\mathrm{rem}\neq \emptyset$ when $\mu^\mathrm{rem}<0$ and $d^2\neq 4\alpha$.
\end{itemize}
 \end{thm}
Proofs can be found in Sections~\ref{sec:EssSpec} and~\ref{sec:AbsSpec}.

\paragraph{Outline.} The remainder of the paper is organized as follows. In Section~\ref{sec:EssSpec}, we discuss  essential and weighted essential spectra of the linear asymptotic operators and characterize the boundary of the remnant instability region. In Section~\ref{sec:Proof}, we present our main stability result and its proof. In Section~\ref{sec:AbsSpec}, we present explicit formulas for the absolute spectra of the Fisher-KPP as well as Swift-Hohenberg components in isolation, characterize the absolute spectrum associated to~\eqref{eq:main} near the unstable state, and characterize (non empty) regions of parameters where our main nonlinear stability result holds. Finally, we provide some numerical simulations in Section~\ref{sec:Num} that show the emergence of pointwise growth modes. 
 
\paragraph{Acknowledgments.} The research of M.\ H.  was partially supported by the National Science Foundation through DMS-$2007759$.  A.\ S. was supported by the National Science Foundation through NSF DMS-1907391  L.\ S. acknowledges support by the Deutsche Forschungsgemeinschaft (German Research Foundation), Projektnummer 281474342/GRK2224/1.

\section{Exponential weights and remnant instabilities}\label{sec:EssSpec}

We review essential and weighted essential spectra and identify parameters of remnant instability in Proposition \ref{prop:boundaryremnant}, below, where exponential  weights cannot stabilize the essential spectrum of (\ref{eq:main}) linearized near the traveling front. In Section~\ref{sec:DispersionRelation} we introduce the dispersion relation and its roots.  In Section~\ref{sec:ExponentialWeights} we discuss weighted essential spectra and derive Proposition~\ref{prop:boundaryremnant}. In Section~\ref{sec:ptw}, we discuss pointwise and other stability concepts. 

While most of this section is relevant for Theorem \ref{thm:paraminformal}, only, we also introduce notation that will be required in subsequent sections including the weighted essential spectrum $\Sigma_{\mathrm{ess}}^\eta(\mathcal{L})$, the dispersion relation $D(\lambda,\nu)$ and its roots ($\nu_j(\lambda)$ for the $v$ component, $\nu_\pm^0(\lambda)$ and $\nu_\pm^1(\lambda)$ for the $u$ component near zero and one, respectively), and finally pinched double roots that lead to singularities of the pointwise Green's function ($\lambda_v^{\mathrm{bp}}$ for the $v$ component and $\lambda^{\mathrm{rp}}$ for resonance poles coming from the coupling of the two equations).  We also provide some background on the signifigance of each of these terms, although we refer the reader to \cite{holzerscheel14} for a more in-depth treatment.

\subsection{The dispersion relation and essential spectra}\label{sec:DispersionRelation}
We study the asymptotic linearized operators in a comoving frame,
\[\calL^+\coloneqq\begin{pmatrix}
\calL^+_u & \beta\\ 0&\calL^+_v
\end{pmatrix},     \quad \qquad  \quad\calL^-\coloneqq\begin{pmatrix}
\calL^-_u & \beta\\ 0&\calL^-_v
\end{pmatrix},\]
where plus and minus indicate the asymptotic rest states $(0,0)$ and $(1,0)$ at $+\infty$ and $-\infty$, respectively, and where 
\[\calL^+_u\coloneqq d\partial_{xx}+s\partial_x +\alpha, \quad \calL^-_u\coloneqq d\partial_{xx} + s\partial_x -2\alpha, \quad \calL^\pm_v\coloneqq -(\partial_{xx}+1)^2+s\partial_x+\mu,\]
are the componentwise linearizations.

The dispersion relation $D(\lambda,\nu)$ relates spatial modes $\rme^{\nu x}$ to their exponential temporal growth rates $\rme^{\lambda t}$. Due to the skew-product nature, the dispersion relation associated to the linearization near the zero state is given by the following product
\begin{equation}\label{eq:DispersionRelation}
D^0(\lambda,\nu)=D^0_u(\lambda,\nu)\cdot D_v(\lambda,\nu)=(d\nu^2+s\nu+\alpha-\lambda)\cdot(-\nu^4-2\nu^2+s\nu-1+\mu-\lambda ).
\end{equation}
The dispersion relation has six spatial roots. The two roots of $D^0_u(\lambda,\nu)$ have analytic expressions which we denote as 
\begin{equation}\label{eq:SpatialRootsU}
\nu_\pm^0(\lambda)=-\frac{s}{2d}\pm\frac{1}{2d}\sqrt{s^2-4d\alpha+4d\lambda}
\end{equation}
while $D_v(\lambda,\nu)$ has four roots which we denote by $\nu_j(\lambda)$ with the specifications that  
\bqq \lim_{\mathrm{Re}(\lambda)\to +\infty} \mathrm{Re}(\nu_{1,2}(\lambda))<0, \quad \lim_{\mathrm{Re}(\lambda)\to +\infty} \mathrm{Re}(\nu_{3,4}(\lambda))>0.\label{eq:pinching} \eqq
When discussing the absolute spectrum it will be more convenient to label roots according to their real part. So, we will also  denote the four spatial eigenvalues of $\calL^\pm_v$ by $\rho_j(\lambda)$ defined so that for all  $\lambda\in\C$,
\bqq\Re(\rho_1(\lambda))\le\Re(\rho_2(\lambda))\le\Re(\rho_3(\lambda))\le\Re(\rho_4(\lambda)).
\label{eq:rootsvorderingabs}\eqq
By Fourier transform, the spectrum of $\calL^+$ as a linear operator on $L^2(\mathbb{R})\times L^2(\mathbb{R})$ consists of the set of spectral values $\lambda\in\mathbb{C}$ for which there exists a solution $D^0(\lambda,\mbi k)=0$ for some $k\in\mathbb{R}$, that is, $\Sess(\calL^+)=\Sess(\calL^+_u)\cup \Sess(\calL^+_v)$, with 
\[  \Sess(\calL^+_u)= \left\{ -dk^2+s\mbi k+\alpha, \ k\in\mathbb{R}\right\} , \quad \Sess(\calL^+_v) = \left\{ -k^4+2k^2+s\mbi k-1+\mu, \ k\in\mathbb{R}\right\} \]
Since throughout $\mu<0$, $\Re(\Sess(\calL^+_v))<0$. However, since $\alpha>0$, $\Sess(\calL^+_u)$ extends into the right half of the complex plane reflecting instability of the invaded state. 

\subsection{Exponential weights}\label{sec:ExponentialWeights}
\begin{figure}[t!]
\centering
\subfigure[$\Sess^\eta(\calL^+_u)$ as $\eta$ varies.]{\includegraphics[width=0.325\textwidth]{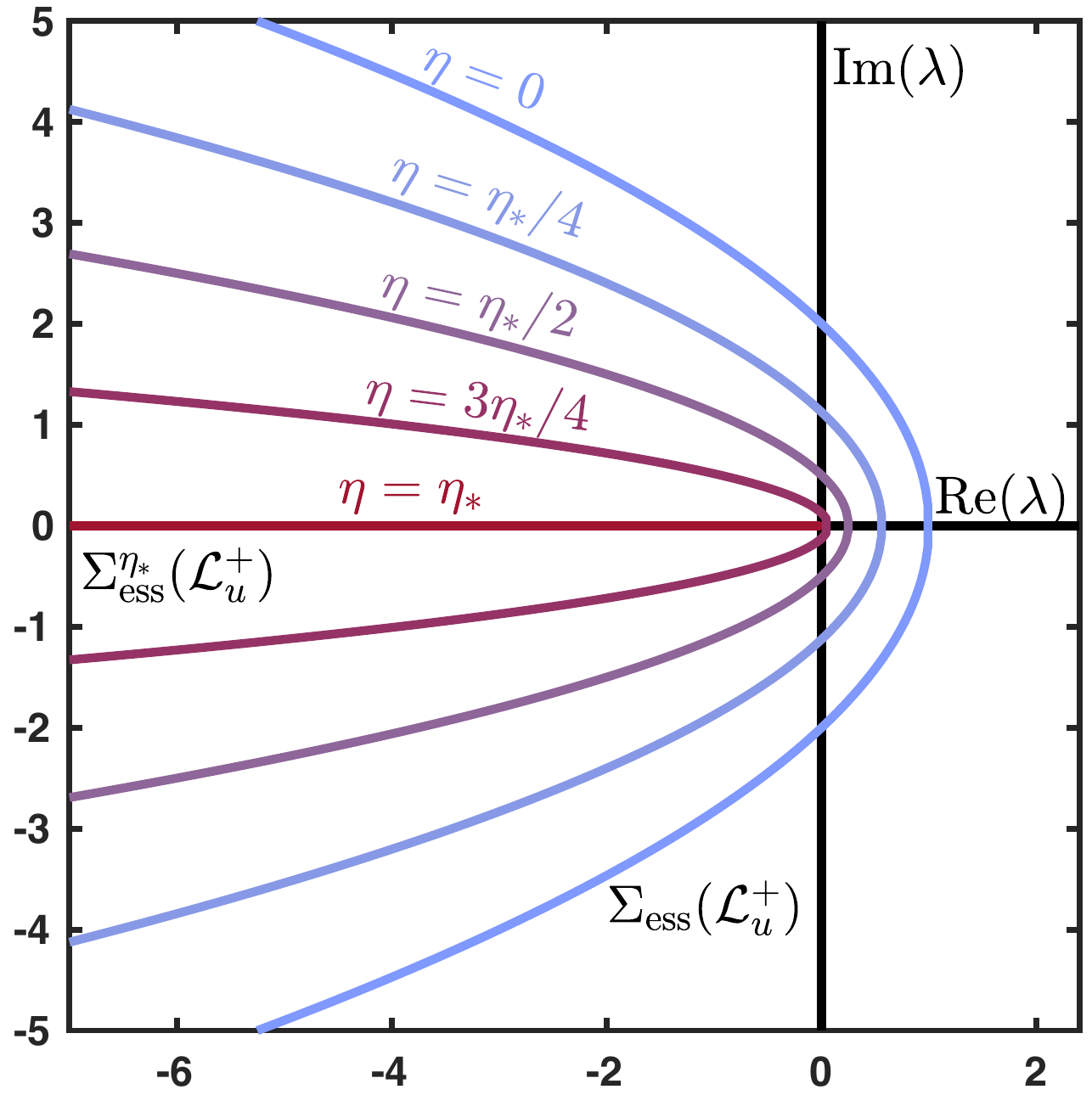}}\hspace*{.7in}
\subfigure[$\Sess^\eta(\calL^+_v)$ as $\eta$ varies.]{\includegraphics[width=0.325\textwidth]{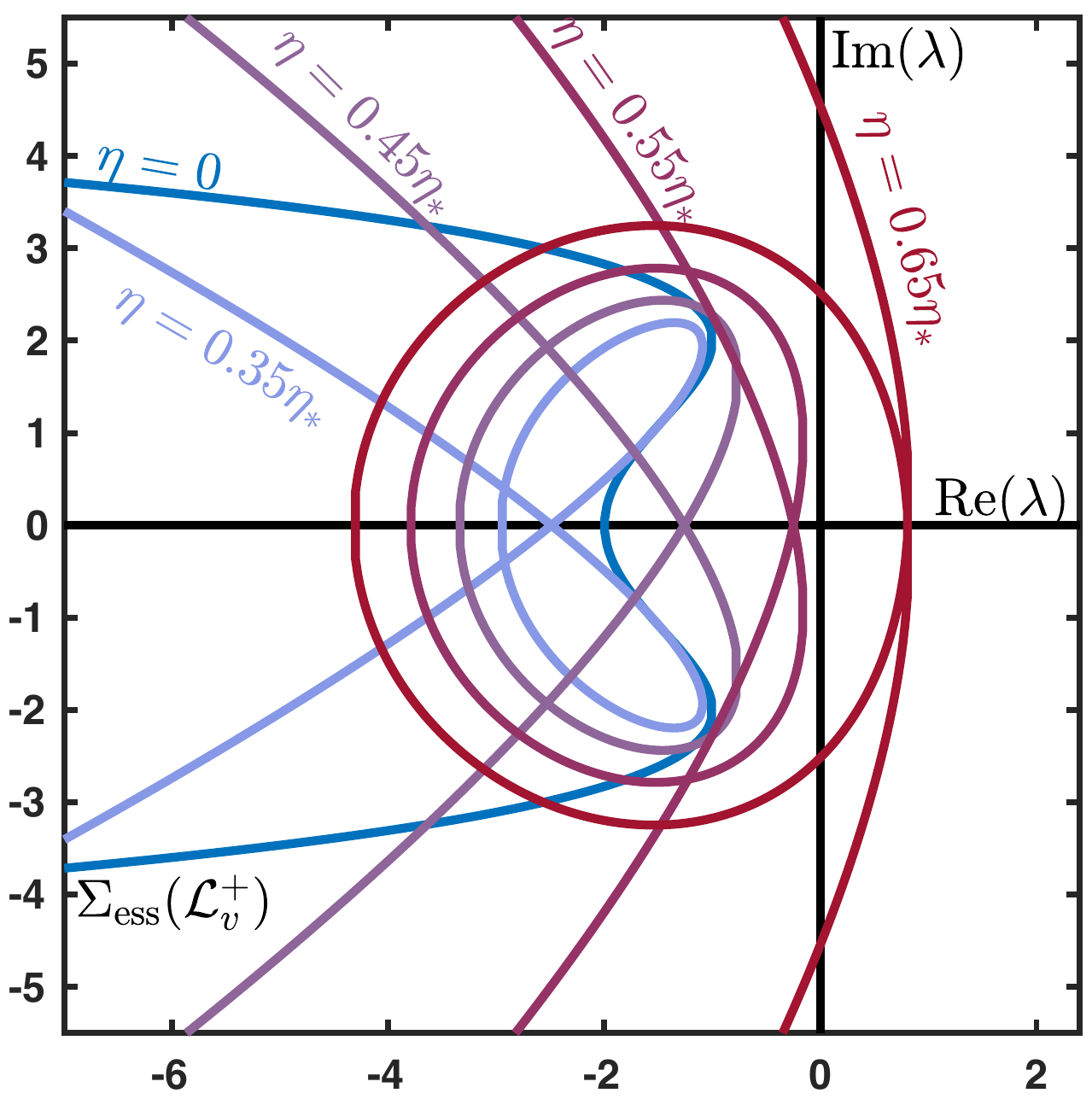}}
\caption{Weighted essential spectra from the $u$ (left) and $v$ (right) component in the complex plane for various choices of weights $\eta$ and $s=s_*$ and other parameters are fixed such that~\eqref{eq:BothWeightedSpectraStableCondition} fails. Note that $\Sess(\calL^+_u)$ is  marginally stabilized using the optimal weight $\eta_*$ while destabilizing $\Sess^\eta(\calL^+_v)$ for $\eta<\eta_\mathrm{c}\sim 0.5695 \eta_*$.}
\label{fig:shiftedES}
\end{figure}

In order to recover stability of the invasion process, one utilizes exponential weights that stabilize the essential spectrum. Therefore, consider the weighted Sobolev space $L^2_\eta(\R)\subset L^2_\mathrm{loc}$ with norm
\[\|u(\cdot)\|_{L^2_\eta}=\|u(\cdot)\rme^{\eta \cdot}\|_{L^2}, \quad \eta \in\R.\]
Clearly, by measuring perturbations in $L^2_\eta$ we penalize mass at $+\infty$  for $\eta>0$. On the other hand, measuring solutions, the norm neglects growth at $-\infty$. Stability in a weighted norm of an otherwise unstable system therefore reflects unidirectional transport. Again using Fourier transform, we obtain the weighted essential spectra through parameterized curves
\begin{align*}
 \Sess^\eta(\calL^+_u)= &\left\{ \sigma_u(k;\eta), \ k\in\mathbb{R}\right\} , \qquad 
 \Sess^\eta(\calL^+_v) = \left\{ \sigma_v(k;\eta), \ k\in\mathbb{R}\right\},\\
 \sigma_u(k;\eta)\coloneqq& -dk^2 + \mbi(s+2d\eta)k+d\eta^2+s\eta+\alpha\\
 \sigma_v(k;\eta)\coloneqq &-k^4+4\mbi \eta k^3 +(2+6\eta^2)k^2+\mbi(s-4\eta^3-4\eta)k-(1+\eta^2)^2+s \eta+\mu.
\end{align*}
Note that rightmost  point of $\Sess^\eta(\calL^+_u)$ occurs  for zero wavenumber  $k=0$. The maximal real part of $\Sess^\eta(\calL^+_u)$ is minimized for  $\eta=-s/(2d)$, when
\begin{equation}\label{eq:OptimalWeightedSpecU}
\Sess^\eta(\calL^+_u)\big\vert_{\eta=-s/(2d)}=\left(-\infty, -\frac{s^2}{4d}+\alpha\right].
\end{equation}
In case $s=s_*=2\sqrt{d\alpha}$ we also observe that this optimal weight $\eta_*\coloneqq -\sqrt{\alpha/d}$ implies marginal stability, $\mathrm{max}\,\Sess^{\eta_*}(\calL^+_u)=0$; see Figure~\ref{fig:shiftedES}.
For the $v$-component, a short calculation shows that the real part of $\Sess^\eta(\calL^+_v)$ is maximized at $k=\pm\sqrt{3\eta^2+1}$ with 
\[\Re(\sigma_v(\pm\sqrt{3\eta^2+1};\eta))= 8\eta^4+4\eta^2+s \eta+\mu.\]
and note that $\Re(\sigma_v(0;\eta))<\Re(\sigma_v(\pm\sqrt{3\eta^2+1};\eta)),\forall \, \eta\in\R$. 

For $s=s_*$, we conclude that there exists an exponential weight that stabilizes both $u$ and $v$ linearizations simultaneously if and only if $\Re(\sigma_v(\pm\sqrt{3\eta_*^2+1};\eta_*))\leq 0$, a region given by 
\begin{equation}\label{eq:BothWeightedSpectraStableCondition}
\mu< \min(\mu^\textnormal{rem},0);
\end{equation}
see Figure~\ref{fig:shiftedES} (b) for an illustration. These considerations are summarized in the following proposition.

\begin{prop}[Boundary of remnant instability]\label{prop:boundaryremnant}
With $\mu^\textnormal{rem}(\alpha,d)\coloneqq 2\alpha-4\alpha/d-8\alpha^2/d^2$, we have that, within $\{\mu<0,\alpha>0,d>0\}$, 
\begin{align*}
 \mathcal{R}_{\mathrm{st}}&=\left\{ (d,\alpha,\mu) ~|~ \mu < \min(\mu^\textnormal{rem},0)  \text{ and } \alpha,d>0 \right\},\\
\mathcal{R}_{\mathrm{pw}}\cup \mathcal{R}_{\mathrm{rem}}\cup \mathcal{R}_{\mathrm{abs}} &= \left\{ (d,\alpha,\mu) ~|~ \mu^\textnormal{rem}< \mu<0 \text{ and } \alpha,d>0 \right\}.
\end{align*}
In particular, there are instabilities with $\mu<0$, $\mathcal{R}_{\mathrm{pw}}\cup \mathcal{R}_{\mathrm{rem}}\cup \mathcal{R}_{\mathrm{abs}}\neq \emptyset$, whenever $\mu^\textnormal{rem}(\alpha,d)<0$, that is, when $\alpha> d^2/4-d/2$.
\end{prop}

\subsection{Dispersion relations and pointwise instabilities}\label{sec:ptw}
We briefly discuss the role of branch points and resonance poles in pointwise instabilities. 

A possible reason for remnant instability is pointwise instability, when localized initial conditions grow in any fixed finite window of observation. Since the effect of the exponential weight in a bounded window is irrelevant, exponential weights cannot stabilize the system.  Pointwise instabilities are typically caused by unstable branch points of the dispersion relation which lead to  resonance poles of the corresponding pointwise Green's function.  Instability of fronts propagating into states exhibiting a pointwise instability in the frame moving with the speed of the front is expected (see however~\cite{holzer14, holzerscheel14, holzerscheelLV} for subtleties pertaining to the definition of branch point) and Theorem \ref{thm:maininformal} therefore is concerned with remnant but not pointwise instabilities.

Recall the dispersion relation for the  linearization at the origin, $D^0(\lambda,\nu)$ from \eqref{eq:DispersionRelation}. We solve $D^0(\lambda,\nu)=0$ for spatial modes $\nu$ as functions of temporal growth $\lambda$. For $\mathrm{Re}(\lambda)\gg1$ one finds 3 roots with $\Re\nu<0$ and 3 roots with $\Re\nu>0$. Decreasing $\mathrm{Re}(\lambda)$, one of the $\nu$'s with positive real part will eventually cross the imaginary axis at a location $\mbi k$, precisely when $\lambda\in \Sess(\mathcal{L}^+)$ where $D^0(\lambda,\mbi k)=0$. 
Still, the roots $\nu(\lambda)$ can be continued  analytically in $\lambda$ as long as they are simple. Deforming Fourier transform integrals in the complex plane, one readily notices that locations $\lambda$ of double roots, where two roots $\nu_{I/II}(\lambda)$ collide are the only candidates for singularities of the Green's function to $\mathcal{L}^+-\lambda$.  One can further verify that only \emph{pinched double roots}, where $\Re\nu_{I}(\lambda)\to +\infty$ and  $\Re\nu_{II}(\lambda)\to -\infty$, can cause such singularities; see~\cite{bers84, brevdo96, briggs, holzerscheel14, huerre90, vansaarloos03}. Thinking in terms of stabilizing the system with exponential weights, one sees immediately that pinched double roots are lower bounds on essential spectra in any exponentially weighted space, since weights $\eta$ need to separate $\nu_{I/II}(\lambda)$ by real part for all $\lambda$: pointwise instabilities imply remnant instabilities.

Thinking in terms of pointwise stability, we can construct the heat kernel to the linearized PDE via inverse Laplace transform of the Green's function to $\mathcal{L}^+-\lambda$. Deformation of the path integrals in the inverse Laplace transform then exhibit that singularities of the pointwise Green's function encode exponential growth or decay in time of the heat kernel. 
In the notation of \eqref{eq:DispersionRelation}--\eqref{eq:pinching}, we find  pinched double roots when $\nu_+^0(\lambda)=\nu_-^0(\lambda)$, and when $\nu_j(\lambda)=\nu_k(\lambda)$ for $j=1$ or $2$ and $k=3$ or $4$ describing pointwise growth rates of $u$- and $v$-components without coupling.  Additional double roots that lead to singularities occur due to the coupling of the two components where $\nu_0^+(\lambda)=\nu_j(\lambda)$ for $j=1$ or $2$. Due to the skew-product nature of (\ref{eq:main}) this lead to poles, rather than branch points, of the pointwise Green's function and we will refer to these singularities as resonance poles.

We finally conclude this section by introducing the notations that will be used in Section~\ref{sec:Proof} (see also Remark~\ref{rmk:singularities}). The two complex conjugates pinched double roots which are branch points of the $v$ dispersion relation $D_v(\lambda,\nu)=0$ will be denoted by $(\lambda^\textnormal{bp}_v,\nu^\textnormal{bp}_v)$ with complex conjugate. The two aforementioned pinched resonance poles will be denoted $(\lambda^\textnormal{rp},\nu^\textnormal{rp})$ with complex conjugate; see  Section~\ref{sec:AbsSpec} for explicit formulas and more details. 

\section{Proof of nonlinear stability for Fisher-KPP invasion front}\label{sec:Proof}
In this section, we present and prove the main result  on nonlinear stability.   Section~\ref{sec:argoutline} presents set-up and outline,  Section~\ref{sec:pointwise} develops estimates on the pointwise Green's function, and Section~\ref{sec:temporal} uses those to obtain estimates for the temporal Green's function via the inverse Laplace Transform.  Finally, Section \ref{sec:nonlinear} contains a precise statement and proof of the nonlinear stability result. 

\subsection{Outline of the argument}\label{sec:argoutline}
We look for solutions to \eqref{eq:main} of the form $(u(t,x),v(t,x))=(Q_*(x)+P(t,x),v(t,x))$, which yields
\bqs
\partial_t \left(\begin{matrix} P \\  v\end{matrix}\right)=\left( \begin{matrix} \cL_u & \beta \\ 0 & \cL_v \end{matrix}\right)\left(\begin{matrix} P \\ v\end{matrix}\right)+\left(\begin{matrix} \mathcal{N}(P) \\  0\end{matrix}\right), \quad t>0, \quad x\in\R.\label{eq:perturb}
\eqs
where
$
\cL_u \coloneqq d \partial_{xx}+s_*\partial_x+f'(Q_*)$, $\cL_v\coloneqq -(\partial_{xx}+1)^2+s_*\partial_x+\mu,
$
and 
\bqs
\mathcal{N}(P)\coloneqq f(Q_*+P)-f(Q_*)-f'(Q_*)P=-3\alpha Q_*P^2-\alpha P^3.
\eqs
For the linear system 
\bqq
\partial_t \left(\begin{matrix} P \\  v\end{matrix}\right)=\left( \begin{matrix} \cL_u & \beta \\ 0 & \cL_v \end{matrix}\right)\left(\begin{matrix} P \\ v\end{matrix}\right), \quad t>0, \quad x\in\R,
\label{eqlin}
\eqq
we define its associated pointwise $2\times2$ matrix Green's function $\bG_\lambda(x,y)$ whose components satisfy 
\bqq
\left( \begin{matrix} \cL_u-\lambda & \beta \\ 0 & \cL_v-\lambda \end{matrix}\right)\left( \begin{matrix} \bG_\lambda^{11}(x,y) &\bG_\lambda^{12}(x,y) \\ \bG_\lambda^{21}(x,y) &\bG_\lambda^{22}(x,y) \end{matrix}\right)=\left( \begin{matrix} -\delta(x-y) & 0 \\ 0 & -\delta(x-y) \end{matrix}\right), \quad \forall (x,y)\in\R^2. \label{eq:ptwsGdefn}
\eqq
Inverse Laplace transform formally gives the $2\times2$ temporal matrix Green's function $\cG(t,x,y)$ with
\bqs
\cG^{ij}(t,x,y)\coloneqq \frac{1}{2\pi \mbi}\int_\Gamma \me^{\lambda t}\bG_\lambda^{ij}(x,y) \md \lambda,
\eqs
for an appropriate  contour $\Gamma\subset \C$. Note that $\bG_\lambda^{21}(x,y)\equiv 0$ by the upper triangular structure. 
Solutions of  \eqref{eq:perturb} with initial condition $(P_0,v_0)$ can be expressed implicitly with Duhamel's formula,
\bqs
\left(\begin{matrix} P \\  v\end{matrix}\right)(t,x)=\int_\R \cG(t,x,y) \left(\begin{matrix} P_0 \\  v_0\end{matrix}\right)(y)\md y+\int_0^t \int_\R \cG(t-s,x,y) \left(\begin{matrix} \mathcal{N}(P) \\  0\end{matrix}\right)(s,y)\md y \md s.
\eqs
Since $\cG^{ij}\equiv 0$, again from the triangular structure, 
\bqs
v(t,x)=\int_\R \cG^{22}(t,x,y)v_0(y)\md y.
\eqs
As a consequence, the perturbation $P(t,x)$ satisfies the integral equation
\bqq
P(t,x)=\int_\R \cG^{11}(t,x,y)P_0(y)\md y+\int_0^t\int_\R \cG^{11}(t-s,x,y)\mathcal{N}(P)(s,y)\md y\md s+\int_\R \cG^{12}(t,x,y)v_0(y)\md y.
\label{eqintsol}
\eqq
Suppose for the moment that the last summand vanishes, $v_0(y)\equiv 0$. We can then follow the analysis in \cite{faye19}, where we have studied the linear problem 
\bqs
\partial_t P = \cL_u P, \quad P(0)=P_0,
\eqs
on weighted spaces with weight function $\omega(x)$ from \eqref{eqomega}. In the weighted variables $P=\omega p$,  $P_0=\omega p_0$, we find
\bqs
\partial_t p = \omega^{-1} \cL_u (\omega p)\eqqcolon \cL p, \quad p(0)=p_0.
\eqs
with solution represented by the associated temporal Green's function $\widetilde{\cG}^{11}(t,x,y)$ as 
\bqs
p(t,x)=\int_\R \widetilde{\cG}^{11}(t,x,y) p_0(y)\md y, \qquad \cG^{11}(t,x,y)=  \widetilde{\cG}^{11}(t,x,y) \frac{\omega(x)}{\omega(y)}, \quad t>0, \quad \forall (x,y)\in\R^2.
\eqs 
Adding nonlinear terms, but still setting $v_0(y)\equiv 0$, one finds 
\bqs
\underset{x\in\R}{\sup}~\left(\frac{1}{1+|x|}\frac{|P(t,x)|}{\omega(x)}\right)\leq \frac{C\epsilon}{(1+t)^{3/2}}, \quad t>0,
\eqs
for some constants $C>0$ and $\epsilon>0$ small enough, controlling a weighted norm of $P_0$; see \cite{faye19}.
The main difficulty here stems from that the fact that the inhomogeneous source term 
\bqs
H(t,x)\coloneqq\int_\R \cG^{12}(t,x,y)v_0(y)\md y,
\eqs
is not bounded in time when multiplied with the weight $\omega(x)$. In fact,  $H(t,x)$ we will decompose
\bqs
H(t,x)\eqqcolon\omega(x) h(t,x)+E(t,x),
\eqs
where $E(t,x)$ vanishes for small time $E(t,x)\equiv 0$ for $t\leq 1$, has weak spatial decay, in the sense that  $\omega^{-1}E$ cannot be controlled in time, but has exponential temporal decay, uniformly in space. It then follows that
\bqs
P(t,x)-E(t,x)=\int_\R \cG^{11}(t,x,y)P_0(y)\md y+\omega h(x)+\int_0^t\int_\R \cG^{11}(t-s,x,y)\mathcal{N}(P)(s,y)\md y\md s.
\eqs
We can therefore redefine our weighted perturbation  $P(t,x)-E(t,x)\eqqcolon\omega(x) p(t,x)$, solving
\bqs
p(t,x)=\int_\R \widetilde{\cG}^{11}(t,x,y)p_0(y)\md y+h(t,x)+\int_0^t\int_\R \widetilde{\cG}^{11}(t-s,x,y)\omega(y)^{-1}\mathcal{N}(\omega p+E)(s,y)\md y\md s.
\eqs
We now sketch how each term in the above equality can be estimated for $t>1$, ignoring the less pertinent difficulties for small times for now. Recall from \cite{faye19}, that
\bqq
\left|\int_\R \widetilde{\cG}^{11}(t,x,y)p_0(y)\md y\right| \lesssim  \frac{1+|x|}{(1+t)^{3/2}} \int_\R (1+|y|)|p_0(y)|\md y, \quad t>1, \quad x\in\R,
\label{eq:CrucialEstimate}
\eqq
for any $p_0$.  With these estimates, we can easily obtain stability following the analysis in \cite{faye19}.

The key ingredient to the analysis here is then to show that $h(t,x)$ satisfies an estimate similar to \eqref{eq:CrucialEstimate}, 
\bqs
|h(t,x)|\lesssim  \frac{1+|x|}{(1+t)^{3/2}} \|v_0\|_\X, \quad t>1, \quad x \in\R, 
\eqs
which we shall obtain in Lemma~\ref{lem:Estimate_h}. Establishing such bounds, we rely on pointwise estimates for $\cG^{12}(t,x,y)$, which in turn will be deduced from pointwise bounds on $\bG_\lambda^{12}(x,y)$, solution  to
\bqs
(\cL_u-\lambda)\bG_\lambda^{12}(x,y)+\beta \bG_\lambda^{22}(x,y)=0, \quad \forall (x,y)\in\R^2,
\eqs
for $\lambda \in \C$ to the right of the essential spectrum of $\cL_u$ and $\cL_v$. In the next section, we obtain estimates on the pointwise Green's function $\bG_\lambda^{12}(x,y)$ and $\bG_\lambda^{22}(x,y)$ for small, medium and large values of $\lambda$. Then in Section~\ref{sec:temporal}, we use these estimates to derive pointwise bounds on $\cG^{12}(t,x,y)$. Finally, we prove our main nonlinear stability result in the last section.

\subsection{Estimates on the pointwise Green's function $\bG_\lambda(x,y)$}\label{sec:pointwise}

\begin{figure}[!t]
\centering
\includegraphics[width=0.75\textwidth]{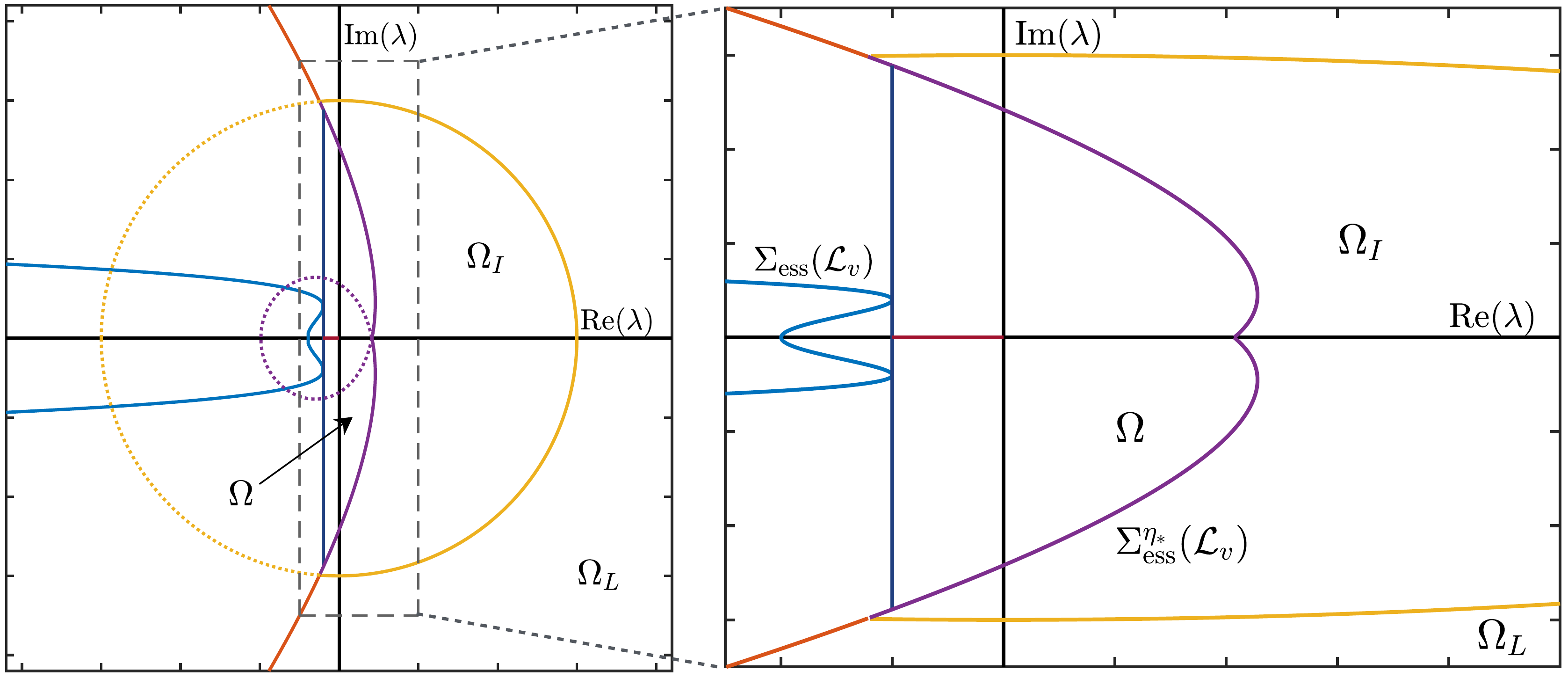}
\caption{Regions $\Omega$, $\Omega_I$ and $\Omega_L$ from  \eqref{e:Omega}, \eqref{e:omegai}, and \eqref{e:omegal}, respectively, showing the ball $\{|\lambda|\leq M_l\}$ (yellow), the boundary between  $\Omega_L$ and $\Omega_I$ (yellow solid),    $\Sigma^{\eta_*}_{\mathrm{ess}}(\mathcal{L}_v)$ (purple and red) separating $\Omega$ and $\Omega_I$ (purple solid). On the left, $\Omega$ is bounded by $\{\Re\lambda=\mu\}$ (dark blue) which bounds the real part of  $\Sigma_{\mathrm{ess}}(\mathcal{L}_v)$ (blue).}
\label{fig:Omega}
\end{figure}

We obtain estimates on the pointwise Green's function $\bG_\lambda(x,y)$ from (\ref{eq:ptwsGdefn}).  Estimates will be divided into three types: ``large'' $\lambda$ estimates valid for $|\lambda|$ sufficiently large, ``small'' $\lambda$ estimates valid for $\lambda$ near the origin and ``intermediate'' $\lambda$ estimates valid in between.  Estimates in the ``large'' $\lambda$ regime are obtained in Lemma~\ref{lem:largelambda} and are obtained by a rescaling of (\ref{eq:ptwsGdefn}) following the approach in \cite{zumbrun98} and are valid on the following subset of the complex plane,
\begin{equation} \Omega_{L}=\left\{\lambda\in\mathbb{C}\ | \ |\lambda|\geq M_l, \ |\mathrm{arg}(\lambda)|<\frac{\pi}{2}+\delta \right\}, \label{e:omegal}\end{equation}
for some $M_l$ sufficiently large.  The intermediate $\lambda$ regime will consist of $\lambda$ in the set
\begin{equation}\Omega_{I}=\left\{\lambda\in\mathbb{C}\ | \ |\lambda|\leq M_l , \ \lambda \ \text{to the right of } \ \Sigma^{\eta_*}_{\text{ess}}(\mathcal{L}_v)  \right\}. \label{e:omegai}\end{equation}
The majority of our effort will be towards estimates for the small $\lambda$ regime, concerned with $\lambda\in\Omega\subset\C$ defined through
\begin{equation}\label{e:Omega}
 \Omega=\{\lambda\in\C|\, |\mathrm{arg}(\lambda)|<\pi,\ \mathrm{Re}(\lambda)\geq  \mu, \ |\lambda|<  M_l, \ \lambda\not\in\Omega_I\}.
\end{equation}
In words, $\Omega$ contains points to the right of $\mu$, to the left of $\Omega_I$, and off the negative real axis. These regions are depicted in Figure~\ref{fig:Omega}.  Note that the ``small'' $\lambda$ regime will actually encompass a rather sizable portion of the complex plane, rather than  a ``small'' neighborhood of the origin. 

Next, let
\bqq  -\gamma_v = \max_{\lambda\in\Omega} \max_{j=1,2} \mathrm{Re}(\nu_j(\lambda)) <0. \label{eq:gammav} \eqq
Our main result will hold for the following set of parameters.

\begin{defi}[Region of validity]\label{defiPi}Let $\Pi\subset \mathcal{R}_{\mathrm{rem}}\cup \mathcal{R}_{\mathrm{abs}}$ denote the set of parameters $(d,\alpha,\mu)$ such that the decay rate condition
\begin{equation}
  3\gamma_v >\frac{s_*}{2d}, \label{e:decayrate}
\end{equation}
holds, where we recall that $s_*=2\sqrt{d\alpha}$ and $\gamma_v$ is defined in (\ref{eq:gammav}).
\end{defi}

\begin{figure}[ht!]
\centering
\includegraphics[width=0.31\textwidth]{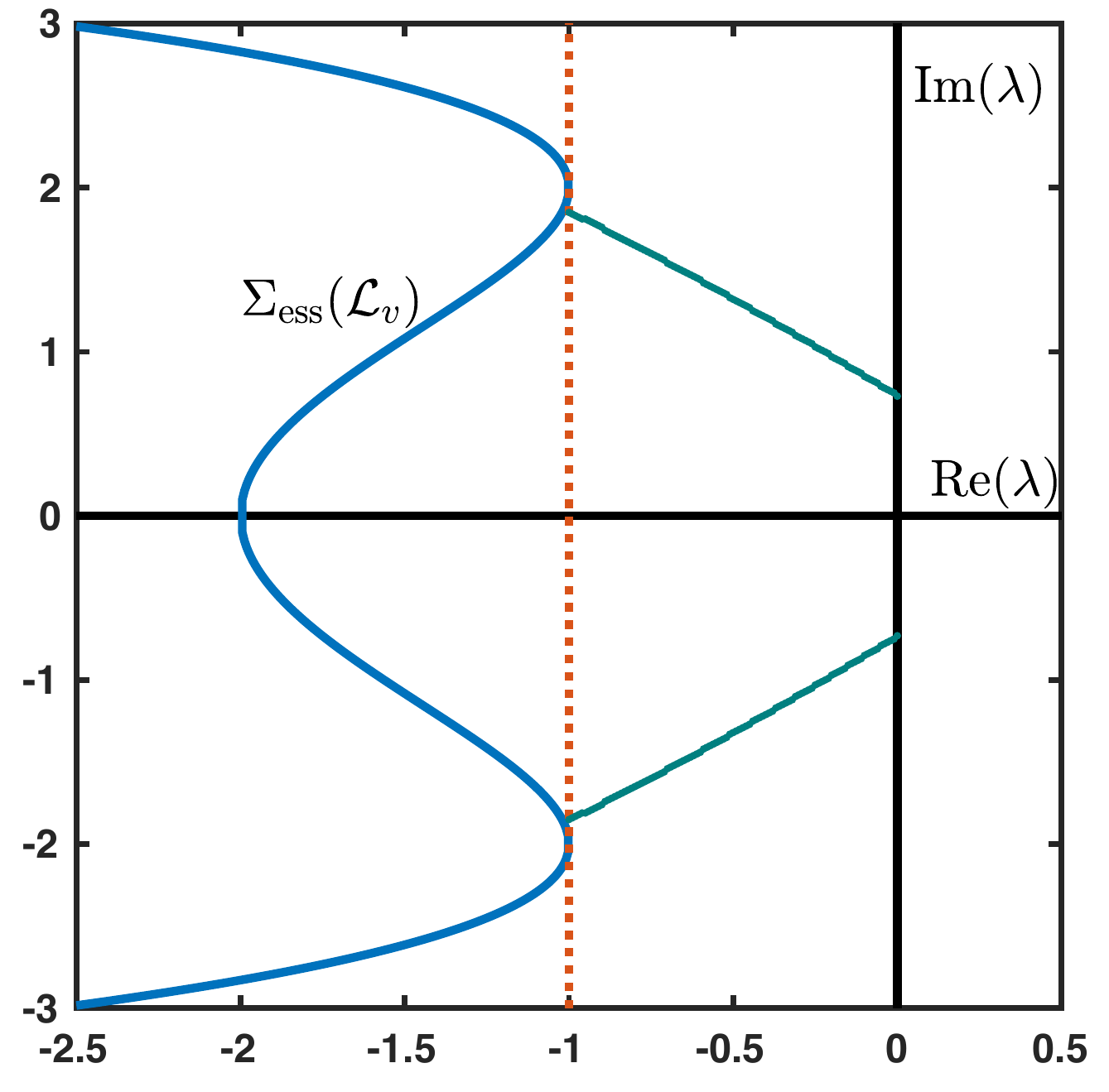}\qquad\quad
\includegraphics[width=0.3\textwidth]{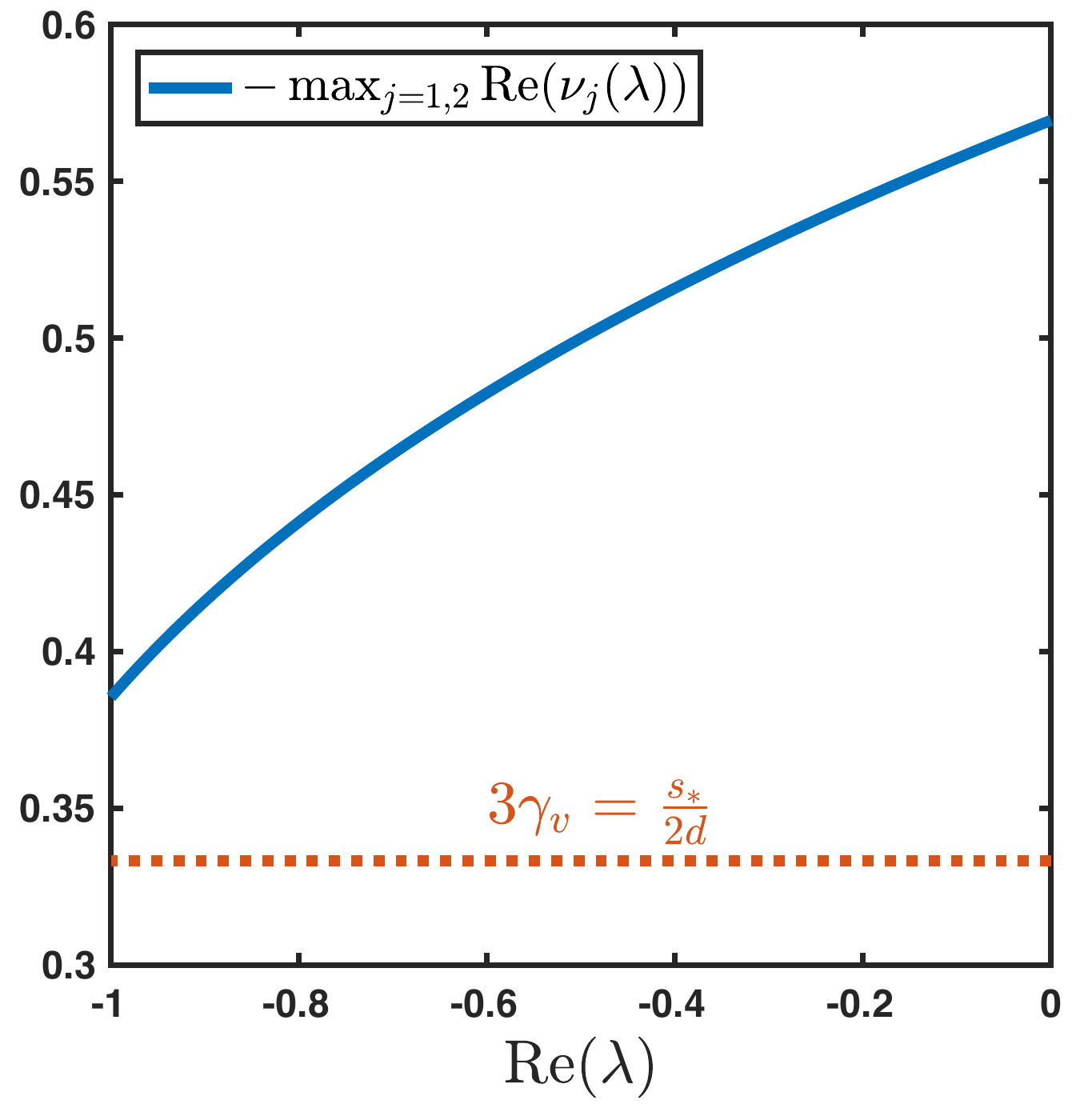}
\caption{Checking \eqref{e:decayrate} numerically. Left: Values of $\lambda$ where  $\max_{\Im\lambda}\max_{j=1,2}\mathrm{Re}(\nu_j(\lambda))$ is attained for fixed $-1<\mathrm{Re}(\lambda)<0$ (dark cyan). Right: Values of $-\max_{\Im\lambda}\max_{j=1,2}\mathrm{Re}(\nu_j(\lambda))$ (blue), always above  $3\gamma_v>\frac{s_*}{2d}$ (dotted orange) as required in \eqref{e:decayrate}; parameters fixed at $d=\alpha=1$ and $\mu=-1$.}
\label{fig:CheckingHypothesis}
\end{figure}

\begin{rmk}[Comments on region of validity] Inspecting the proof, this condition is needed to control nonlinear terms. It is conceivable that a violation of this condition induces faster spreading speed based on resonant interaction mechanisms studied in \cite{faye17}. On the other hand, we show in Section~\ref{sec:absoluteformain} that $\mathcal{R}_{\mathrm{rem}}\cup \mathcal{R}_{\mathrm{abs}}$ have nonempty interior. The (open) decay rate condition is more difficult to verify analytically but easy to verify for specific parameter values; see Figure~\ref{fig:CheckingHypothesis}. On the other hand, taking $\mu\approx \mu^{abs}_0$, we shall see in  Section~\ref{sec:absoluteformain} that the decay rate condition holds by continuity and that therefore so we are guaranteed that the set $\Pi$ of validity is nonempty.
\end{rmk}

\subsubsection{Construction of $\bG_\lambda^{11}(x,y)$ in $\Omega$}

Recall the two dispersion relations for the $u$ component given by
\bqs
D_u^0(\lambda,\nu)\coloneqq d\nu^2+s_*\nu+f'(0)-\lambda, \quad D_u^1(\lambda,\nu)\coloneqq d\nu^2+s_*\nu+f'(1)-\lambda,
\eqs
with corresponding roots $\nu_\pm^0(\lambda)\coloneqq -\frac{s_*}{2d}\pm \sqrt{\frac{\lambda}{d}}$ and $\nu_\pm^1(\lambda)\coloneqq -\frac{s_*}{2d}\pm \frac{\sqrt{s_*^2-4d(f'(1)-\lambda)}}{2d}$.
\begin{lem}\cite[Lemma 2.2]{faye19} \label{lem:phipsi} For $\lambda\in\Omega$, there exists solutions $\vp^\pm$ and $\psi^\pm$ of $\cL_u p=\lambda p$ with
\begin{align*}
\vp^+(x)&=\me^{\nu_-^0(\lambda)x}\left(1+\theta_+(x,\lambda)\right), \ x>0, \\
\psi^+(x)&=\me^{\nu_+^0(\lambda)x}\left(1+\kappa_+(x,\lambda)\right), \ x>0,\\
\vp^-(x)&=\me^{\nu^1_+(\lambda)x}\left(1+\theta_-(x,\lambda)\right), \  x<0, \\
\psi^-(x)&=\me^{\nu^1_-(\lambda)x}\left(1+\kappa_-(x,\lambda)\right), \ x<0,
\end{align*}
where $\theta^\pm(x)$ and $\kappa^\pm(x)$ are both $\O(\me^{-\vartheta |x|})$.
\end{lem}

Thus, $\vp^\pm(x)$ represent bounded solutions of $\cL p=\lambda p$ on either half-line while $\psi^\pm$ represent choices of unbounded (on $\mathbb{R}^-$) or weakly decaying (on $\mathbb{R}^+$) solutions.  The pointwise Green's function is then given by 
\bqq 
\bG_\lambda^{11}(x,y)=\left\{ \begin{array}{cc} \frac{\vp^+(x)\vp^-(y)}{\W_\lambda(y)} & x>y \\ \frac{\vp^-(x)\vp^+(y)}{\W_\lambda(y)} & x<y \end{array}\right. , 
\eqq 
where $\W_\lambda(y)$ is the Wronskian of $\vp^+(x)$ and $\vp^-(x)$, known as the Evans function. An important ingredient to \cite{faye19} is the fact that $\W_0(y)\neq 0$ reflecting  absence of an embedded translational eigenvalue.

\subsubsection{Construction of $\bG_\lambda^{22}(x,y)$ to the right of $\Sigma_{\mathrm{ess}}(\cL_v)$}

Recall that $\bG_\lambda^{22}(x,y)$ satisfies, in the distributional sense, the equation
\bqs
(\cL_v-\lambda)\bG_\lambda^{22}(x,y)=-\delta(x-y), \quad \forall (x,y)\in\R^2.
\eqs
For $\lambda\in\C$ to the right of $\Sigma_{\mathrm{ess}}(\cL_v)$,  we have that
\bqq \bG_\lambda^{22}(x,y)=\left\{ \begin{array}{cc} c_1(\lambda) \me^{\nu_1(\lambda)(x-y)} +c_2(\lambda) \me^{\nu_2(\lambda)(x-y)}, & x>y, \\  -c_3(\lambda) \me^{\nu_3(\lambda)(x-y) }-c_4(\lambda) \me^{\nu_4(\lambda)(x-y)}, & x<y, \end{array}\right.  \label{eq:G22} \eqq
where 
\[ c_i(\lambda)=\frac{-1}{\prod_{j\neq i} (\nu_i(\lambda)-\nu_j(\lambda)) },  \qquad \mathrm{Re}~\nu_{1,2}(\lambda)<0,\quad \mathrm{Re}~\nu_{3,4}(\lambda)>0,\]
and 
\[
 D_v(\lambda,\nu_j(\lambda))=0, \qquad D_v(\lambda,\nu)=-(1+\nu^2)^2+s_*\nu+\mu-\lambda.
\]

\subsubsection{Outline of approach for obtaining estimates on $\bG_\lambda^{12}(x,y)$ for $\lambda\in\Omega$ }
The component $\bG_\lambda^{12}(x,y)$ obeys the following second order inhomogeneous differential equation
\bqq \left(d\partial_{xx}+s_*\partial_x+f'(Q_*(x))-\lambda\right)\bG_\lambda^{12}(x,y)=-\beta \bG_{\lambda}^{22}(x,y), \quad \forall x\in\R. \label{eq:inhomo} \eqq
To motivate our approach, let us first note that we would like to express the solution using the variation of constants formula using the bounded solutions of the homogeneous equation $\vp^\pm(x)$ (see Lemma~\ref{lem:phipsi}) as follows, 
\[ \bG_{\lambda}^{12} (x,y)=-\beta \vp^+(x) \int_{-\infty}^x \frac{\vp^-(\tau)}{\W_\lambda(\tau)} \bG_\lambda^{22}(\tau,y)\md \tau +\beta\vp^-(x) \int_x^\infty  
\frac{\vp^+(\tau)}{\W_\lambda(\tau)} \bG_\lambda^{22}(\tau,y)\md \tau. \]
However, for such a representation to be valid we would need $\bG_\lambda^{22}(x,y)$ to decay faster than $\me^{-\frac{s^*}{2d}x}$ so that the second integral would converge.  This condition fails exactly as  $\lambda$ enters the weighted essential spectrum $\Sigma^{\eta_*}_{\text{ess}}(\mathcal{L}_v)$.  We overcome this issue by first extracting a leading order description of $\bG_\lambda^{12}(x,y)$ in which terms with insufficient decay can be isolated.  Subsequently, the correction term can be derived using a variation of constants formula.    

To begin this process, we define
\bqs
\cL_u^\infty\coloneqq\left\{ \begin{array}{cc}\cL_u^+= d\partial_{xx}+s_*\partial_x+f'(0), & x>0, \\ \cL_u^-= d\partial_{xx}+s_*\partial_x+f'(1), & x<0, \end{array}\right.
\eqs
and then decompose 
\bqs
\bG_\lambda^{12}(x,y)=\bG_\lambda^{12,\infty}(x,y)+\widetilde{\bG_\lambda^{12}}(x,y),
\eqs
where
\[ \left (\mathcal{L}_u^\infty-\lambda\right)\bG_\lambda^{12,\infty}(x,y)=-\beta \bG_\lambda^{22}(x,y). \]
The remaining term in the decomposition is then defined by 
\[ \left (\mathcal{L}_u-\lambda\right)\widetilde{\bG_\lambda^{12}}(x,y)= -\left(f'(Q_*(x)-f_\infty(x)\right) \bG_\lambda^{12,\infty}(x,y).\]
Comparing this equation to (\ref{eq:inhomo}) we see that the exponential decay is stronger here due to the prefactor $(f'(Q_*(x))-f_\infty(x))$.   We will show that solutions of this equation can be represented using  the variation of constants formula,
\bqq \widetilde{\bG_{\lambda}^{12}} (x,y)= \vp^+(x) \int_{-\infty}^x \frac{\vp^-(\tau)}{\W_\lambda(\tau)} h(\tau)\bG_\lambda^{12,\infty}(\tau,y)\md \tau -\vp^-(x) \int_x^\infty  
\frac{\vp^+(\tau)}{\W_\lambda(\tau)}h(\tau) \bG_\lambda^{12,\infty}(\tau,y)\md \tau, \label{eq:varcontgooddecay} \eqq
where 
\[
h(x)=-\left(f'(Q_*(x)-f_\infty(x)\right), \qquad |h(x)|\leq C \me^{-\vartheta |x|}
. \]
The exponential decay rate of $h(x)$ is in fact equal to that of the traveling front.  For the sake of motivation, assume for the moment that $y<x$ with $x>0$.  Inspection of the second integral in (\ref{eq:varcontgooddecay}) reveals that it includes terms bounded by
\[ C\int_x^\infty \me^{-\nu_+^0(\lambda)\tau}\me^{-\vartheta \tau} \me^{\nu_1(\lambda)(\tau-y)}\md \tau \]
for which integrability holds if
\[ \mathrm{Re} \left(\nu_1(\lambda)-\nu_+^0(\lambda)\right)<\vartheta \]
which is the familiar condition from the Gap Lemma; see \cite{gardner98,kapitula98} requiring that the negative spectral gap of the spatial eigenvalues must not exceed the decay rate of the front.  

\subsubsection{The asymptotic Green's function $\bG_\lambda^{12,\infty}(x,y)$}

\begin{lem}\label{lem:G12infinity} There exists constants $b_j(\lambda)$ and $h_j(\lambda)$, $j=1,2,3,4$,  such that $\bG_\lambda^{12,\infty}(x,y)$ admits the following representation
\bqs
\bG_\lambda^{12,\infty}(x,y)=\left\{
\begin{array}{ll}
b_1(\lambda,y)\me^{\nu_-^0(\lambda)x}-\sum_{j=1}^2 \frac{\beta}{D_u^0(\lambda,\nu_j(\lambda))}c_j(\lambda) \me^{\nu_j(\lambda)(x-y)} ,&  x>y>0, \\
-b_2(\lambda,y)\me^{\nu_-^0(\lambda)x}-b_3(\lambda,y)\me^{\nu_+^0(\lambda)x}+\sum_{j=3}^4 \frac{\beta}{D_u^0(\lambda,\nu_j(\lambda))}c_j(\lambda) \me^{\nu_j(\lambda)(x-y)},& y>x>0, \\
b_4(\lambda,y)\me^{\nu_+^1(\lambda)x}+\sum_{j=3}^4 \frac{\beta}{D_u^1(\lambda,\nu_j(\lambda))}c_j(\lambda) \me^{\nu_j(\lambda)(x-y)} ,& x<0<y, \\
h_1(\lambda,y)\me^{\nu_-^0(\lambda)x}-\beta\sum_{j=1}^2\frac{c_j(\lambda)}{D_u^0(\lambda,\nu_j(\lambda))}\me^{\nu_j (\lambda)(x-y)},&y<0<x,\\
-h_2(\lambda,y)\me^{\nu_-^1(\lambda)x}-h_3(\lambda,y)\me^{\nu_+^1(\lambda)x}-\beta\sum_{j=1}^2\frac{c_j(\lambda)}{D_u^1(\lambda,\nu_j(\lambda))}\me^{\nu_j (\lambda)(x-y)},&y<x<0,\\
h_4(\lambda,y)\me^{\nu_+^1(\lambda)x}+\beta\sum_{j=3}^4\frac{c_j(\lambda)}{D_u^1(\lambda,\nu_j(\lambda))}\me^{\nu_j (\lambda)(x-y)},&x<y<0.
\end{array}
\right. 
\eqs

\end{lem}
The ordinary differential equation defining $\bG_\lambda^{12,\infty}(x,y)$ is piecewise constant-coefficient and therefore explicit solution formulas are available; we present these calculations in Appendix~\ref{sec:G12infinityproof}. We remark that the dependence of the various constants on $y$ and $\lambda$ will be important in the subsequent analysis and we provide formulas for them in (\ref{eq:hs}) and (\ref{eq:bs}).

\begin{rmk}\label{rmk:singularities} To obtain temporal bounds for the asymptotic system $\partial_t u=\mathcal{L}_u^\infty u+\beta v$, the inverse Laplace transform of the  pointwise Green's function presented in Lemma~\ref{lem:G12infinity} are required.  Here singularities are paramount.  Observe that singularities of $\bG_\lambda^{12,\infty}(x,y)$ arise in several places.  First, the coefficients $c_j(\lambda)$ have singularities at branch points of $D_v(\lambda,\nu)$.  We denote the branch point as $\lambda_v^\textnormal{bp}$ (together with its complex conjugate) explicit formulas for which are provided in (\ref{eq:vbranch}).  Additional singularities come from resonance poles.  These arise from pinched double roots where $\nu_{1,2}(\lambda^\textnormal{rp})=\nu_0^+(\lambda^\textnormal{rp})$ so that $D_u^0(\lambda^\textnormal{rp},\nu_{1,2}(\lambda^\textnormal{rp}))=0$ (see the formulas for $\lambda_{\pm}^{\textnormal{rp}_-}$ from Section~\ref{sec:absoluteformain}).  Singularities stemming from non-pinched double roots are removable; see Remark~\ref{rmk:nosingularity}.  

Also note that the absolute spectrum, aside from contributions at double roots, does not lead to any singularities in $\bG_\lambda^{12,\infty}(x,y)$.  In the next lemma, we will show that the absolute spectrum is not seen in the singularities of the the full Green's function $\bG_\lambda^{12}(x,y)$.
\end{rmk}

\begin{rmk}\label{rmk:zeta} 
In a neighborhood of $\lambda_v^\textnormal{bp}$, roots of $D_v(\lambda,\nu)=0$ have leading order expansions 
\[ \left(\nu-\nu_v^\textnormal{bp}\right)^2=-\frac{2 \partial_\lambda D_v(\lambda_v^{\textnormal{bp}},\nu_v^\textnormal{bp})}{\partial_{\nu\nu} D_v(\lambda_v^{\textnormal{bp}},\nu_v^\textnormal{bp})} \ \left (\lambda-\lambda_v^\textnormal{bp}\right). \]
Let 
\[ \zeta\coloneqq-\frac{2 \partial_\lambda D_v(\lambda_v^{\textnormal{bp}},\nu_v^\textnormal{bp})}{\partial_{\nu\nu} D_v(\lambda_v^{\textnormal{bp}},\nu_v^\textnormal{bp})}=\frac{-1}{6\left(\nu_v^{\textnormal{bp}}\right)^2+2}, \] 
and note that $\Re (\zeta) >0$ as $\lambda_v^\textnormal{bp}$ is the rightmost part of $\Sigma_{\textnormal{abs}}(\cL_v^+)$.
\end{rmk}

\subsubsection{Estimates of $\bG_\lambda^{12}(x,y)$ of $\lambda\in\Omega$ }
With formulas for $\bG_\lambda^{12,\infty}(x,y)$ in hand, we proceed to derive estimates on $\widetilde{\bG_\lambda^{12}}(x,y)$ using
\[ \left (\mathcal{L}_u-\lambda\right)\widetilde{\bG_\lambda^{12}}(x,y)= -\left(f'(Q_*(x)-f_\infty(x)\right) \bG_\lambda^{12,\infty}(x,y).\]  
In the following lemma, we compile estimates for $\bG_\lambda^{12}(x,y)$, demonstrating that the terms coming from $\widetilde{\bG_\lambda^{12}}(x,y)$ are in fact higher order.

\begin{lem}\label{lem:G12full} For $\lambda\in\Omega$ with $\Im(\lambda)\geq 0$, there exists functions $\boldH_\lambda(x,y)$, $\bI_{\lambda,j}(x,y)$ and $\boldJ_{\lambda,j}(x,y)$, bounded uniformly in $(x,y)$ such that the pointwise Green's function $\bG_\lambda^{12}(x,y)$ obeys the following bounds for each of the six arrangements of $x$, $y$ and $0$;
\begin{itemize}
\item  $y<0<x$ 
\begin{align*}
\bG_\lambda^{12}(x,y) & = \frac{1}{\sqrt{\zeta(\lambda-\lambda_v^\textnormal{bp})}}\me^{\nu_-^0(\lambda)x-\nu_-^1(\lambda)y} \boldH_\lambda(x,y) \\
&~~~+\sum_{j=1}^2 \frac{1}{(\lambda-\lambda^{\textnormal{rp}})\sqrt{\zeta(\lambda-\lambda_v^{\textnormal{bp}})}} \left[\me^{\nu_-^0(\lambda)x-\nu_j(\lambda)y}\bI_{\lambda,j}(x,y)+\me^{\nu_j(\lambda)(x-y)}\boldJ_{\lambda,j}(x,y) \right], 
\end{align*}

\item  $0<y<x$ 
\begin{align*}
\bG_\lambda^{12}(x,y) & = \frac{1}{\RPBP}\me^{\nu_-^0(\lambda)(x-y)} \boldH_\lambda(x,y) \\
%&~~~+\sum_{j=3}^4 \frac{1}{\BP} \me^{\nu_-^0(\lambda)x-\nu_j(\lambda)y}\bI_{\lambda}(x,y) \\
&~~~+\sum_{j=1}^2 \frac{1}{\RPBP} \me^{\nu_j(\lambda)(x-y)}\boldJ_{\lambda,j}(x,y),
\end{align*}

\item $x<y<0$
\begin{align*}
\bG_\lambda^{12}(x,y) & = \frac{1}{\BP}\me^{\nu_+^1(\lambda)(x-y)}\mathbf{H}_\lambda(x,y) \\
&~~~+\sum_{j=1}^2\frac{1}{\RPBP}\me^{\nu_+^1(\lambda)x-\nu_{j}(\lambda)y}\bI_{\lambda,j}(x,y) \\
&~~~+\sum_{j=3}^4 \frac{1}{\BP} \me^{\nu_j(\lambda)(x-y)}\boldJ_{\lambda,j}(x,y),
\end{align*}
\item  $y<x<0$ 
\begin{align*}
\bG_\lambda^{12}(x,y) & = \frac{1}{\BP}\me^{\nu_-^1(\lambda)(x-y)}\mathbf{H}_\lambda(x,y) \\
&~~~+\sum_{j=1}^2 \frac{1}{\RPBP} \me^{\nu_+^1(\lambda)x-\nu_j(\lambda) y}\mathbf{I}_{\lambda,j}(x,y) \\
&~~~+\sum_{j=1}^2 \frac{1}{\BP} \me^{\nu_j(\lambda)(x-y)}\boldJ_{\lambda,j}(x,y),
\end{align*}
\item $x<0<y$
\begin{align*}
\bG_\lambda^{12}(x,y) & = \frac{1}{\RPBP}\me^{\nu_+^1(\lambda)x-\nu_+^0(\lambda) y} \mathbf{H}_\lambda(x,y) \\
&~~~+\sum_{j=3}^4 \frac{1}{\BP}\left[\me^{\nu_+^1(\lambda)x-\nu_j(\lambda) y}\mathbf{I}_{\lambda,j}(x,y)+ \me^{\nu_j(\lambda)(x-y)}\boldJ_{\lambda,j}(x,y)\right],
\end{align*}
\item $0<x<y$
\begin{align*}
\bG_\lambda^{12}(x,y) & = \frac{1}{\RPBP}\me^{\nu_+^0(\lambda)(x- y)} \mathbf{H}_\lambda(x,y) \\
&~~~ \sum_{j=3}^4 \frac{1}{\BP}\left[\me^{\nu_-^0(\lambda)x-\nu_j(\lambda) y}\mathbf{I}_{\lambda,j}(x,y)+ \me^{\nu_j(\lambda)(x-y)}\boldJ_{\lambda,j}(x,y)\right] ,
\end{align*}
\end{itemize}
where $\zeta$ is defined in Remark~\ref{rmk:zeta}.  For $\Im(\lambda)<0$, the same estimates hold but with the complex conjugates of $\lambda_v^{\textnormal{bp}}$ and $\lambda^{\textnormal{rp}}$. 
\end{lem}

We present the proof of this lemma in Appendix~\ref{sec:G12full}.

\subsubsection{Large and intermediate $\lambda$ estimates for $\bG_\lambda^{12}(x,y)$}

We now turn our attention to estimates for $\bG_\lambda^{12}(x,y)$ for $\lambda$ in $\Omega_L$ and $\Omega_I$.  

\begin{lem} \label{lem:largelambda} There exists $M_l>0$ and $\delta>0$ such that for any $\lambda$ satisfying  $|\lambda|>M_l$ and for $|\mathrm{arg}(\lambda)|<\frac{\pi}{2}+\delta$ the pointwise Green's functions admit the following bounds for some $\eta>0$,
\begin{align*}
 \left| \bG_\lambda^{11}(x,y)\right|&\leq \frac{C}{\sqrt{|\lambda|}}\me^{-\sqrt{|\lambda|} \eta |x-y|}, \\
 \left| \bG_\lambda^{22}(x,y)\right|&\leq \frac{C}{|\lambda|^{3/4}}\me^{-|\lambda|^{1/4} \eta |x-y|} , \\
  \left| \bG_\lambda^{12}(x,y)\right|&\leq \frac{C}{|\lambda|^{7/4}}\me^{-|\lambda|^{1/4} \eta |x-y|}.
\end{align*}
\end{lem}
\begin{Proof} The estimate for $\bG_\lambda^{11}(x,y)$ was stated in \cite{faye19} and obtained by rescaling of the spatial variable when $|\lambda|$ is large following \cite{zumbrun98}.  A similar argument works for $\bG_\lambda^{22}(x,y)$, although we note that since the linearized operator for $v$ is constant coefficient the expression for $\bG_\lambda^{22}(x,y)$ is exact and the estimate can be obtained directly from (\ref{eq:G22}).   For $\bG_\lambda^{12}(x,y)$ we recall the inhomogenous relation
\[ (\mathcal{L}_u-\lambda)\bG_\lambda^{12}(x,y)=-\beta \bG_\lambda^{22}(x,y).\]
Note that as $|\lambda|$ becomes large the homogeneous equation has spatial eigenvalues with asymptotics $\pm\sqrt{\lambda}$ while the $\nu_j(\lambda)$ scale as $|-\lambda|^{1/4}$.  Therefore as the modulus of $\lambda$ becomes large, the decay of the solution is dominated by the inhomogeneous terms and the estimate follows.
\end{Proof}

\begin{lem}\label{lem:lambdaint} Consider $\Omega_{I}$, the region of the  complex plane to the right of $\Sigma_{\text{ess}}^{-\frac{s_*}{2d}}(\mathcal{L}_v)$ with $|\lambda|\leq M_l$ and with $|\mathrm{arg}(\lambda)|<\frac{\pi}{2}+\delta$.  The pointwise Green's functions obey the following bounds
\bqq \left| \bG_\lambda^{11}(x,y)\right|\leq C\frac{\omega(x)}{\omega(y)} , \quad \left|\bG_\lambda^{22}(x,y)\right|\leq C , \quad \left|\bG_\lambda^{12}(x,y)\right|\leq C\frac{\omega(x)}{\omega(y)} \eqq
\end{lem}
\begin{Proof}
The estimate for $\bG_\lambda^{11}(x,y)$ again follows as in \cite{faye19}.  The estimate for $\bG_\lambda^{22}(x,y)$ holds since $\lambda$-values  considered are contained in a compact set to the right of $\Sigma_{\text{ess}}(\mathcal{L}_v)$.   Finally, for $\bG_\lambda^{12}(x,y)$ we observe that, since $\lambda$ lies in a compact set, the definitions of $\varphi^\pm(x)$ can be extended to this set as we can write the bounded solution satisfying $(\mathcal{L}_u-\lambda)\bG_{\lambda}^{12}(x,y)=-\beta \bG_\lambda^{22}(x,y)$ as 
\[\bG_{\lambda}^{12}(x,y)=  -\beta\varphi^+(x) \int_{-\infty}^x \frac{\varphi^-(\tau)}{\W_\lambda(\tau)} \bG_\lambda^{22}(\tau,y)\md \tau -\beta\varphi^-(x) \int_x^\infty  
\frac{\varphi^+(\tau)}{\W_\lambda(\tau)} \bG_\lambda^{22,\infty}(\tau,y)\md \tau, \]
In contrast to the case of $\lambda\in\Omega$, here the set of $\lambda$ under consideration excludes curves of absolute spectrum.  Thus, the integrals converge.   Inspecting formulas for $\varphi^\pm(x)$, we see that  $\omega(x)$ can be factored from the solution.  
\end{Proof}

\subsection{Estimates on the temporal Green's function $\cG(t,x,y)$} \label{sec:temporal}
In this section, we derive estimates on the temporal Green's function as the inverse Laplace Transform of the pointwise Green's function, 
\bqs
\cG(t,x,y)=\frac{1}{2\pi\mbi}\int_\Gamma \me^{\lambda t}\bG_{\lambda}(x,y)\md\lambda
\eqs
for a suitable contour $\Gamma \subset \C$. We recall estimates on $\cG^{11}(t,x,y)$ from \cite{faye19}.

\begin{lem}\cite[Prop. 4.1]{faye19} \label{timeGreen11}
There exists constants $\kappa_u>0$, $r>0$ and $C>0$, such that the Green's function $\cG^{11}(t,x,y)$  satisfies the following estimates.
\begin{itemize}
\item[(i)] For $|x-y|\geq Kt$ or $0<t<1$,  with $K$ sufficiently large, 
\bqs
|\cG^{11}(t,x,y)|\leq C\frac{1}{t^{1/2}}\frac{\omega(x)}{\omega(y)}\me^{-\frac{|x-y|^2}{\kappa_u t}}.
\eqs
\item[(ii)] For $|x-y|\leq Kt$ and $t\geq 1$, with $K$ as above, 
\bqs
|\cG^{11}(t,x,y)|\leq C\frac{\omega(x)}{\omega(y)}\left(\frac{1+|x-y|}{t^{3/2}} \me^{-\frac{|x-y|^2}{\kappa_u t}}+\me^{-rt}\right).
\eqs
\end{itemize}

\end{lem}
The estimates not sharp but sufficient for our purpose here. We next give estimates on $\cG^{22}(t,x,y)$.
\begin{lem}\label{timeGreen22} There exists positive constants $C$, $\theta$, $\kappa_v$, $\gamma_v$ and $\delta_v$ such that
\begin{align*}
 |\cG^{22}(t,x,y)|&\leq \left\{ \begin{array}{cc} C \me^{-\theta t} \me^{-\gamma_v (x-y)}, & x>y, \\
 C \me^{-\theta t} \me^{\delta_v (x-y)}, & x<y, \end{array}\right.\qquad t\geq 1,\\
 |\cG^{22}(t,x,y)|&\leq C\frac{1}{t^{1/4}}\me^{-\frac{|x-y|^{4/3}}{\kappa_vt^{1/3}}},\qquad 0<t<1.
\end{align*}
\end{lem}
\begin{Proof} Recall the explicit  expression for $\bG^{22}_\lambda(x,y)$ in (\ref{eq:G22}) valid to the right of $\Sigma_{\text{ess}}(\mathcal{L}_v)$.  Note that $\mathcal{L}_v$ is a sectorial operator and generates an analytic semigroup.  The steps required to obtain estimates on the resolvent kernel $\bG_\lambda^{22}(x,y)$ mimic those required for estimates of the semigroup; see \cite{lunardi}. We reproduce those estimates to prepare for the adaptation to estimates on $\bG_\lambda^{12}(x,y)$.  Consider the sectorial contour in the upper half plane,
\[ \Gamma= -\theta +\me^{\mbi \eta}\ell, \quad  \frac{\pi}{2}<\eta<\pi, \quad \ell>0  \]
Singularities of $\bG_\lambda^{22}(x,y)$ occur at double roots of the dispersion relation, none of which are located in $\Omega$.  As a result, this contour can be chosen such that it is contained in the left half of the complex plane but to the right of $\Sigma_{\text{ess}}(\mathcal{L}_v)$.  Along $\Gamma$, there are constants $\gamma_v>0$ and $\delta_v>0$ such that
\[ \mathrm{Re}(\nu_1(\lambda))<-\gamma_v, \quad  \mathrm{Re}(\nu_2(\lambda))<-\gamma_v, \quad  \mathrm{Re}(\nu_3(\lambda))>\delta_v,  \quad \mathrm{Re}(\nu_4(\lambda))>\delta_v. \]
Changing the integration variable to $\ell$, we consider the case of $x>y$ and focus on the term in (\ref{eq:G22}) involving $\nu_1(\lambda)$ where 
\bqs \left|\frac{1}{2\pi\mbi} \int_\Gamma \me^{\lambda t} c_1(\lambda) \me^{\nu_1(\lambda)(x-y)} \md \lambda \right| \leq C \me^{-\theta t} \int_0^\infty \me^{\cos(\eta)\ell t}\me^{-\gamma_1 (x-y)} \md \ell \leq \frac{C}{t} \me^{-\theta t} \me^{-\gamma_v(x-y)},  \eqs
which yields the desired estimate when $t\geq 1$. The case for $x<y$ is analogous.

The estimates for $t<1$ are dominated by the large $\lambda$ behavior and closely resembles  the analysis in \cite[Sec. 5.1]{howard12} for the Cahn-Hilliard equation. We omit details here. 
\end{Proof}

We now turn our attention to $\cG^{12}(t,x,y)$ with small-time,  $\frac{|x-y|}{t}$ large,  estimates in Proposition~\ref{timeGreen12} and with large time,  $\frac{|x-y|}{t}$ bounded, estimates in Proposition~\ref{prop:timegreen12}, our key ingredient. 

\begin{prop}\label{timeGreen12}
There exists positive constants $K>0$, $C>0$ and $\widetilde{\kappa}>0$ such that for $|x-y|\geq Kt$ or $0<t<1$ we have
\[ |\cG^{12}(t,x,y) |\leq C\frac{1}{t^{1/4}}\frac{\omega(x)}{\omega(y)}\me^{-\frac{|x-y|^{4/3}}{\widetilde{\kappa} t^{1/3}}}. \]
\end{prop}

\begin{Proof}
We first with the case $|x-y|\geq Kt$, where $K$ sufficiently large which will be chosen in the proof. Recall the large $\lambda$ estimate for $\bG_\lambda^{12}(x,y)$:
\bqs
\left| \bG_\lambda^{12}(x,y)\right|\leq \frac{C}{|\lambda|^{7/4}}\me^{-|\lambda|^{1/4} \eta |x-y|}, \quad \forall(x,y)\in\R^2,
\eqs
valid for all $|\lambda|>M_l$. Possibly restricting to larger values of $\lambda$, we may assume 
\bqs
\left| \bG_\lambda^{12}(x,y)\right|\leq \frac{C}{|\lambda|^{7/4}}\frac{\omega(x)}{\omega(y)}\me^{-|\lambda|^{1/4} \eta |x-y|}, \quad \forall(x,y)\in\R^2.
\eqs
We now consider a contour that consists of two parts:
\begin{itemize}[itemsep=-0pt,topsep=-0pt]
\item a polynomial contour $\Gamma_p$ parametrized by
\bqs
\lambda_p(\ell)=(\rho+z\ell)^4, \quad \ell \in\R, \quad z\coloneqq \me^{\mbi \frac{\pi}{4}}
\eqs
which lies outside  the ball of radius $\rho^4>M_l$ for some $\rho>0$ which will be chosen later;
\item a linear contour $\Gamma_l$ parametrized by
\bqs
\lambda_l(\ell)=-\delta |\ell| +\mbi \ell, \quad \ell \in\R,\quad \text{ for some }\delta>0.
\eqs
\end{itemize}
We first focus on $\Gamma_p$ and  write  $I_p=\{\ell|\,|\lambda_p(\ell)|\geq \rho^4\}\subset\R$.
Assume for now that $\rho^4>M_l$ such that large $\lambda$ estimates hold. We also note that $|\lambda'_p(\ell)|=4|\lambda(\ell)|^{3/4}$ so that
\bqs
\left|\int_{\Gamma_p} \me^{\lambda t}\bG_\lambda^{12}(x,y) \md \lambda\right| 
\lesssim \frac{\omega(x)}{\omega(y)} \me^{-\rho  \eta |x-y|} \int_{I_p} \me^{\Re(\lambda_p(\ell))t} \frac{\md \ell}{|\lambda_p(\ell)|}.
\eqs
We next use the fact that 
\bqs
\Re(\lambda_p(\ell))=\rho^4-\ell^4+2\sqrt{2}(\rho^3\ell-\rho \ell^3)\leq 2 \rho^4, \quad \ell \in \R,
\eqs
which ensures the existence of $c_0>0$ such that
\bqs
\Re(\lambda_p(\ell))-3\rho^4 \leq -c_0 \ell^4.
\eqs
As a consequence, 
\bqs
\left|\int_{\Gamma_p} \me^{\lambda t}\bG_\lambda^{12}(x,y) \md \lambda\right| 
\lesssim \frac{\omega(x)}{\omega(y)} \me^{3\rho^4t-\rho  \eta |x-y|} \int_{I_p} \me^{-c_0 \ell^4t}\md \ell,
\eqs
where we also used the fact that
\bqs
\frac{1}{{|\lambda_p(\ell)|}}\leq \frac{1}{M_l}, \quad \ell \in I_p.
\eqs
Notice that the remaining integral can easily be estimated as 
\bqs
\int_{I_p} \me^{-c_0 \ell^4t}\md \ell \leq \int_{\R} \me^{-c_0 \ell^4t}\md \ell  \lesssim \frac{1}{t^{1/4}}.
\eqs
We now choose $\rho$ as
\bqs
\rho\coloneqq \frac{1}{L} \left(\frac{|x-y|}{t} \right)^{1/3},
\eqs
where the constant $L>0$ is chosen such that
\bqs
3\rho^4t-\rho  \eta |x-y| = \frac{|x-y|^{3/4}}{t^{1/3}}\left(\frac{3}{L^4}-\frac{\eta}{L} \right)=-\frac{\eta}{2L}\frac{|x-y|^{3/4}}{t^{1/3}},
\eqs
that is, $L\coloneqq(6/\eta)^{1/3}$. Then the condition $\rho^4\geq M_l$ becomes
\bqs
\frac{|x-y|}{t} \geq M_l^{3/4}L^3.
\eqs
As a consequence, $K\coloneqq M_l^{3/4}L^3$ and our assumption that $|x-y|\geq Kt$ guarantees $\rho^4\geq M_l$. 

Next, for the contour $\Gamma_l$ we denote by $\ell_*>0$ the largest value for which $\lambda_p(\ell_*)=\lambda_l(\ell_*)$. We define $I_l\coloneqq(-\infty,-\ell_*)\cup(\ell_*,+\infty)$ such that
\bqs
\Gamma_l=\left\{\lambda_l(\ell)~|~ \ell \in I_l \right\},
\eqs
and  we have for all $\ell\in I_l$
\bqs
|\lambda_l(\ell)|=|\ell| \sqrt{1+\delta^2}, \quad |\lambda_l(\ell)|\geq |\lambda_l(\ell_*)|=|\lambda_p(\ell_*)|\geq \rho^4>M_l.
\eqs
As a consequence, using the symmetry of $\Gamma_l$, we obtain
\begin{align*}
\left|\int_{\Gamma_l} \me^{\lambda t}\bG_\lambda^{12}(x,y) \md \lambda\right| 
&\lesssim \frac{\omega(x)}{\omega(y)} \me^{-\rho \eta |x-y|} \int_{I_l} \me^{\Re(\lambda_l(\ell))t}\frac{\md\ell}{|\lambda_l(\ell)|^{7/4}},\\
&\lesssim \frac{\omega(x)}{\omega(y)} \me^{-\frac{\eta}{L} \frac{|x-y|^{4/3}}{t^{1/3}}} \int_{\ell_*}^{+\infty} \me^{-\delta \ell t}\frac{\md \ell}{\ell}\\
&\lesssim \frac{\omega(x)}{\omega(y)} \me^{-\frac{\eta}{L} \frac{|x-y|^{4/3}}{t^{1/3}}} \me^{-\delta \ell_* t}.
\end{align*}

Finally, combining our estimates along $\Gamma_p$ and $\Gamma_l$, we have proved that for all $|x-y|\geq Kt$,
\bqs
 |\cG^{12}(t,x,y) |\leq C\frac{1}{t^{1/4}}\frac{\omega(x)}{\omega(y)}\me^{-\frac{\eta}{2L}\frac{|x-y|^{4/3}}{\widetilde{\kappa} t^{1/3}}}.
\eqs

We now turn our attention to the small time estimate where $0<t<1$. We first remark that if $|x-y|\geq Kt$ then the previous pointwise bound holds, and we are done. As a consequence, from now on we assume that $0<t<1$ and $|x-y|<Kt$. In that case, we consider a contour composed of two parts:
\begin{itemize}[itemsep=-0pt,topsep=-0pt]
\item a portion of a circle $\Gamma_c$ of radius $R>M_l$ parametrized by
\bqs
\lambda_c(\theta) = R \me^{\mbi \theta}, \quad \theta \in [-\theta_*,\theta*],\quad
\text{for some }\pi/2<\theta_*<\pi/2+\epsilon \text{ with } \epsilon>0;
\eqs
\item a linear contour $\Gamma_l$ parametrized by
\bqs
\lambda_l(\ell)=-\delta |\ell| +\mbi \ell, \quad \ell \in\R,
\quad \text{ for some }\ \delta>0.
\eqs
\end{itemize}
Note that both contours intersect at $\lambda_*$ and $\overline{\lambda_*}$ with
\bqs
\lambda_* = R \me^{\mbi \theta_*} = -\delta \ell_* +\mbi \ell_*, \quad \ell_*=\sqrt{1+\delta^2}/R.
\eqs
Along the contour $\Gamma_c$, we have that
\bqs
\left|\int_{\Gamma_c} \me^{\lambda t}\bG_\lambda^{12}(x,y) \md \lambda\right| 
\lesssim \frac{\omega(x)}{\omega(y)} \me^{-R\eta |x-y|} \int_{\Gamma_c} \me^{\Re(\lambda)t} \frac{\md \lambda }{|\lambda|^{7/4}}.
\eqs
As $0<t<1$, the above integral is bounded by some fixed constant which depends on $R$. On the other hand, we have that $|x-y|<Kt$ so that
\bqs
\frac{R\eta}{K^{1/3}}\frac{|x-y|^{1/3}}{t^{1/3}}|x-y|<R\eta |x-y|,
\eqs
and 
\bqs
\left|\int_{\Gamma_c} \me^{\lambda t}\bG_\lambda^{12}(x,y) \md \lambda\right| 
\lesssim \frac{\omega(x)}{\omega(y)} \me^{-\frac{R\eta}{K^{1/3}}\frac{|x-y|^{4/3}}{t^{1/3}}}\leq \frac{1}{t^{1/4}}  \frac{\omega(x)}{\omega(y)} \me^{-\frac{R\eta}{K^{1/3}}\frac{|x-y|^{4/3}}{t^{1/3}}}.
\eqs
Along $\Gamma_l$, we obtain a slightly better estimate which can be subsumed into the previous one.

We finally note that 
\bqs
\frac{K^{1/3}}{R}=\frac{M_l^{1/4}L}{R}<\frac{L}{M_l^{3/4}}<2L,
\eqs 
provided that $M_l$ is large enough and so we let $\widetilde{\kappa}\coloneqq 2L/\eta$. This concludes the proof of the proposition.
\end{Proof}

We now turn our attention to the crucial case $\frac{|x-y|}{t}\leq K$, where we decompose the Green's function into a part that is well behaved in the weighted space and one induced by the inhomogenous coupling with weak spatial but strong temporal decay.

\begin{figure}[!t]
\centering
\includegraphics[width=0.35\textwidth]{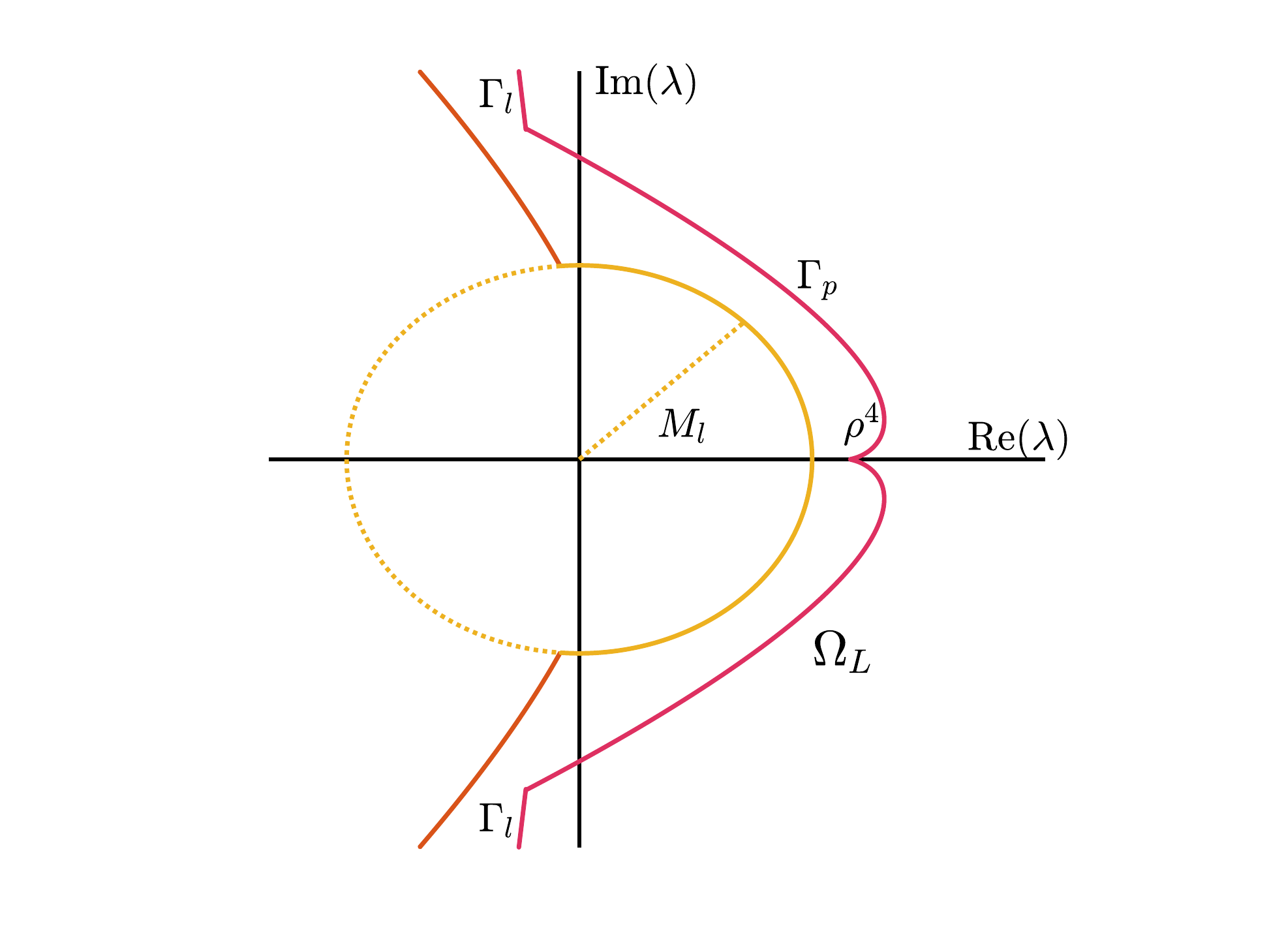}\hspace*{0.8in}\includegraphics[width=0.35\textwidth]{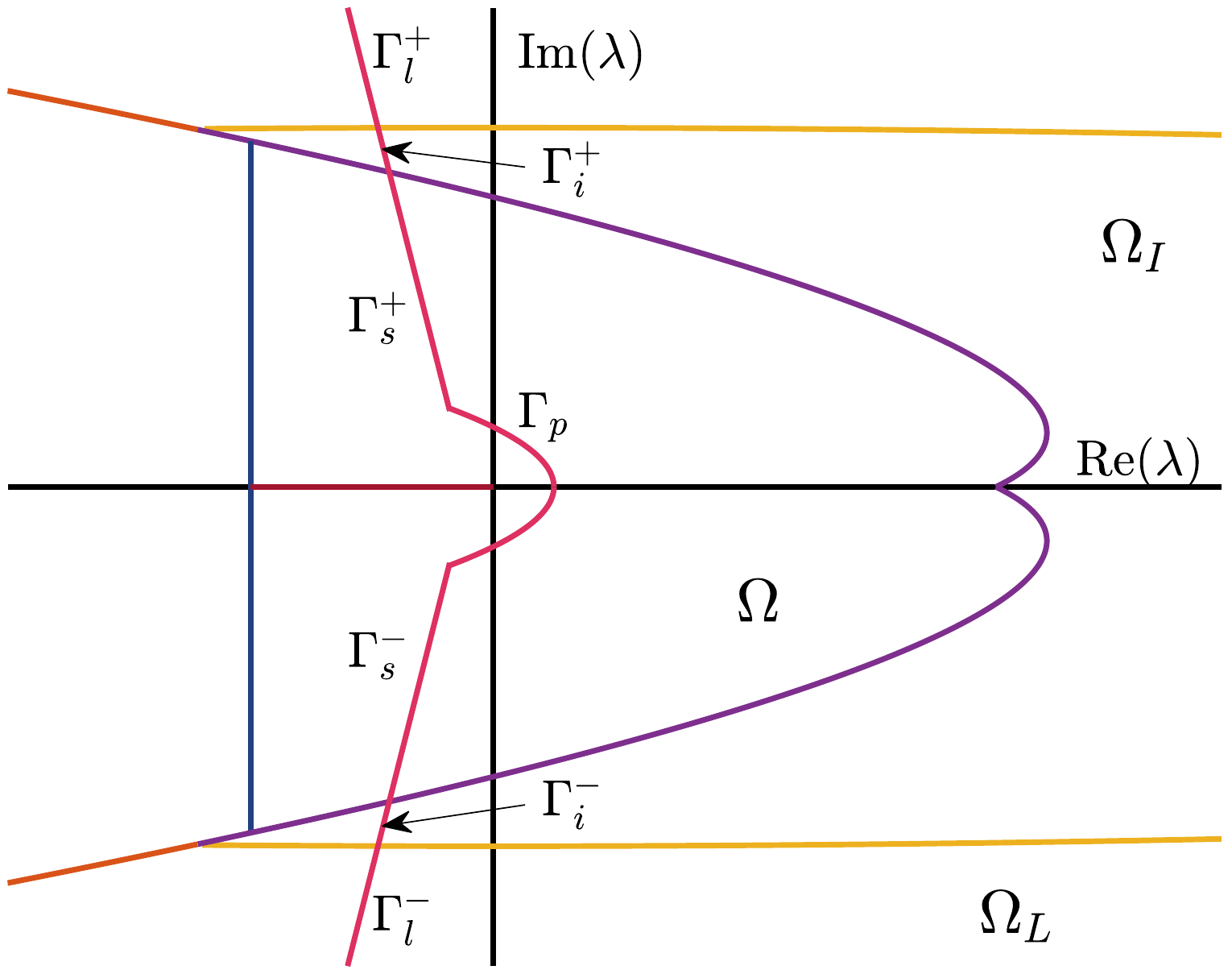}

\caption{Left: The contours $\Gamma_p$ and $\Gamma_l$ used in the proof of Proposition~\ref{timeGreen12} in the regime $|x-y|\geq Kt$. Right: Contours used in the proof of Proposition~\ref{prop:timegreen12} in the regime $|x-y|\leq Kt$ and $t\geq1$.}
\label{fig:LargeRegion}\label{fig:SmallRegion}
\end{figure}

\begin{prop}\label{prop:timegreen12} Let $K$ be as in Proposition~\ref{timeGreen12}.  There exists constants $C,\,\kappa,\, \theta>0$ such that for $t>1$ and $|x-y|\leq Kt$, the temporal Green's function $\cG^{12}(t,x,y)$ has the  decomposition
\bqs
\cG^{12}(t,x,y)=\cG^{12}_\omega(t,x,y)+\widetilde{\cG^{12}}(t,x,y)
\eqs
with estimates
\begin{align*}
\left|\cG^{12}_\omega(t,x,y)\right| &\leq C \frac{\omega(x)}{\omega(y)}\left(\frac{1+|x-y|}{t^{3/2}}\me^{-\frac{|x-y|^2}{\kappa t}}+\me^{-\theta t}  \right),\\
\left| \widetilde{\cG^{12}}(t,x,y) \right| &\leq   \left\{ \begin{array}{cc} C \me^{-\theta t} \me^{-\gamma_v (x-y)}, & x>y, \\
 C \me^{-\theta t} \me^{\delta_v (x-y)}, & x<y. \end{array}\right.
\end{align*}

\end{prop}

\begin{Proof} Recall that $\cG^{12}(t,x,y)$ is obtained via the following inverse Laplace transform formula
\bqq
\cG^{12}(t,x,y)=\frac{1}{2\pi \mbi}\int_\Gamma \me^{\lambda t}\bG_\lambda^{12}(x,y)\md \lambda,
\label{eq:LaplaceG12temp}
\eqq
for some well-chosen contour $\Gamma$ in the complex plane that may depend on the relative position between $x$ and $y$; see  Figure~\ref{fig:SmallRegion} for a typical contour. Our contour integral will be decomposed using
\begin{itemize}[itemsep=-0pt,topsep=-0pt]
\item{\it large} $\lambda$ estimates: this corresponds to the part of $\lambda\in\Gamma$ for which $\lambda\in\Omega_L$ where we will use our large $\lambda$ estimate from Lemma~\ref{lem:largelambda};
\item {\it intermediate} $\lambda$ estimates: this corresponds to the part of $\lambda\in\Gamma$ for which $\lambda\in\Omega_I$ where we will use our intermediate $\lambda$ estimate from Lemma~\ref{lem:lambdaint};
\item {\it small} $\lambda$ estimates: this corresponds to the part of $\lambda\in\Gamma$ for which $\lambda \in \Omega$ where we will use our small $\lambda$ estimate from Lemma~\ref{lem:G12full}.
\end{itemize}

We first decompose $\Gamma=\Gamma_{\text{out}}\cup \Gamma_\Omega$, with $\Gamma\subset\Omega$ and $\Gamma_{\text{out}} \subset \Omega_I\cup\Omega_L$. For constants $\delta_{0,1}>0$, $\ell_{l,i}>0$, specified later, we decompose further $\Gamma_{\text{out}}=\Gamma_l\cup \Gamma_i$ with
\bqs
\Gamma_l\coloneqq\left\{\lambda(\ell)=-\delta_0-\delta_1|\ell|+\mbi\ell ~|~ |\ell|>\ell_l \right\} \ \text{ and } \ \Gamma_i\coloneqq\left\{\lambda(\ell)=-\delta_0-\delta_1|\ell|+\mbi\ell ~|~ \ell_i<|\ell|<\ell_l \right\}.
\eqs
We also define $\Gamma_{i,l}\coloneqq\Gamma_{i,l}^+\cup\Gamma_{i,l}^-$ with $\Gamma_{i,j}^\pm$ located in $\Im(\lambda)>0$ and  $\Im(\lambda)<0$, respectively. The contour $\Gamma_\Omega$ which joins the two points $\lambda(\ell_i)$ and $\overline{\lambda(\ell_i)}$ will be chosen later in the proof.

Now decomposing \eqref{eq:LaplaceG12temp}, 
\bqs
\cG^{12}(t,x,y)=\frac{1}{2\pi \mbi}\int_{\Gamma_\Omega} \me^{\lambda t}\bG_\lambda^{12}(x,y)\md \lambda+\frac{1}{2\pi \mbi}\int_{\Gamma_{\text{out}}}\me^{\lambda t}\bG_\lambda^{12}(x,y)\md \lambda.
\eqs
we can easily estimate the second integral using Lemma~\ref{lem:largelambda} and Lemma~\ref{lem:lambdaint} as follows.

\textbf{Integral along $\Gamma_l$.} Along $\Gamma_l$, we use the large $\lambda$ estimate
\bqs
\left| \bG_\lambda^{12}(x,y)\right|\leq \frac{C}{|\lambda|^{7/4}}\frac{\omega(x)}{\omega(y)}\me^{-|\lambda|^{1/4} \eta |x-y|}, \quad \forall(x,y)\in\R^2.
\eqs
and obtain, in a similar fashion as in the proof of Proposition~\ref{timeGreen12}, 
\bqs
\left| \frac{1}{2\pi \mbi}\int_{\Gamma_l}\me^{\lambda t}\bG_\lambda^{12}(x,y)\md \lambda \right| \leq C_0 \frac{\omega(x)}{\omega(y)} \me^{-M_l^{1/4}\eta |x-y|}\me^{-\delta_0 t} \int_{\ell_l}^\infty \me^{-\delta_1 \ell t} \md \ell,
\eqs
where we used the fact that for $\lambda(\ell)\in \Omega_L$ we have that
\bqs
\frac{1}{|\lambda(\ell)|^{7/4}}\leq C_1,
\eqs
for some constant $C_1>0$. As a consequence, we obtain an estimate of the form
\bqs
\left| \frac{1}{2\pi \mbi}\int_{\Gamma_l}\me^{\lambda t}\bG_\lambda^{12}(x,y)\md \lambda \right| \lesssim \frac{\omega(x)}{\omega(y)} \me^{-\tilde{\eta} |x-y|}\me^{-\theta t},
\eqs
with positive constants $\tilde{\eta}>0$ and $\theta>0$.

\textbf{Integral along $\Gamma_i$.} Along $\Gamma_i$ we use the intermediate $\lambda$ estimate 
\bqs
\left| \bG_\lambda^{12}(x,y)\right|\leq C\frac{\omega(x)}{\omega(y)},
\eqs
to find
\bqs
\left| \frac{1}{2\pi \mbi}\int_{\Gamma_i}\me^{\lambda t}\bG_\lambda^{12}(x,y)\md \lambda \right| \leq C_0 \frac{\omega(x)}{\omega(y)} \me^{-\delta_0 t} \int_{\ell_i}^{\ell_l} \me^{-\delta_1 \ell t} \md\ell \lesssim \frac{\omega(x)}{\omega(y)} \me^{-\theta t}.
\eqs
Summarizing, we have obtained the estimate 
\bqs
\left|\frac{1}{2\pi \mbi}\int_{\Gamma_{\text{out}}}\me^{\lambda t}\bG_\lambda^{12}(x,y)\md \lambda\right| \lesssim \frac{\omega(x)}{\omega(y)} \me^{-\theta t}.
\eqs

\textbf{Integral along $\Gamma_\Omega$.} In this case, we distinguish several cases depending on the location of $x$ and $y$, detailing only one and outlining the differences for the other cases. 

{\bf Case 1: $y<0<x$.}  From Lemma~\ref{lem:G12full}, we have the decomposition of $\bG_\lambda^{12}(x,y)$, for $\lambda\in\Omega$,
\begin{align*}
\bG_\lambda^{12}(x,y) & = \frac{1}{\sqrt{\zeta(\lambda-\lambda_v^{\textnormal{bp}})}}\me^{\nu_-^0(\lambda)x-\nu_-^1(\lambda)y} \boldH_\lambda(x,y) \\
&~~~+\sum_{j=1}^2 \frac{1}{(\lambda-\lambda^{\textnormal{rp}})\sqrt{\zeta(\lambda-\lambda_v^{\textnormal{bp}})}} \left[\me^{\nu_-^0(\lambda)x-\nu_j(\lambda)y}\bI_{\lambda,j}(x,y)+\me^{\nu_j(\lambda)(x-y)}\boldJ_{\lambda,j}(x,y) \right].
\end{align*}
Within the region $\Omega$, we use different contours when handling the first and second terms in the expansion involving $\boldH_\lambda(x,y)$ and $\bI_{\lambda,j}(x,y)$ versus  handling the third term involving $\boldJ_{\lambda,j}(x,y)$. We first treat the integral involving $\boldH_\lambda(x,y)$. We decompose $\Gamma_\Omega$ into 
\begin{itemize}[itemsep=0pt,topsep=0pt]
\item a parabolic contour $\Gamma_p$ near the origin depending on $\rho>0$ to be specified later,
\bqs
\Gamma_p=\left\{ \lambda\in\C~|~ \sqrt{\frac{\lambda(k)}{d}} = \rho +\mbi k \text{ for } -k_p\leq k \leq k_p \right\}\subset \Omega.
\eqs
\item two straight, symmetric contours $\Gamma_s^\pm$, continuations of the rays coming from $\Gamma_{\text{out}}$, 
\bqs
\Gamma_s\coloneqq\Gamma_s^+\cup\Gamma_s^-=\left\{\lambda(\ell)=-\delta_0-\delta_1|\ell|+\mbi \ell~|~ \ell_p<|\ell|<\ell_i\right\}.
\eqs
\end{itemize}
The constants $k_p>0$ and $\ell_p>0$ are defined according to 
\bqs
 k_p = \delta_1 \rho + \sqrt{(\delta_1\rho)^2+\rho^2+\frac{\delta_0}{d}} \text{ and }\ell_p = 2d\rho k_p,
\eqs
such that
\bqs
d(\rho^2-k_p^2)+2d\rho k_p\mbi = -\delta_0-\delta_1 \ell_p+\mbi \ell_p.
\eqs

Along the parabolic contour, we note that
\bqs
\frac{1}{\sqrt{\zeta(\lambda-\lambda_v^{\textnormal{bp}})}}\me^{\nu_-^0(\lambda)x-\nu_-^1(\lambda)y} \boldH_\lambda(x,y)=\me^{-\frac{s_*}{2d}(x-y)-\sqrt{\frac{\lambda}{d}}(x-y)}\me^{\left(\frac{\sqrt{s_*^2-4d(f'(1)-\lambda)}}{2d}-\sqrt{\frac{\lambda}{d}}\right)y}\frac{\boldH_\lambda(x,y)}{\sqrt{\zeta(\lambda-\lambda_v^{\textnormal{bp}})}},
\eqs
and we select 
\bqs
\rho=\frac{|x-y|}{Lt},
\eqs
with $L>d$ sufficiently large so that the parabolic contour lies within $\Omega$ and the region for which 
\bqs
\Re\left(\frac{\sqrt{s_*^2-4d(f'(1)-\lambda)}}{2d}-\sqrt{\frac{\lambda}{d}}\right)>0.
\eqs
As a consequence, along $\Gamma_p$, the function
\bqs
\widetilde{\boldH_\lambda}(x,y)\coloneqq \me^{\left(\frac{\sqrt{s_*^2-4d(f'(1)-\lambda)}}{2d}-\sqrt{\frac{\lambda}{d}}\right)y}\frac{\boldH_\lambda(x,y)}{\sqrt{\zeta(\lambda-\lambda_v^{\textnormal{bp}})}},
\eqs
is uniformly bounded in $(x,y)$ and analytic in $\sqrt{\lambda}$, and we have the identity 
\bqs
\frac{1}{2\pi\mbi}\int_{\Gamma_p}\me^{\lambda t}  \frac{1}{\sqrt{\zeta(\lambda-\lambda_v^{\textnormal{bp}})}}\me^{\nu_-^0(\lambda)x-\nu_-^1(\lambda)y} \boldH_\lambda(x,y)\md \lambda = \frac{1}{2\pi\mbi} \me^{-\frac{s_*}{2d}(x-y) }\int_{\Gamma_p}\me^{\lambda t} \me^{-\sqrt{\frac{\lambda}{d}}(x-y)}  \widetilde{\boldH_\lambda}(x,y) \md \lambda.
\eqs
From there, we proceed along similar lines as in \cite{faye19}. We have that
\bqs
\md \lambda(k) = 2d\mbi \left(\rho+\mbi k \right)\md k, \quad \frac{\lambda(k)}{d} = \rho^2-k^2+2\rho\mbi k,
\eqs
such that
\bqs
\int_{\Gamma_p}\me^{\lambda t} \me^{-\sqrt{\frac{\lambda}{d}}(x-y)}  \widetilde{\boldH_\lambda}(x,y) \md \lambda=2d \mbi \me^{d\rho^2t-\rho (x-y)}\int_{-k_p}^{k_p}\me^{-dk^2t}\me^{2d \rho \mbi kt - \mbi k(x-y)} \widetilde{\boldH_{\lambda(k)}}(x,y) (\rho+\mbi k)\md k.
\eqs
Next, we separate $\me^{2d \rho \mbi kt - \mbi k(x-y)} \widetilde{\boldH_{\lambda(k)}}(x,y)$ into its real and imaginary parts as
\bqs
\me^{2d \rho \mbi kt - \mbi k(x-y)} \widetilde{\boldH_{\lambda(k)}}(x,y)\coloneqq\boldH_r(x,y,k)+\mbi \boldH_i(x,y,k),
\eqs
and note that $\boldH_r(x,y,k)$ is even in $k$ while $\boldH_i(x,y,k)$ is odd in $k$. As a consequence, we obtain
\bqs
\int_{\Gamma_p}\me^{\lambda t} \me^{-\sqrt{\frac{\lambda}{d}}(x-y)}  \widetilde{\boldH_\lambda}(x,y) \md \lambda=2d\mbi \int_{-k_p}^{k_p}\me^{-dk^2t} \left(\rho\boldH_r(x,y,k)-k\boldH_i(x,y,k)  \right) \md k.
\eqs
Boundedness of $\boldH_r(x,y,k)$ implies that
\bqs
\left| \int_{-k_p}^{k_p}\me^{-dk^2t}\boldH_r(x,y,k)\md k\right| \leq \frac{C}{\sqrt{t}},
\eqs
while oddness of $\boldH_i(x,y,k)$ leads to
\bqs
\left| \int_{-k_p}^{k_p}\me^{-dk^2t}k\boldH_i(x,y,k)\md k\right| \leq \frac{C}{t^{3/2}}.
\eqs
Combining all the above estimates and recalling that $\rho=\frac{|x-y|}{Lt}$, we arrive at
\bqs
\left|\int_{\Gamma_p}\me^{\lambda t} \me^{-\sqrt{\frac{\lambda}{d}}(x-y)}  \widetilde{\boldH_\lambda}(x,y) \md \lambda\right| \lesssim \frac{1+|x-y|}{t^{3/2}}\me^{-\frac{|x-y|^2}{\kappa t}},
\eqs
with $\kappa=\max\left(\frac{L^2}{L-d},\frac{L^2}{2d\delta_1^2}\right)>0$. To complete this part of the analysis, we need to obtain estimates along the remaining contour $\Gamma_s$. For this, we note that 
\bqs
\left|\int_{\Gamma_s}\me^{\lambda t} \me^{-\sqrt{\frac{\lambda}{d}}(x-y)}  \widetilde{\boldH_\lambda}(x,y) \md \lambda \right| \leq C \me^{-\delta_0 t} \int_{\ell_p}^{\ell_i} \me^{-\delta_1 \ell t} \me^{-\Re\left(\sqrt{\frac{\lambda(\ell)}{d}} \right) (x-y)} \md \ell \leq C \me^{-\delta_0 t} \frac{\me^{-\delta_1 \ell_p t}}{t}.
\eqs
Next, note that $\ell_p=2d\rho k_p>2d\delta_1\rho^2$ such that
\bqs
\left|\int_{\Gamma_s}\me^{\lambda t} \me^{-\sqrt{\frac{\lambda}{d}}(x-y)}  \widetilde{\boldH_\lambda}(x,y) \md \lambda \right| \leq \frac{C}{t^{3/2}} \me^{-\frac{|x-y|^2}{\kappa t}}.
\eqs
In summary, we have have shown that
\bqs
\left| \frac{1}{2\pi\mbi}\int_{\Gamma_\Omega} \me^{\lambda t}  \frac{1}{\sqrt{\zeta(\lambda-\lambda_v^{\textnormal{bp}})}}\me^{\nu_-^0(\lambda)x-\nu_-^1(\lambda)y} \boldH_\lambda(x,y) \md\lambda\right| \lesssim \frac{1+|x-y|}{t^{3/2}} \me^{-\frac{s_*}{2d}(x-y)} \me^{-\frac{|x-y|^2}{\kappa t}}.
\eqs
We next address the integral term involving $\bI_{\lambda,j}(x,y)$ using the same contour $\Gamma_\Omega=\Gamma_s\cup\Gamma_p$. Since $\Re(\lambda^{\textnormal{rp}})<0$, we can choose $\delta_{0,1}$ small enough and $L$ large such that the contour $\Gamma_\Omega$ is located to the right of the resonance pole $\lambda^{\textnormal{rp}}$ and its complex conjugate. Along the parabolic contour,
\bqs
\me^{\nu_-^0(\lambda)x-\nu_j(\lambda)y} =\me^{-\frac{s_*}{2d}x+\gamma_v y-\sqrt{\frac{\lambda}{d}}(x-y)}\me^{\left(-\gamma_v-\nu_j(\lambda)-\sqrt{\frac{\lambda}{d}}\right)y}, \quad j=1,2,
\eqs
and we choose  
$\rho=|x-y|/(Lt)
$
with $L>d$ sufficiently large such that $\Gamma_p\subset\Omega$ and 
\bqs
\Gamma_p\subset\left\{\lambda\left|\Re\left(-\gamma_v-\nu_j(\lambda)-\sqrt{\lambda/d} \right)\right.>0, \quad j=1,2\right\},
\eqs
which is always possible with our definition of $\gamma_v>0$ in \eqref{eq:gammav}. As a consequence, 
along $\Gamma_p$, %the following function
\bqs
\widetilde{\bI_{\lambda,j}}(x,y)\coloneqq \me^{\left(-\gamma_v-\nu_j(\lambda)-\sqrt{\frac{\lambda}{d}}\right)y}\frac{\bI_{\lambda,j}(x,y)}{(\lambda-\lambda^{\textnormal{rp}})\sqrt{\zeta(\lambda-\lambda_v^{\textnormal{bp}})}}, \quad j=1,2,
\eqs
is uniformly bounded in $(x,y)$ and analytic in $\sqrt{\lambda}$. An analysis analogous to the case of $\boldH_\lambda(x,y)$, gives
\bqs
\left| \frac{1}{2\pi\mbi}\int_{\Gamma_\Omega} \me^{\lambda t} \sum_{j=1}^2  \frac{1}{(\lambda-\lambda^{\textnormal{rp}})\sqrt{\zeta(\lambda-\lambda_v^{\textnormal{bp}})}} \me^{\nu_-^0(\lambda)x-\nu_j(\lambda)y}\bI_{\lambda,j}(x,y)\md\lambda\right| \lesssim \frac{1+|x-y|}{t^{3/2}} \me^{-\frac{s_*}{2d}x+\gamma_vy} \me^{-\frac{|x-y|^2}{\kappa t}}.
\eqs
To complete this case, we address the term involving $\bJ_{\lambda,j}(x,y)$. Since there is no singularity at the origin and we may choose the contour 
\bqs
\Gamma_\Omega=\left\{\lambda(\ell)=-\delta_0-\delta_1|\ell|+\mbi \ell~|~ 0\leq |\ell|<\ell_i\right\},
\eqs 
with $\delta_{0,1}>0$ small such that $\Gamma_\Omega$ is to the right of potential resonance poles. Next, notice that
\bqs
\left| \frac{1}{(\lambda-\lambda^{\textnormal{rp}})\sqrt{\zeta(\lambda-\lambda_v^{\textnormal{bp}})}} \me^{\nu_j(\lambda)(x-y)}\bJ_{\lambda,j}(x,y) \right|\leq C \me^{-\gamma_v(x-y)}, \quad j=1,2,
\eqs
for all $\lambda \in \Gamma_\Omega$ uniformly in $(x,y)$. Altogether, this yields
\bqs
\left| \frac{1}{2\pi\mbi}\int_{\Gamma_\Omega} \me^{\lambda t} \sum_{j=1}^2  \frac{1}{(\lambda-\lambda^{\textnormal{rp}})\sqrt{\zeta(\lambda-\lambda_v^{\textnormal{bp}})}} \me^{\nu_j(\lambda)(x-y)}\bJ_{\lambda,j}(x,y)\md\lambda\right| \lesssim \me^{-\theta t} \me^{-\gamma_v(x-y)}.
\eqs
Returning to the definition of $\cG^{12}(t,x,y)$ we  established that the decomposition
\bqs
\cG^{12}(t,x,y)=\cG^{12}_\omega(t,x,y)+\widetilde{\cG^{12}}(t,x,y), \qquad \cG^{12}_\omega(t,x,y)\coloneqq\cG^{12}(t,x,y)-\widetilde{\cG^{12}}(t,x,y),
\eqs
with
\bqs
\widetilde{\cG^{12}}(t,x,y)\coloneqq \frac{1}{2\pi\mbi}\int_{\Gamma_\Omega} \me^{\lambda t} \sum_{j=1}^2  \frac{1}{(\lambda-\lambda^{\textnormal{rp}})\sqrt{\zeta(\lambda-\lambda_v^{\textnormal{bp}})}} \me^{\nu_j(\lambda)(x-y)}\bJ_{\lambda,j}(x,y)\md\lambda,
\eqs
gives estimates 
\bqs
\left| \cG^{12}_\omega(t,x,y)\right| \leq C \me^{-\frac{s_*}{2d}x-\delta y} \left( \frac{1+|x-y|}{t^{3/2}}  \me^{-\frac{|x-y|^2}{\kappa t}} +\me^{-\theta t} \right),
\text{ and }\left|\widetilde{\cG^{12}}(t,x,y)\right| \leq C \me^{-\theta t} \me^{-\gamma_v(x-y)}.
\eqs

The remaining cases can be handled in a similar fashion.

{\bf Case 2: $0<y<x$.} Estimates for the term involving $\boldH_\lambda(x,y)$ are obtained as in  Case 1.  Note  again that the resonance poles occurring in this term are of no concern since all contour integrations can take place to the right of them. Treating the second term involving $\bJ_{\lambda,j}(x,y)$, we follow the reasoning in above. Since $y>0$, we find
\bqs
\left| \cG^{12}_\omega(t,x,y)\right| \leq C \me^{-\frac{s_*}{2d}(x-y)} \left( \frac{1+|x-y|}{t^{3/2}}  \me^{-\frac{|x-y|^2}{\kappa t}} +\me^{-\theta t} \right),\text{ and }
\left|\widetilde{\cG^{12}}(t,x,y)\right| \leq C \me^{-\theta t} \me^{-\gamma_v(x-y)}.
\eqs

{\bf Case 3: $0<x<y$.}  The terms involving $\boldH_\lambda(x,y)$ and $\bI_{\lambda,j}(x,y)$ are analogous to Case 2.  For the term involving $\boldJ_{\lambda,j}(x,y)$ note the lack of a resonance pole.  Therefore, the analysis is similar to that of Lemma~\ref{timeGreen22} for the case $x<y$.  We thus obtain
\bqs
\left| \cG^{12}_\omega(t,x,y)\right| \leq C \me^{-\frac{s_*}{2d}(x-y)} \left( \frac{1+|x-y|}{t^{3/2}}  \me^{-\frac{|x-y|^2}{\kappa t}} +\me^{-\theta t} \right),
\eqs
together with
\bqs
\left|\widetilde{\cG^{12}}(t,x,y)\right| \leq C \me^{-\theta t} \me^{\delta_v(x-y)}.
\eqs
We now turn our attention to the cases where $x<0$. 

{\bf Case 4: $x<y<0$.}  Estimates for the term involving $\boldH_\lambda(x,y)$ can be obtained by writing 
\bqs
\me^{\nu_+^1(\lambda)(x-y)}=\me^{\delta(x-y)}\me^{\sqrt{\frac{\lambda}{d}}(x-y)}\me^{\left(\nu_+^1(\lambda)-\delta-\sqrt{\frac{\lambda}{d}}\right)(x-y)}
\eqs
such that along a parabolic contour near the origin one can ensure that
\bqs
\Re\left(\nu_+^1(\lambda)-\delta-\sqrt{\frac{\lambda}{d}} \right)>0,
\eqs
which is always possible since $\delta>0$ can be chosen arbitrarily small. One finds
\bqs
\left| \frac{1}{2\pi\mbi}\int_{\Gamma_\Omega} \me^{\lambda t}  \frac{1}{\sqrt{\zeta(\lambda-\lambda_v^{\textnormal{bp}})}}\me^{\nu_+^1(\lambda)(x-y)} \boldH_\lambda(x,y) \md\lambda\right| \lesssim \frac{1+|x-y|}{t^{3/2}} \me^{\delta(x-y)} \me^{-\frac{|x-y|^2}{\kappa t}}.
\eqs
The contribution from $\bI_{\lambda,j}(x,y)$ is handled similarly and we obtain
\bqs
\left| \frac{1}{2\pi\mbi}\int_{\Gamma_\Omega} \me^{\lambda t} \sum_{j=1}^2  \frac{1}{(\lambda-\lambda^{\textnormal{rp}})\sqrt{\zeta(\lambda-\lambda_v^{\textnormal{bp}})}}\me^{\nu_+^1(\lambda)x-\nu_j(\lambda)y} \bI_{\lambda,j}(x,y) \md\lambda\right| \lesssim \frac{1+|x-y|}{t^{3/2}} \me^{\delta x + \gamma_v y} \me^{-\frac{|x-y|^2}{\kappa t}}.
\eqs
Finally, the last contribution involving $\boldJ_{\lambda,j}(x,y)$ is treated as in the proof of Lemma~\ref{timeGreen22} to obtain
\bqs
\left| \frac{1}{2\pi\mbi}\int_{\Gamma_\Omega} \me^{\lambda t} \sum_{j=3}^4  \frac{1}{\sqrt{\zeta(\lambda-\lambda_v^{\textnormal{bp}})}}\me^{\nu_j(\lambda)(x-y)} \bJ_{\lambda,j}(x,y) \md\lambda\right| \lesssim \me^{-\theta t} \me^{\delta_v(x-y)}.
\eqs
{\bf Case 5: $y<x<0$.} This case is almost identical to Case 4. We only note that the contribution from $\boldJ_{\lambda,j}(x,y)$ gives
\bqs
\left| \frac{1}{2\pi\mbi}\int_{\Gamma_\Omega} \me^{\lambda t} \sum_{j=1}^2  \frac{1}{\sqrt{\zeta(\lambda-\lambda_v^{\textnormal{bp}})}}\me^{\nu_j(\lambda)(x-y)} \bJ_{\lambda,j}(x,y) \md\lambda\right| \lesssim \me^{-\theta t} \me^{-\gamma_v(x-y)}.
\eqs
{\bf Case 6: $x<0<y$.} We only comment on the contribution from $\boldH_{\lambda}(x,y)$ as the terms involving $\bI_{\lambda,j}(x,y)$ and $\boldJ_{\lambda,j}(x,y)$ do not present new difficulties. For $\boldH_{\lambda}(x,y)$, we note that
\bqs
\me^{\nu_+^1(\lambda)x-\nu_+^0(\lambda)y}=\me^{\delta x +\frac{s_*}{2d}y}\me^{\sqrt{\frac{\lambda}{d}}(x-y)}\me^{\left(\nu_+^1(\lambda)-\delta -\sqrt{\frac{\lambda}{d}}\right)x},
\eqs
and that, along a parabolic contour near the origin, one can always ensure that
\bqs
\Re\left(\nu_+^1(\lambda)-\delta -\sqrt{\frac{\lambda}{d}} \right)>0,
\eqs
which is always possible since $\delta>0$ is arbitrarily small. Concluding the proof, one finds
\bqs
\left| \frac{1}{2\pi\mbi}\int_{\Gamma_\Omega} \me^{\lambda t}  \frac{1}{(\lambda-\lambda^{\textnormal{rp}})\sqrt{\zeta(\lambda-\lambda_v^{\textnormal{bp}})}}\me^{\nu_+^1(\lambda)x-\nu_+^0(\lambda)y} \boldH_\lambda(x,y) \md\lambda\right| \lesssim \frac{1+|x-y|}{t^{3/2}} \me^{\delta x +\frac{s_*}{2d}y} \me^{-\frac{|x-y|^2}{\kappa t}}.
\eqs
\end{Proof}

\subsection{Nonlinear stability}\label{sec:nonlinear}
We are now ready to state precisely and prove our main result on nonlinear stability, Theorem~\ref{thm:maininformal}. Consider the smooth and bounded weight function $w_v>0$ 
\[ w_v(x) \coloneqq \left\{ \begin{array}{cc} \me^{-\gamma_v x}, & x>1, \\
1, & x=0,\\ \me^{\delta_v x}, & x<-1, \end{array}\right. 
\]
and recall the definition of $\omega$ from \eqref{eqomega} with $\delta>0$ small enough such that $0<\delta<\min\left(\nu_+^1(0),2\delta_v\right)$.

\begin{thm}[Main result]\label{thmMain} For $(d,\alpha,\mu)\in\Pi$ as in Definition~\ref{defiPi} and that $\beta\neq0$,  consider \eqref{eq:main} with initial condition $(u(0,\cdot),v(0,\cdot))=(Q_*(\cdot),0)+(P_0(\cdot),v_0(\cdot))$. There exist $C,\,\epsilon,\,\theta>0$ such that, if $(P_0,v_0)$ satisfies 
\begin{equation}\label{e:ismall}
\left\|\frac{P_0}{\omega}\right\|_\infty+\left\|\frac{v_0}{\omega}\right\|_\infty+\int_\R (1+|y|)\left(\frac{|v_0(y)|}{\omega(y)}+\frac{|v_0(y)|}{\omega(y)}\right)\md y+\left\| \frac{v_0}{\omega_v} \right\|_\infty+\left\| \frac{v_0}{\omega_v} \right\|_1 <\epsilon,
\end{equation}
then the solution $(u(t,x),v(t,x))$ is defined for all time and the critical pulled front is nonlinearly stable in the sense that the perturbation $P(t,x)\coloneqq u(t,x)-Q_*(x)$ admits the  decomposition
\bqs
P(t,x)=E(t,x)+\omega(x)p(t,x), \quad x\in\R,
\eqs
such that for all $t\geq 1$, 
\begin{align*}
\underset{x\in\R}{\sup}~\frac{|E(t,x)|}{\omega_v(x)} &\leq Ce^{-\theta t} \left(\left\| \frac{v_0}{\omega_v} \right\|_\infty+\left\| \frac{v_0}{\omega_v} \right\|_1 \right),\\
\underset{x\in\R}{\sup}~\frac{|p(t,x)|}{1+|x|}&\leq \frac{C\epsilon}{(1+t)^{3/2}},\\
\underset{x\in\R}{\sup}~\frac{|v(t,x)|}{\omega_v(x)}&\leq C\me^{-\theta t}\left\| \frac{v_0}{\omega_v} \right\|_\infty.
\end{align*}
On the other hand, there exists $\epsilon_0>0$  such that for any $\epsilon>0$ and initial conditions $(P_0,v_0)$ satisfying \eqref{e:ismall} with this choice of  $\epsilon$ arbitrarily small we have for some $t>0$, 
\[
  \underset{x\in\R}{\sup}~\frac{|P(t,x)|}{\omega(x)(1+|x|)} \geq \varepsilon,
\]
that is, the decay in $p$ does not hold for $E/\omega$. 
\end{thm}

The remainder of this section consists of the proof of this theorem. We introduce the decomposition in Section~\ref{s:decomp}, and carry out the nonlinear iteration scheme to prove the nonlinear estimates on $E$ and $p$ in Section~\ref{s:nonlp}. We conclude in Section~\ref{s:inst} demonstrating that stability in a fixed weight cannot be achieved. 

\subsubsection{Identification of $E(t,x)$ and $h(t,x)$}\label{s:decomp}
Following the outline of the proof in Section~\ref{sec:argoutline}, we wish to separate terms whose spatial decay is too weak to  yield temporal decay in the exponentially weighted space required for the stability proof of the front, but who nevertheless decay  exponentially in time. For $t>0$ and $x\in\R$, we set
\begin{equation}\label{e:defH}
H(t,x)\coloneqq\int_\R \cG^{12}(t,x,y)v_0(y)\md y.
\end{equation}
We first note that for $0<t\leq1$, Proposition~\ref{timeGreen12} shows that
\bqs
\left|\int_\R \cG^{12}(t,x,y)v_0(y)\md y\right| \leq C \omega(x) \left\|\frac{v_0}{\omega}\right\|_\infty.
\eqs
We then therefore find short-time bounds on $h(t,x)\coloneqq{H(t,x)}/{\omega(x)}, \quad x \in\R. $
\begin{lem}\label{lem:shorttimeh} There exists $C>0$ such that for all $0<t\leq1$ we have
\bqs
\|h(t,\cdot)\|_\infty\leq C \left\|\frac{v_0}{\omega}\right\|_\infty, \quad x \in\R.
\eqs
\end{lem}
For $t>1$, we decompose
\begin{equation}
H(t,x)=\int_{-\infty}^{x-Kt} \cG^{12}(t,x,y)v_0(y)\md y+\int_{x-Kt}^{x+Kt} \cG^{12}(t,x,y)v_0(y)\md y+\int_{x+Kt}^{+\infty} \cG^{12}(t,x,y)v_0(y)\md y,
\end{equation}
and further split the second integral using Proposition~\ref{prop:timegreen12},
\bqs
\int_{x-Kt}^{x+Kt} \cG^{12}(t,x,y)v_0(y)\md y=\int_{x-Kt}^{x+Kt} \cG^{12}_\omega(t,x,y)v_0(y)\md y+\int_{x-Kt}^{x+Kt} \widetilde{\cG^{12}}(t,x,y)v_0(y)\md y.
\eqs
The next result gives bounds on the second term  $E(t,x)$ in this splitting, 
\bqs
E(t,x)\coloneqq\int_{x-Kt}^{x+Kt} \widetilde{\cG^{12}}(t,x,y)v_0(y)\md y, \quad x\in\R.
\eqs
\begin{lem}\label{lemEstimE} There exists $C>0$ such that for all $t\geq 1$ and all $x\in\R$ we have
\bqq |E(t,x)|\leq C \me^{-\theta t}\omega_v(x) \left( \left\| \frac{v_0}{\omega_v} \right\|_\infty+\left\| \frac{v_0}{\omega_v} \right\|_1\right).
\label{eqEstimE}
\eqq
\end{lem}

\begin{Proof}
We use the estimate on $\widetilde{\cG^{12}}(t,x,y)$ from Proposition~\ref{prop:timegreen12}, valid for $t\geq1$ and $|x-y|\leq Kt$
\begin{align*}
\left|\int_{x-Kt}^{x+Kt} \widetilde{\cG^{12}}(t,x,y)v_0(y)\md y\right| &\leq  \int_{x-Kt}^{x} \left|\widetilde{\cG^{12}}(t,x,y)\right| \left|v_0(y)\right|\md y +\int_{x}^{x+Kt} \left|\widetilde{\cG^{12}}(t,x,y)\right| \left|v_0(y)\right|\md y\\
&\leq C \me^{-\theta t} \left( \me^{-\gamma_v x} \int_{x-Kt}^{x} \me^{\gamma_v y} \left|v_0(y)\right|\md y + \me^{\delta_v x} \int_{x}^{x+Kt} \me^{-\delta_v y}\left|v_0(y)\right|\md y \right).
\end{align*}
For $x\geq0$, the second integral is bounded by
\bqs
\int_{x}^{x+Kt} \me^{-\delta_v y}\left|v_0(y)\right|\md y \leq \int_{x}^{+\infty} \me^{-\delta_v y}\left|v_0(y)\right|\md y\leq  C \me^{-(\delta_v+\gamma_v)x} \left\| \frac{v_0}{\omega_v} \right\|_\infty.
\eqs
On the other hand, the first integral can be evaluated as
\begin{align*}
 \int_{x-Kt}^{x} \me^{\gamma_v y} \left|v_0(y)\right|\md y &\leq \int_{-\infty}^0 \me^{(\gamma_v+\delta_v)y} \me^{-\delta_v y} \left|v_0(y)\right|\md y + \int_0^\infty \me^{\gamma_v y} \left|v_0(y)\right|\md y \\
&\leq C\left(\left\| \frac{v_0}{\omega_v} \right\|_\infty+\left\| \frac{v_0}{\omega_v} \right\|_1 \right).
\end{align*}
Similar arguments for $x\leq0$ yield
\bqs
\left|\int_{x-Kt}^{x+Kt} \widetilde{\cG^{12}}(t,x,y)v_0(y)\md y\right| \leq C \me^{-\theta t} \me^{\delta_v x}\left(\left\| \frac{v_0}{\omega_v} \right\|_\infty+\left\| \frac{v_0}{\omega_v} \right\|_1 \right),
\eqs
which concludes the proof.
\end{Proof}
Returning to the definition of $H(t,x)$ in \eqref{e:defH}, we now define
\begin{multline}\label{e:hdef}
h(t,x)\coloneqq\frac{1}{\omega(x)}\left(\int_{-\infty}^{x-Kt} \cG^{12}(t,x,y)v_0(y)\md y+\int_{x-Kt}^{x+Kt} \cG^{12}_\omega(t,x,y)v_0(y)\md y\right.\\ \left.+\int_{x+Kt}^{+\infty} \cG^{12}(t,x,y)v_0(y)\md y \right), 
\end{multline}
for all $x\in\R$ and $t\geq 1$. We then have the following estimates on $h(t,x)$.
\begin{lem}\label{lem:Estimate_h}
There exists $C>0$ such that for all $t\geq1$ and $x\in\R$
\bqq
|h(t,x)|\leq C \frac{1+|x|}{(1+t)^{3/2}} \int_\R (1+|y|)\frac{|v_0(y)|}{\omega(y)}\md y.
\label{eq:estimateh}
\eqq
\end{lem}
\begin{Proof}
The first and third integral in \eqref{e:hdef} are estimated using Proposition~\ref{timeGreen12} such that
\begin{align*}
\frac{1}{\omega(x) }\left| \int_{-\infty}^{x-Kt} \cG^{12}(t,x,y)v_0(y)\md y\right| &\leq \frac{C}{t^{1/4}} \int_{-\infty}^{x-Kt} \me^{-\frac{|x-y|^{4/3}}{\tilde{\kappa}t^{1/3}}}\frac{|v_0(y)|}{\omega(y)}\md y \leq \frac{C}{t^{1/4}} \me^{-\frac{K^{4/3}}{\tilde{\kappa}}t}\left\| \frac{v_0}{\omega} \right\|_1,\\
\frac{1}{\omega(x) }\left| \int_{x+Kt}^{+\infty} \cG^{12}(t,x,y)v_0(y)\md y\right| &\leq \frac{C}{t^{1/4}}\int_{x+Kt}^{+\infty} \me^{-\frac{|x-y|^{4/3}}{\tilde{\kappa}t^{1/3}}}\frac{|v_0(y)|}{\omega(y)}\md y \leq \frac{C}{t^{1/4}} \me^{-\frac{K^{4/3}}{\tilde{\kappa}}t}\left\| \frac{v_0}{\omega} \right\|_1,
\end{align*}
where we noticed that $\me^{-\frac{|x-y|^{4/3}}{\tilde{\kappa}t^{1/3}}}$ is maximized at the boundary.
Estimates on the second integral in \eqref{e:hdef} use Proposition~\ref{prop:timegreen12} such that
\bqs
\frac{1}{\omega(x) }\left| \int_{x-Kt}^{x+Kt} \cG^{12}_\omega(t,x,y)v_0(y)\md y \right| \leq C \me^{-\theta t} \left\| \frac{v_0}{\omega} \right\|_1 + C\frac{1+|x|}{(1+t)^{3/2}} \int_\R (1+|y|) \frac{|v_0(y)|}{\omega(y)}\md y,
\eqs
where we used $1+|x-y|\leq 1+|x|+|y|+|x||y|$ and sent the limits of integration to infinity. Combining the above estimates  gives \eqref{eq:estimateh}.
\end{Proof}
We may impose that $E(t,x)=0$ for all $0<t\leq 1$ and $x\in\R$ such that we have 
\bqs
H(t,x) = \omega(x)h(t,x)+E(t,x), \quad t>0, \quad x\in\R.
\eqs

\subsubsection{Nonlinear decay estimates}\label{s:nonlp}

We recall from the discussion in Section~\ref{sec:argoutline} that we look for solutions of \eqref{eq:main} that can be decomposed  $(u(t,x),v(t,x))=(Q_*(x)+P(t,x),v(t,x))$ where $P(t,x)$ stands for a perturbation around the critical front such that
\bqq
\left\{
\begin{array}{ll}
\partial_t P&= \cL_u P+\beta v +\mathcal{N}(P), \\
\partial_t v&= \cL_v v,
\end{array} \qquad t>0, \quad x\in\R.
\right.
\label{Nonlin}
\eqq
We introduce another smooth bounded weight functions $\rho_v>0$ through
\[ \rho_v(x) \coloneqq \left\{ \begin{array}{cc} \me^{-\frac{s_*}{2d} x}, & x>1, \\ 
1, & x=0,\\
\me^{\delta_v x}, & x<-1, \end{array}\right. \]
For brevity, we write  $L^q_\rho(\R)\coloneqq\left\{ u\in L^q_\mathrm{loc}(\R)~|~ \frac{u}{\rho}\in L^q(\R)\right\}$ for any $1\leq q \leq +\infty$ for any given positive function $\rho>0$. The Cauchy problem associated to \eqref{Nonlin} with initial condition $(P_0,v_0)$ such that $P_0\in L^1_\omega(\R)\cap L^\infty_\omega(\R)$ and $v_0\in L^1_{\rho_v}(\R)\cap L^\infty_{\rho_v}(\R)$ with 
\bqq
\int_\R |y|\frac{|P_0(y)|}{\omega(y)} +|y|\frac{|v_0(y)|}{\omega(y)}\md y <+\infty,
\label{CondIC}
\eqq
is locally well-posed in $L^\infty_\omega(\R)\times L^\infty_{\rho_v}(\R)$. We also remark that if $v_0\in L^1_{\rho_v}(\R)\cap L^\infty_{\rho_v}(\R)$ then it actually  belongs to both $L^1_{\omega_v}(\R)\cap L^\infty_{\omega_v}(\R)$ and $L^1_\omega(\R)\cap L^\infty_\omega(\R)$. As a consequence,  we let $T_*>0$ be the maximal time of existence of a solution $(P,v)\in L^\infty_\omega(\R)\times L^\infty_{\rho_v}(\R)$ with initial condition $(P_0,v_0)$ such that $P_0\in L^1_\omega(\R)\cap L^\infty_\omega(\R)$ and $v_0\in L^1_{\omega_v}(\R)\cap L^\infty_{\rho_v}(\R)$ further satisfying \eqref{CondIC}. Solutions of this nonlinear system with initial condition $(P_0,v_0)$ can be expressed using Duhamel's formula,
\bqs
P(t,x)=\int_\R \cG^{11}(t,x,y)P_0(y)\md y+\int_\R \cG^{12}(t,x,y)v_0(y)\md y+\int_0^t\int_\R \cG^{11}(t-s,x,y)\mathcal{N}(P)(s,y)\md y\md s,
\eqs
together with
\bqs
v(t,x)=\int_\R \cG^{22}(t,x,y)v_0(y)\md y.
\eqs
Actually, one can easily check that the solution $v$ is globally well-posed in $L^\infty_{\omega_v}(\R)$ and a direct application of Lemma~\ref{timeGreen22} gives us the global in time bound
\begin{equation}
\left\|\frac{ v(t)}{\omega_v}\right\|_\infty \leq C \me^{-\theta t}  \left\| \frac{v_0}{\omega_v} \right\|_\infty, \quad t>0. \label{e:vdecay}
\end{equation}
We now concentrate on $P(t,x)$. From the previous section, we know that
\bqs
H(t,x) = \int_\R \cG^{12}(t,x,y)v_0(y)\md y= \omega(x)h(t,x)+E(t,x), \quad t>0, \quad x\in\R.
\eqs
It then follows that
\bqs
P(t,x)-E(t,x)=\int_\R \cG^{11}(t,x,y)P_0(y)\md y+\omega(x)h(t,x)+\int_0^t\int_\R \cG^{11}(t-s,x,y)\mathcal{N}(P)(s,y)\md y\md s,
\eqs
and so we define $P(t,x)-E(t,x)=\omega(x)p(t,x)$ which solves
\bqq
p(t,x)=\int_\R \widetilde{\cG}^{11}(t,x,y)p_0(y)\md y+h(t,x)+\int_0^t\int_\R \widetilde{\cG}^{11}(t-s,x,y)\omega(y)^{-1}\mathcal{N}(\omega p+E)(s,y)\md y\md s,
\label{eqDefp}
\eqq
with $p(0,x)=p_0(x)=\frac{P_0(x)}{\omega(x)}$ as we have set $E(t,x)=0$ for $0\leq t \leq 1$. For $t\in[0,T_*)$, we define
\bqs
\Theta(t)\coloneqq\underset{0\leq s \leq t}{\sup}~\underset{x\in\R}{\sup}~ (1+s)^{3/2}\frac{|p(s,x)|}{1+|x|},
\eqs
well-defined since $p\in L^\infty(\R)$. We now bound all terms in \eqref{eqDefp} for $0<t\leq 1$ and for $t>1$.

\paragraph{Short time bound $0<t\leq 1$.} For the short time bound, we use 
Proposition~\ref{timeGreen11} to obtain
\bqs
\left|\int_\R\widetilde{\cG}^{11}(t,x,y)p_0(y)\md y\right|\leq C \|p_0\|_\infty = C \left\|\frac{P_0}{\omega}\right\|_\infty,
\eqs
and Lemma~\ref{lem:shorttimeh} to obtain
\bqs
|h(t,x)|\leq C \left\|\frac{v_0}{\omega}\right\|_\infty, \quad x \in\R.
\eqs
For the last  term in \eqref{eqDefp}, we recall that
\bqs
\frac{1}{\omega}\mathcal{N}(\omega p+E)=-\alpha\frac{E^2}{\omega}\left(3Q_*+E\right)-3\alpha\left(2Q_*+E\right)pE-3\alpha\left(Q_*+E\right)\omega p^2-\alpha \omega^2 p^3,
\eqs
and that $E(t,x)=0$ for all $0<t\leq1$ such that for all $0\leq s \leq t \leq 1$,
\bqs
\left| \frac{1}{\omega(y)}\mathcal{N}(\omega(y) p(s,y)+E(s,y)) \right| \leq C \left( \Theta(t)^2 (1+|y|)^2\omega(y) + \Theta(t)^3 (1+|y|)^3 \omega(y)^2\right), \quad y\in\R.
\eqs
As a consequence, we have
\begin{align*}
\left| \int_\R \widetilde{\cG}^{11}(t-s,x,y)\omega(y)^{-1}\mathcal{N}(\omega p+E)(s,y)\md y \right| &\leq \frac{C}{(t-s)^{1/2}} \Theta(t)^2 \int_\R \me^{-\frac{|x-y|^2}{\kappa_u(t-s)}}(1+|y|)^2\omega(y)\md y\\
&~~~+\frac{C}{(t-s)^{1/2}} \Theta(t)^3 \int_\R \me^{-\frac{|x-y|^2}{\kappa_u(t-s)}}(1+|y|)^3\omega(y)^2\md y,
\end{align*}
which finally gives
\bqs
\Theta(t)\leq C\left( \left\|\frac{P_0}{\omega}\right\|_\infty+\left\|\frac{v_0}{\omega}\right\|_\infty\right)+C\Theta(t)^2 \int_\R (1+|y|)^2\omega(y)\md y+C\Theta(t)^3 \int_\R (1+|y|)^3\omega(y)^2\md y.
\eqs

\paragraph{Large time bound $t>1$.} For the large time bound we use 
Proposition~\ref{timeGreen11} to obtain
\bqs
\left|\int_\R\widetilde{\cG}^{11}(t,x,y)p_0(y)\md y\right|\leq C \frac{1+|x|}{(1+t)^{3/2}} \int_\R (1+|y|)\frac{|P_0(y)|}{\omega(y)}\md y,
\eqs
and Lemma~\ref{lem:Estimate_h} then gives
\bqs
|h(t,x)|\leq C \frac{1+|x|}{(1+t)^{3/2}} \int_\R (1+|y|)\frac{|v_0(y)|}{\omega(y)}\md y.
\eqs
To derive bounds for the last term in \eqref{eqDefp}, we  use estimate \eqref{eqEstimE} from Lemma~\ref{lemEstimE}
\bqs
 |E(t,x)|\leq C \me^{-\theta t}\omega_v(x) \left( \left\| \frac{v_0}{\omega_v} \right\|_\infty+\left\| \frac{v_0}{\omega_v} \right\|_1\right),
\eqs
and find
\begin{multline*}
\left| \int_\R \widetilde{\cG}^{11}(t-s,x,y) Q_*(y) \frac{|E(s,y)|^2}{\omega(y)}\md y\right| \lesssim \frac{(1+|x|)\me^{-2\theta s}}{(1+t-s)^{3/2}} \left(\int_\R (1+|y|) Q_*(y) \frac{\omega_v^2(y)}{\omega(y)}\md y\right)\\
\quad \cdot \left( \left\| \frac{v_0}{\omega_v} \right\|_\infty+\left\| \frac{v_0}{\omega_v} \right\|_1\right)^2,
\end{multline*}
\begin{align*}\left| \int_\R \widetilde{\cG}^{11}(t-s,x,y) \frac{|E(s,y)|^3}{\omega(y)}\md y\right| &\lesssim \frac{(1+|x|)\me^{-3\theta s}}{(1+t-s)^{3/2}}\left( \int_\R (1+|y|) \frac{\omega_v^3(y)}{\omega(y)}\md y\right)\left( \left\| \frac{v_0}{\omega_v} \right\|_\infty+\left\| \frac{v_0}{\omega_v} \right\|_1\right)^3,
\end{align*}
where
\bqs
\int_\R (1+|y|) Q_*(y) \frac{\omega_v^2(y)}{\omega(y)}\md y<\infty \text{ and } \int_\R (1+|y|) \frac{\omega_v^3(y)}{\omega(y)}\md y <\infty,
\eqs
thanks to our assumption that $\frac{s_*}{2d}<3\gamma_v$ and the fact that the critical front satisfies $\frac{Q_*(x)}{\omega(x)}\sim x $ as $x\rightarrow+\infty$. As a consequence, we obtain the following bounds
\begin{align*}
\left| \int_0^t\int_\R \widetilde{\cG}^{11}(t-s,x,y) Q_*(y) \frac{|E(s,y)|^2}{\omega(y)}\md y\md s\right| & \lesssim \frac{1+|x|}{(1+t)^{3/2}} \left( \left\| \frac{v_0}{\omega_v} \right\|_\infty+\left\| \frac{v_0}{\omega_v} \right\|_1\right)^2,\\
\left| \int_0^t \int_\R \widetilde{\cG}^{11}(t-s,x,y) \frac{|E(s,y)|^3}{\omega(y)} \md y\md s\right| & \lesssim \frac{1+|x|}{(1+t)^{3/2}} \left( \left\| \frac{v_0}{\omega_v} \right\|_\infty+\left\| \frac{v_0}{\omega_v} \right\|_1\right)^3.
\end{align*}
We now treat the term $-3\alpha(Q_*+E)pE$ in the nonlinearity. Using the fact that both $Q_*$ and $E$ are bounded one only needs estimates for $pE$. We have
\begin{multline*}
\left| \int_\R \widetilde{\cG}^{11}(t-s,x,y)|p(s,y)|E(s,y)|\md y\right| \leq C \Theta(t)\frac{(1+|x|)\left( \int_\R (1+|y|) \omega_v(y)\md y\right)}{(1+t-s)^{3/2}(1+s)^{3/2}}\me^{-\theta s}\\
 \  \qquad \cdot \left( \left\| \frac{v_0}{\omega_v} \right\|_\infty+\left\| \frac{v_0}{\omega_v} \right\|_1\right),
\end{multline*}
such that
\bqs
\left| \int_0^t\int_\R \widetilde{\cG}^{11}(t-s,x,y)|p(s,y)|E(s,y)|\md y\md s\right| \lesssim \frac{1+|x|}{(1+t)^{3/2}}\Theta(t)\left( \left\| \frac{v_0}{\omega_v} \right\|_\infty+\left\| \frac{v_0}{\omega_v} \right\|_1\right).
\eqs
Finally nonlinear terms in $p$ are handled similarly, 
\begin{align*}
\left| \int_\R \widetilde{\cG}^{11}(t-s,x,y)\omega(y)|p(s,y)|^2\md y\right| &\leq C \Theta(t)^2\frac{1+|x|}{(1+t-s)^{3/2}(1+s)^3}\left( \int_\R (1+|y|)^2 \omega(y)\md y\right),\\
\left| \int_\R \widetilde{\cG}^{11}(t-s,x,y)\omega(y)^2|p(s,y)|^3\md y\right| &\leq C \Theta(t)^3\frac{1+|x|}{(1+t-s)^{3/2}(1+s)^{9/2}}\left( \int_\R (1+|y|)^3 \omega^2(y)\md y\right),
\end{align*}
which leads to
\begin{align*}
\left|\int_0^t \int_\R \widetilde{\cG}^{11}(t-s,x,y)\omega(y)|p(s,y)|^2\md y \md s\right| &\lesssim \Theta(t)^2\frac{1+|x|}{(1+t)^{3/2}},\\
\left| \int_0^t \int_\R \widetilde{\cG}^{11}(t-s,x,y)\omega(y)^2|p(s,y)|^3\md y \md s\right| &\lesssim \Theta(t)^3\frac{1+|x|}{(1+t)^{3/2}}.
\end{align*}
Combining all the above estimates, we obtain the bound
\begin{align*}
\Theta(t)&\leq C\int_\R (1+|y|)\left(\frac{|v_0(y)|}{\omega(y)}+\frac{|v_0(y)|}{\omega(y)}\right)\md y+C\left( \left\| \frac{v_0}{\omega_v} \right\|_\infty+\left\| \frac{v_0}{\omega_v} \right\|_1\right)^2+C\left( \left\| \frac{v_0}{\omega_v} \right\|_\infty+\left\| \frac{v_0}{\omega_v} \right\|_1\right)^3\\
&~~~+C\Theta(t)\left( \left\| \frac{v_0}{\omega_v} \right\|_\infty+\left\| \frac{v_0}{\omega_v} \right\|_1\right)+C\Theta(t)^2+C\Theta(t)^3.
\end{align*}

\paragraph{Conclusion of the proof.} We can now combine our small and large time bound to deduce that for all $t\in[0,T_*)$
\bqs
\Theta(t) \leq C_0 \Omega_0+ C_1\Theta(t)\left( \left\| \frac{v_0}{\omega_v} \right\|_\infty+\left\| \frac{v_0}{\omega_v} \right\|_1\right)+C_2\Theta(t)^2+C_3\Theta(t)^3,
\eqs
for some positive constants $C_{j}>0$ and 
\begin{align*}
\Omega_0\coloneqq&\left\|\frac{P_0}{\omega}\right\|_\infty+\left\|\frac{v_0}{\omega}\right\|_\infty+\int_\R (1+|y|)\left(\frac{|v_0(y)|}{\omega(y)}+\frac{|v_0(y)|}{\omega(y)}\right)\md y+\left( \left\| \frac{v_0}{\omega_v} \right\|_\infty+\left\| \frac{v_0}{\omega_v} \right\|_1\right)^2\\
&+\left( \left\| \frac{v_0}{\omega_v} \right\|_\infty+\left\| \frac{v_0}{\omega_v} \right\|_1\right)^3
\end{align*}
which only depends on the initial condition $(P_0,v_0)$. Therefore, assuming that $(P_0,v_0)$ satisfies 
\bqs
C_1\left( \left\| \frac{v_0}{\omega_v} \right\|_\infty+\left\| \frac{v_0}{\omega_v} \right\|_1\right) < \frac{1}{2}, \quad 4C_0\Omega_0<1, \text{ and } 16C_0C_2\Omega_0+64C_0^2C_3\Omega_0^2<1,
\eqs
we claim that
\bqs
\Theta(t)\leq 4C_0\Omega_0<1, \quad t\in[0,T_*).
\eqs
To see this, note that $C_1\left( \left\| \frac{v_0}{\omega_v} \right\|_\infty+\left\| \frac{v_0}{\omega_v} \right\|_1\right)< \frac{1}{2}$ implies
\begin{equation}\label{e:above}
\Theta(t)<2C_0 \Omega_0+2C_2\Theta(t)^2+2C_3\Theta(t)^3.
\end{equation}
Next, by eventually taking $C_0$ larger, we can always assume that
\bqs
\Theta(0)=\underset{x\in\R}{\sup}~\frac{p_0(x)}{1+|x|}\leq \left\|\frac{P_0}{\omega}\right\|_\infty <\Omega_0 <4C_0\Omega_0,
\eqs 
such that the continuity of $\Theta(t)$ implies that, for small times, we have $\Theta(t)<4C_0\Omega_0$. Suppose now that there exists $T>0$ such that $\Theta(T)=4C_0\Omega_0$ for the first time, then from \eqref{e:above} we obtain
\bqs
\Theta(T)\leq 2C_0 \Omega_0+2C_2 (4C_0 \Omega_0)^2+2C_3(4C_0 \Omega_0)^3=2C_0\Omega_0\left(1+16C_0C_2\Omega_0+64C_0^2C_3\Omega_0^2\right)<4C_0\Omega,
\eqs
a contradiction, which proves the claim. As a consequence, the maximal time of existence is $T_*=+\infty$ and the perturbation $p$ satisfies
\bqs
\underset{t\geq0}{\sup}~\underset{x\in\R}{\sup}~ (1+t)^{3/2} \frac{|p(t,x)|}{1+|x|}<4C_0\Omega_0.
\eqs
Decay of the $v$-component was established in \eqref{e:vdecay}, which concludes the proof of the decay estimates in our main result. 
\subsubsection{Instability in fixed weights}\label{s:inst} 
We show the last statement of Theorem \ref{thmMain}.
We argue by contradiction and assume $P$ is small, bounded in the weighted $L^\infty$ space, with vanishing initial conditions. We then find that $P$ solves the inhomogeneous convection-diffusion equation 
\[
\partial_t P=d\partial_{xx}P+s_* \partial_xP + g_p+v, 
\]
where $g_p$ incorporates all remainder terms depending on $P$ and is hence uniformly bounded and small in the weighted space. 
On the other hand, we know that $v$ grows exponentially in the weighted space, with explicit expressions in Fourier space.
One therefore easily constructs solutions to the modified equation
\[
\partial_t  \tilde{P}=d\partial_{xx}\tilde{P}+s_* \partial_x\tilde{P} +v, 
\]
that grow exponentially in the weighted space, solving again explicitly in Fourier space. Since $g_p$ was bounded and small, the solution to 
\[
\partial_t  \hat{P}=d\partial_{xx}\hat{P}+s_* \partial_x\hat{P} +g_p, 
\]
is bounded in time, contradicting the assumption that $P=\tilde{P}+\hat{P}$ is small, bounded. 

\section{Absolute Spectra}\label{sec:AbsSpec}

In order to characterize regions $\mathcal{R}_\mathrm{rem}$ and $\mathcal{R}_\mathrm{abs}$ further, and prepare for the discussion of numerical simulations and Conjecture~\ref{prop:sabs}, we discuss the absolute spectrum associated with $\Sabs(\calL^+)$, the linearization at the unstable state. We start with a brief review of absolute spectra, in  Section~\ref{sec:RemnantInstabilities,AbsoluteSpectrumAndDoubleRoots}, including a characterization of absolute spectra for the $u$- and $v$-component and conclude the proof of Theorem~\ref{thm:paraminformal} in   Section~\ref{sec:absoluteformain}. 

\subsection{Absolute spectra and double roots --- review and the case of KPP and SH}\label{sec:RemnantInstabilities,AbsoluteSpectrumAndDoubleRoots}

Recall from Section~\ref{sec:EssSpec} the dispersion relation~\eqref{eq:DispersionRelation}, given as a product of dispersion relations for $u$- and $v$-components, separately. Associated with the dispersion relation is the Morse index $i_\infty$, the number of roots $\nu$ to $D(\lambda,\nu)=0$, with $\Re\nu>0$ for a fixed $\Re\lambda\gg1$, to the right of the essential spectrum. In our case, one quickly sees $i_\infty=3$, while for the $u$-component alone $i^u_\infty=1$ and for the $v$-component $i_\infty^v=2$. 
We may also, allowing for discontinuities of the labeling in $\lambda$ and some ambiguities, order all roots $\nu=\tilde{\nu}_j(\lambda)$ by increasing real part, repeating by multiplicity when necessary, 
\[
\Re\tilde{\nu}_1(\lambda)\leq \Re\tilde{\nu}_2(\lambda)\leq \ldots                                                                                                                                                                                                                                                               
\]
Following~\cite{sandstede00}, we say 
\begin{equation}\label{e:mi}
 \lambda\in\Sigma_\mathrm{abs}\ \Longleftrightarrow \ \Re\tilde{\nu}_{i_\infty}(\lambda)=\Re\tilde{\nu}_{i_\infty+1}(\lambda). 
\end{equation}
It is crucial to recognize that this ordering  and labeling depends in our case on whether the $u$- and $v$-component are considered separately, or together. In fact, we previously ordered roots from the $v$-equation alone as $\Re\rho_1(\lambda)\leq \Re\rho_2(\lambda)\leq \Re\rho_3(\lambda)\leq \Re\rho_4(\lambda)$ (see \eqref{eq:rootsvorderingabs}), and roots from the $u$-equation $\Re\nu^0_+(\lambda)\geq \Re\nu^0_-(\lambda)$. The ordering in \eqref{e:mi} refers to the combined ordering, and we used $\tilde{\nu}_j,\,1\leq j\leq 6$, to distinguish this ordering from the ordering $\rho_j,\, 1\leq j\leq 4$, to distinguish from earlier notation for the roots in the $v$-equation. A key example of a subtlety related to this combined ordering, which will appear later, is that, considering the $v$-component alone, we may see $\Re\rho_1=\Re\rho_2$, which would not contribute to the absolute spectrum of the $v$-equation since $i_\infty^v=2$. If however the roots $\nu_\pm^0$ from the $u$-equation lie to the left, $\Re\nu_-^0\leq \Re\nu_+^0\leq \Re\rho_1=\Re\rho_2<\Re\rho_3\leq \Re\rho_4$, then $\lambda$ would belong to the absolute spectrum of the combined $u$-$v$-system. 

It turns out that $\Sigma_\mathrm{abs}(\mathcal{L}^+)$ consists of parameterized curves, solving 
\begin{equation}\label{e:abscont}
 D(\lambda,\nu)=0,\qquad D(\lambda,\nu+\mbi k)=0,\qquad \lambda,\nu\in\C, 
\end{equation}
with parameter $k\in\R$ and for which~\eqref{e:mi} holds. Note that $\Sigma_\mathrm{abs}$ is to the left of $\Sigma^\eta_\mathrm{ess}$ for any $\eta$ in the sense that any curve from $\Sigma_\mathrm{abs}$ to $\lambda=+\infty$ needs to cross $\Sigma_\mathrm{ess}$.
We remark in passing that the absolute spectrum was originally introduced as the continuous part of the limit of spectra in bounded domains in the limit of infinite domain size~\cite{sandstede00}, but we will not rely on this property, here.

Generic singularities of these parameterized curves are double roots, where $k=0$, and triple points, where $\Re\tilde{\nu}_{\ell-1}=\Re\tilde{\nu}_{\ell}=\Re\tilde{\nu}_{\ell+1}$, and $\ell=i_\infty$ or $\ell-1=i_\infty$; see~\cite{rademacher07}. 

At double roots, curves of absolute spectrum typically end: continuing in $k$ through 0  simply interchanges the two roots $\nu(\lambda)$ solving \eqref{e:abscont}. Interestingly, curves of absolute spectrum simply pass smoothly through resonance poles, where individual roots $\nu$ can be continued analytically.   We discussed double roots in Section~\ref{sec:ptw}. Loosely speaking, double roots at rightmost points of the absolute spectrum automatically satisfy the pinching condition whenever they satisfy the Morse index condition.

Near a triple point, we can continue the equal real part condition between any of the three possible pairs of roots that have equal real part at the triple point to find 3 curves that cross in the triple point. Following such a curve, the Morse index of the pair with equal real part changes~\cite{rademacher07} and only the part of the curve to one side of the triple point belongs to the absolute spectrum.

\paragraph{Absolute spectra of KPP.} Here, $D(\lambda,\nu)=d\nu^2+s\nu+\alpha-\lambda$, $i_\mathrm{\infty}=1$, and $\Sigma_\mathrm{abs}(\mathcal{L}^+_u)=\{\lambda|\,   \lambda \leq  \alpha  - s^2/(4d)\}$, that is, when the discriminant is negative. Notice that the rightmost point in the absolute spectrum is the double root $\alpha-s^2/(4d)$, which vanishes precisely at the KPP spreading speed. In this sense, the concept of remnant instabilities, instabilities of the absolute spectrum, and pointwise instabilities coincide for this problem. 

\paragraph{Absolute spectra of SH.} Here, $D(\lambda,\nu)=-(\nu^2+1)^2+s \nu + \mu-\lambda$. Solutions to~\eqref{e:abscont} are explicit but rather intractable. In fact, it is not immediately clear that in this case the rightmost points of the absolute spectrum are pinched double roots. However, since the essential spectrum has negative real part, so does the absolute spectrum. Although the exact form of the absolute spectrum is therefore not relevant to our analysis, we present a quick summary of the structure and refer to the appendix for more detailed formulas.

\begin{figure}
    \centering
     \includegraphics[width=0.32\textwidth]{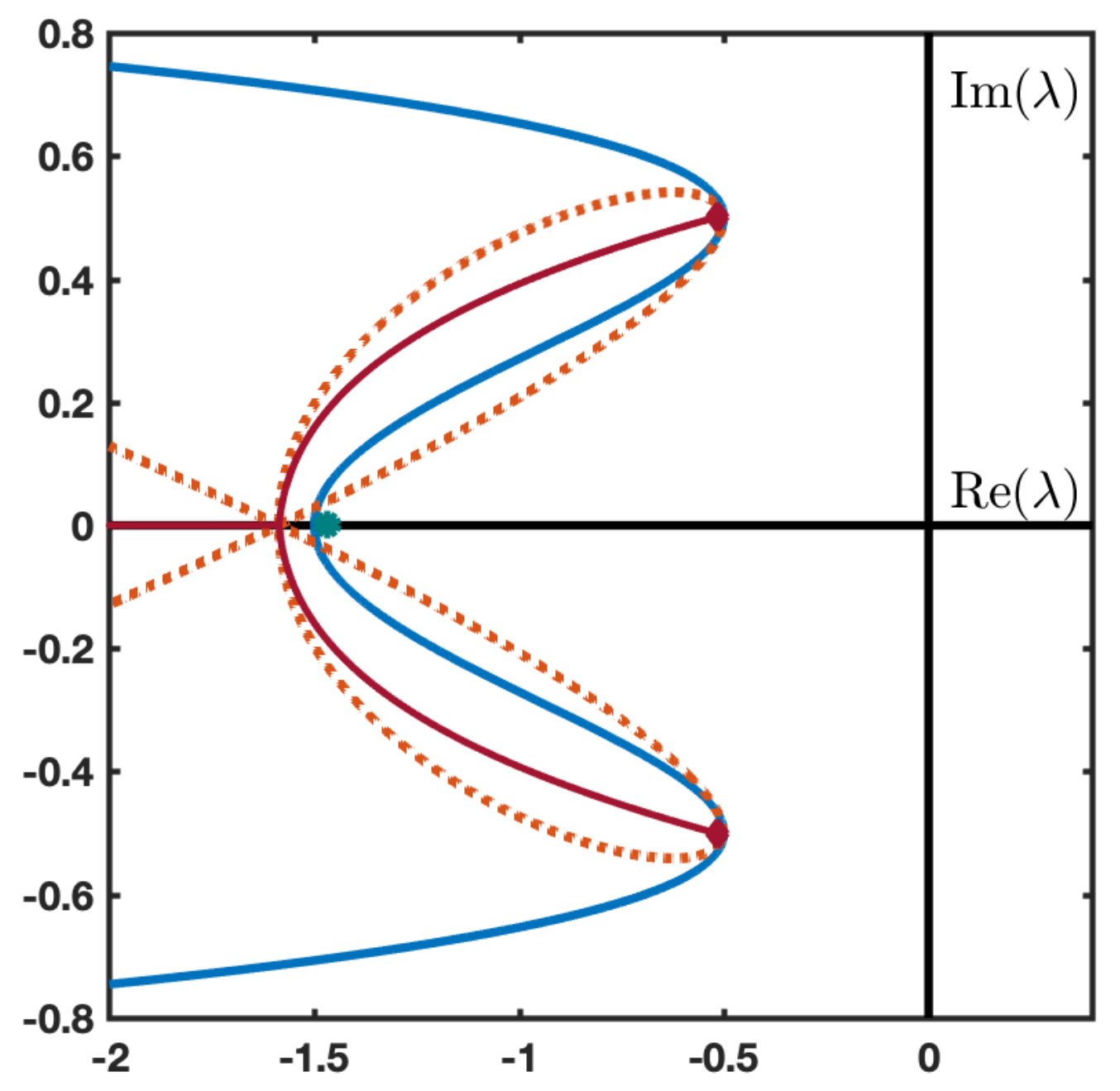}
    \includegraphics[width=0.31\textwidth]{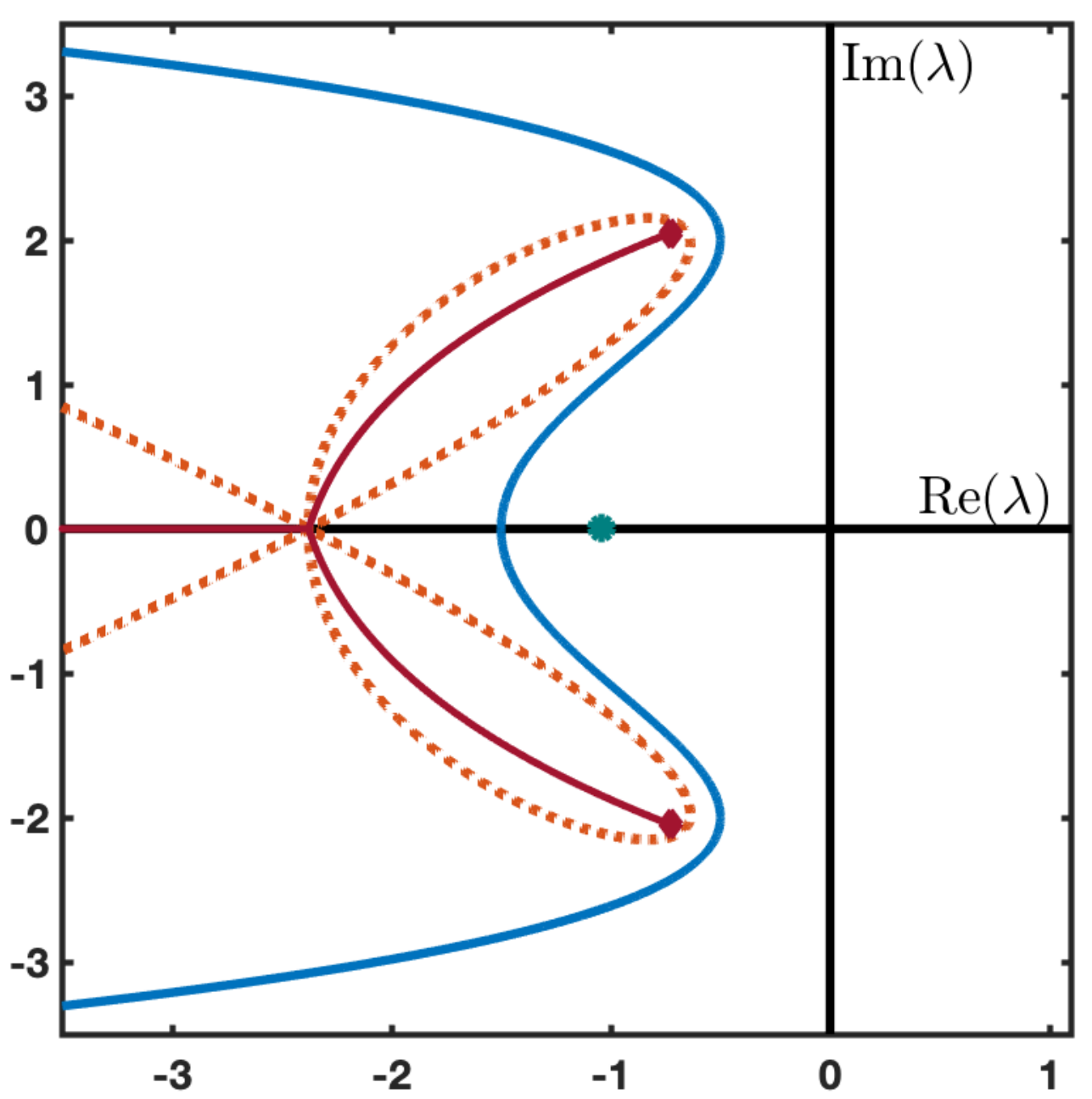}
     \includegraphics[width=0.32\textwidth]{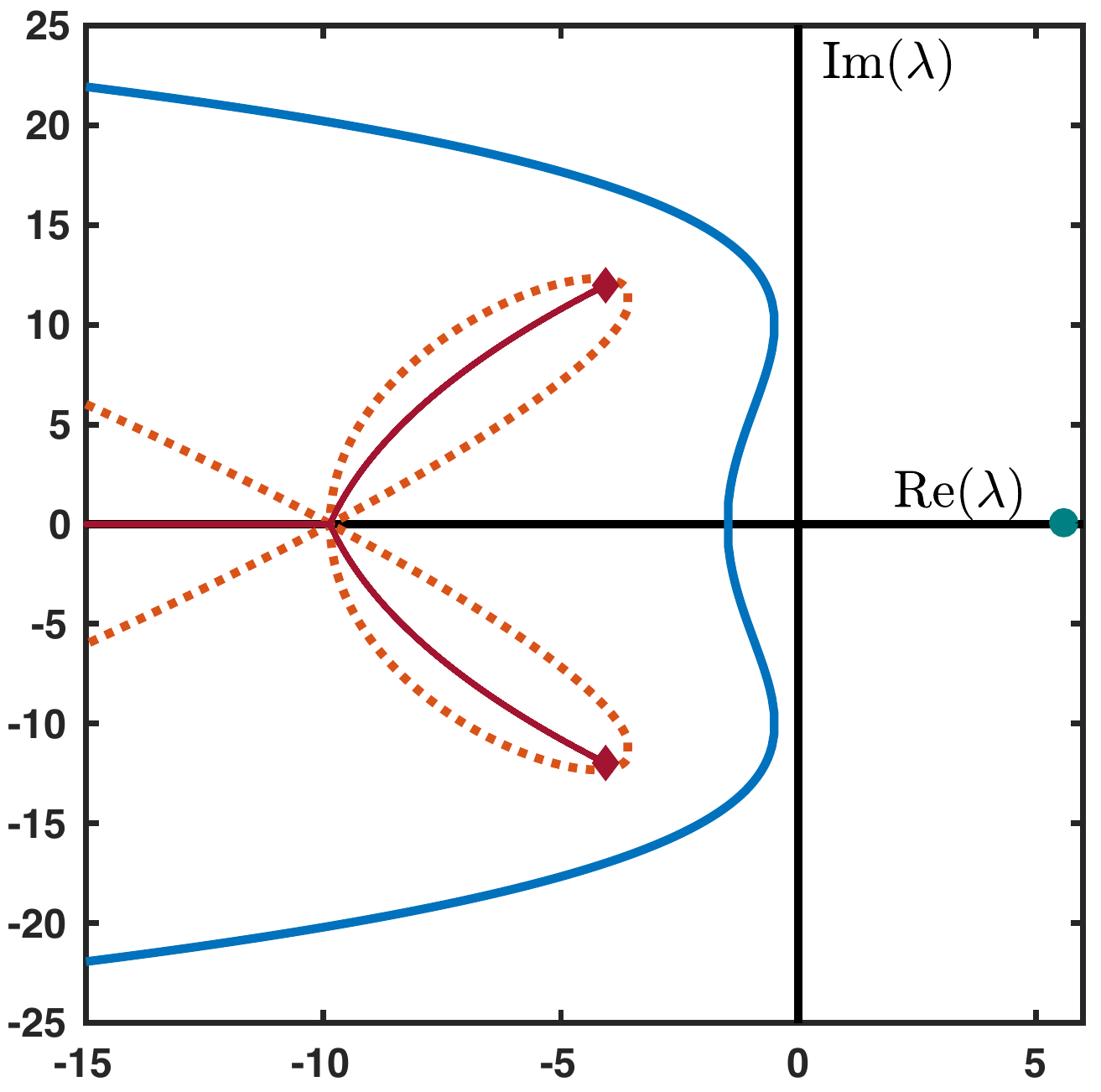}
    \caption{\small Essential (blue solid), weighted essential with weight $\eta=\eta^\textnormal{tr}_v$ (orange dotted), and absolute spectrum (red solid) of $\calL^+_v$ for parameters $d=1,\, \alpha=1,\, \mu=-1/2$, and $s=s_*/4,s_*,5s_*$ (left to right) with $s_*=2\sqrt{\alpha d}=2$; also shown pinched double roots in the absolute spectrum (dark red diamonds) and the non-pinched double root, not part of the absolute spectrum (cyan dot). The triple point is located at $\lambda^v_\textnormal{tr}=-1.586,-2.447,-9.863$ with corresponding weight  $\eta^\textnormal{tr}_v=-0.117,-0.326,-0.7101$. Note the different scales on the coordinate axes. }
    \label{fig: SH abs spec}
\end{figure}

\begin{lem}[Absolute spectrum of SH]\label{l:shabs}For any $s>0$ and $\mu<0$, the absolute spectrum of SH consists of three curves emanating from a triple point located at 
\begin{equation}\label{eq:TriplePointSwiftHohenberg}
\lambda_v^\mathrm{tr}=\mu-1+\frac{s\left(-4\sqrt[3]{225}+\sqrt[3]{15}C_1^2\right)}{30C_1}-\frac{\left(-4\sqrt[3]{225}+\sqrt[3]{15}C_1^2\right)^2}{450C_1^2}-\frac{\left( -4\sqrt[3]{225}+\sqrt[3]{15}C_1^2\right)^4}{810000C_1^4},
\end{equation}
where $C_1\coloneqq(-45s+\sqrt{960+2025s^2})^{1/3}$. The first curve is real and consists of the unbounded interval $\lambda\in (-\infty,\lambda_v^\mathrm{tr}]$.
The other two curves are complex conjugate and connect the triple point to the two pinched double roots 
\begin{equation}\label{eq:DoubleRootsSwiftHohenberg}
 \lambda_v^\mathrm{dr}=\mu-\frac{1}{3}\pm\frac{(27s^2-8)(\mbi\sqrt{3}-1)}{6C_2}+\frac{1}{48}(1\pm\mbi\sqrt{3})C_2,
\end{equation}
with
\[C_2\coloneqq \left(512+4320s^2-729s^4+3\sqrt{3}\sqrt{s^2(64+27s^2)^3}\right)^{\frac{1}{3}},\]
which also maximize the real part of the absolute spectrum. Moreover, remnant instabilities imply pointwise instabilities, that is, values of $\mu$ such that we have a remnant instability, unstable absolute spectrum, or and unstable pinched double root coincide.
\end{lem}
We refer to the appendix, Section~\ref{sec:SwiftHohenbergAbsSpec}, for a proof of this lemma and to Figure \ref{fig: SH abs spec} for an illustration. 
The result also motivates the model system under consideration here as an extension of the Swift-Hohenberg equation necessary to find remnant instabilities that are not caused by pointwise instabilities. 

\subsection{Remnant instabilities, absolute spectra, and pointwise instabilities in the full system}\label{sec:absoluteformain}
We next study the absolute spectrum, also in relation to remnant and pointwise instabilities, for 
\[
 D(\lambda,\nu)=(d\nu^2+s\nu+\alpha-\lambda)\cdot(-\nu^4-2\nu^2+s\nu-1+\mu-\lambda ),
\]
The following three results will together conclude the proof of  Theorem \ref{thm:paraminformal}. 

\begin{prop}[$\mathcal{R}_\mathrm{rem}\neq \emptyset$]\label{p:rrem}
 For $\mu-\mu^\mathrm{rem}(\alpha,d)>0$, sufficiently small, the origin has a remnant instability but the real part of the absolute spectrum is negative. 
\end{prop}

\begin{figure}[!h]
\centering
\subfigure[$(d,\alpha,\mu)=(1,1,-2)$.]{\includegraphics[width=0.24\textwidth]{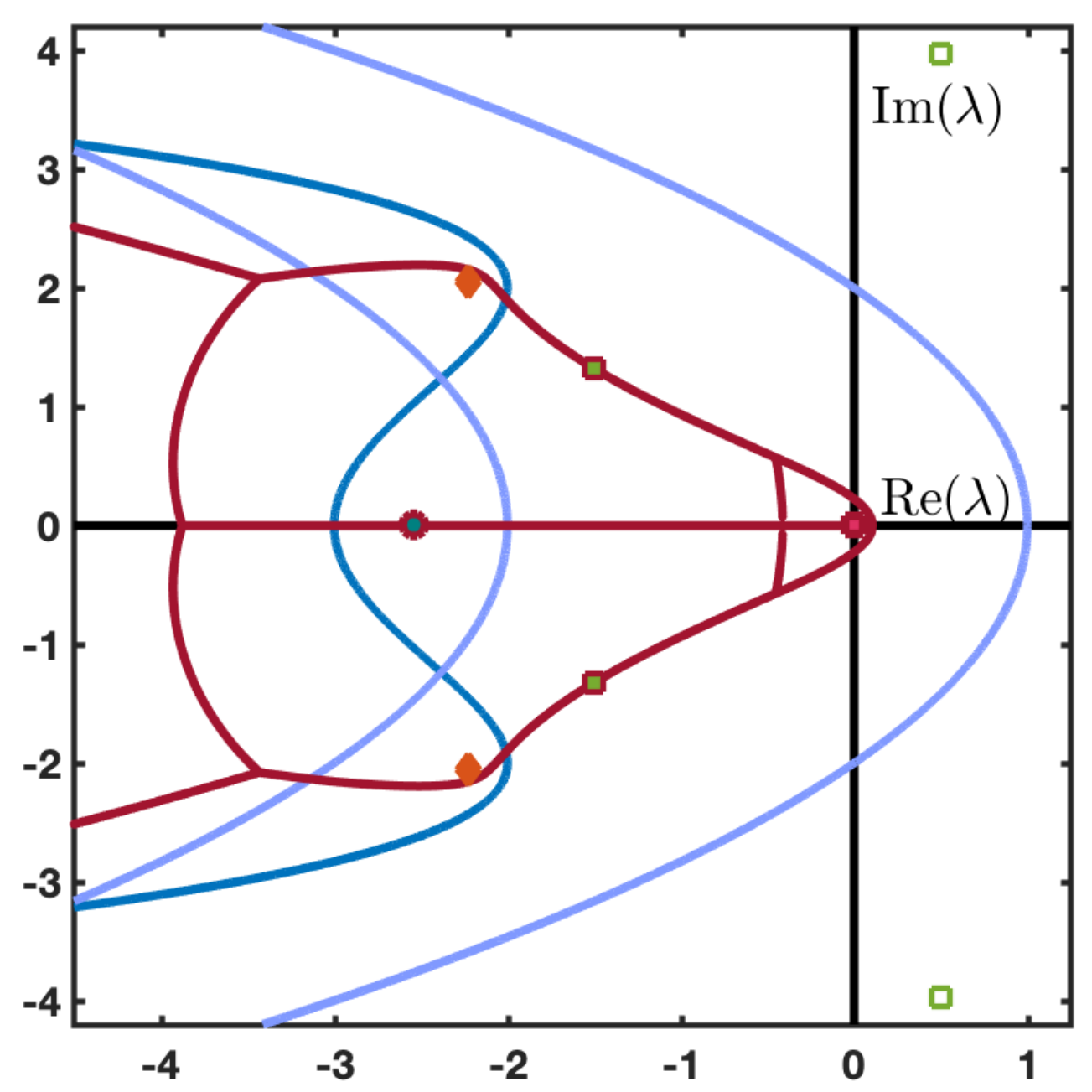}}
\subfigure[$(d,\alpha,\mu)=(1,1,-0.75)$.]{\includegraphics[width=0.245\textwidth]{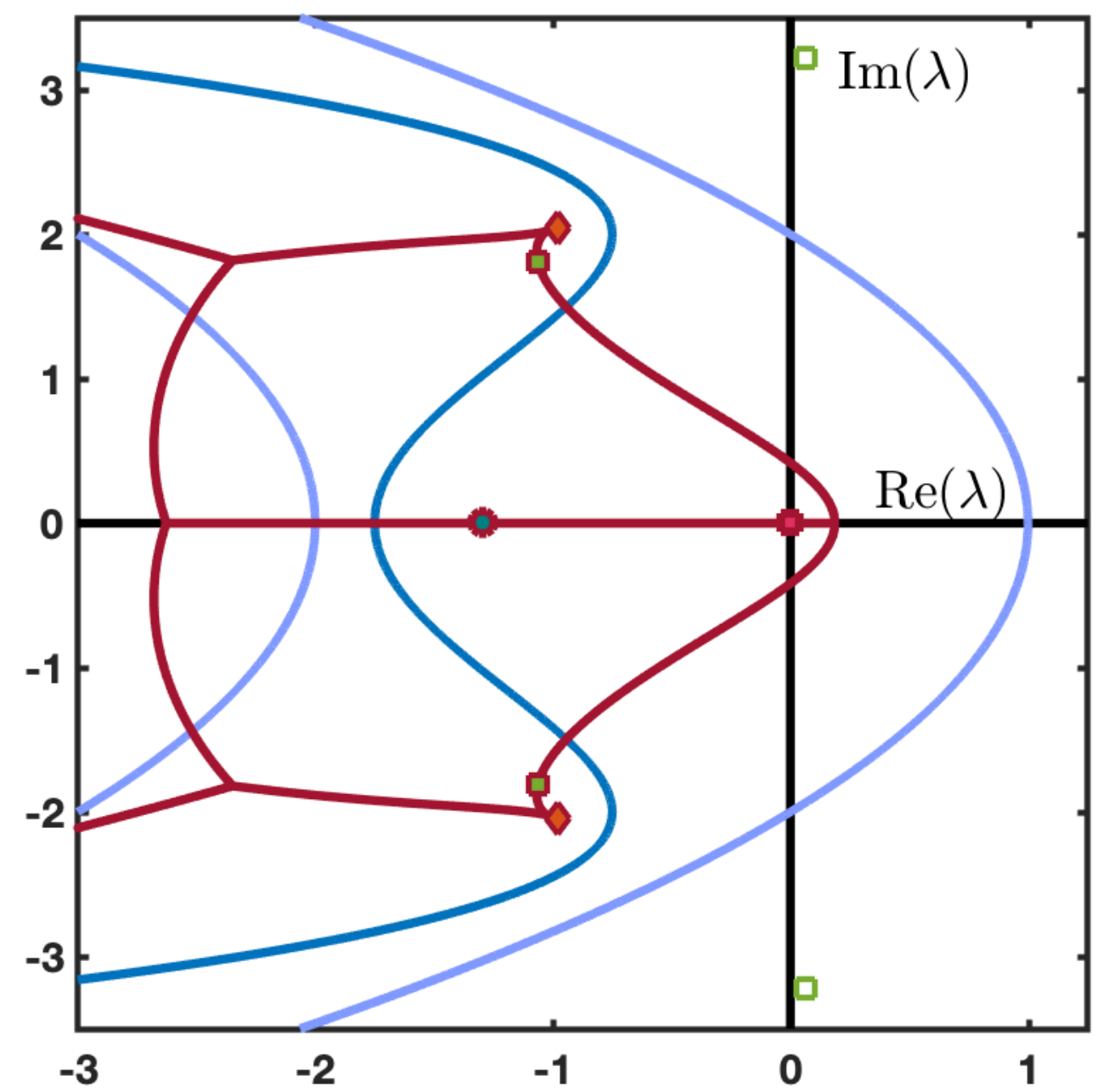}}
 \subfigure[$(d,\alpha,\mu)=(1,1,-0.4)$.]{\includegraphics[width=0.24\textwidth]{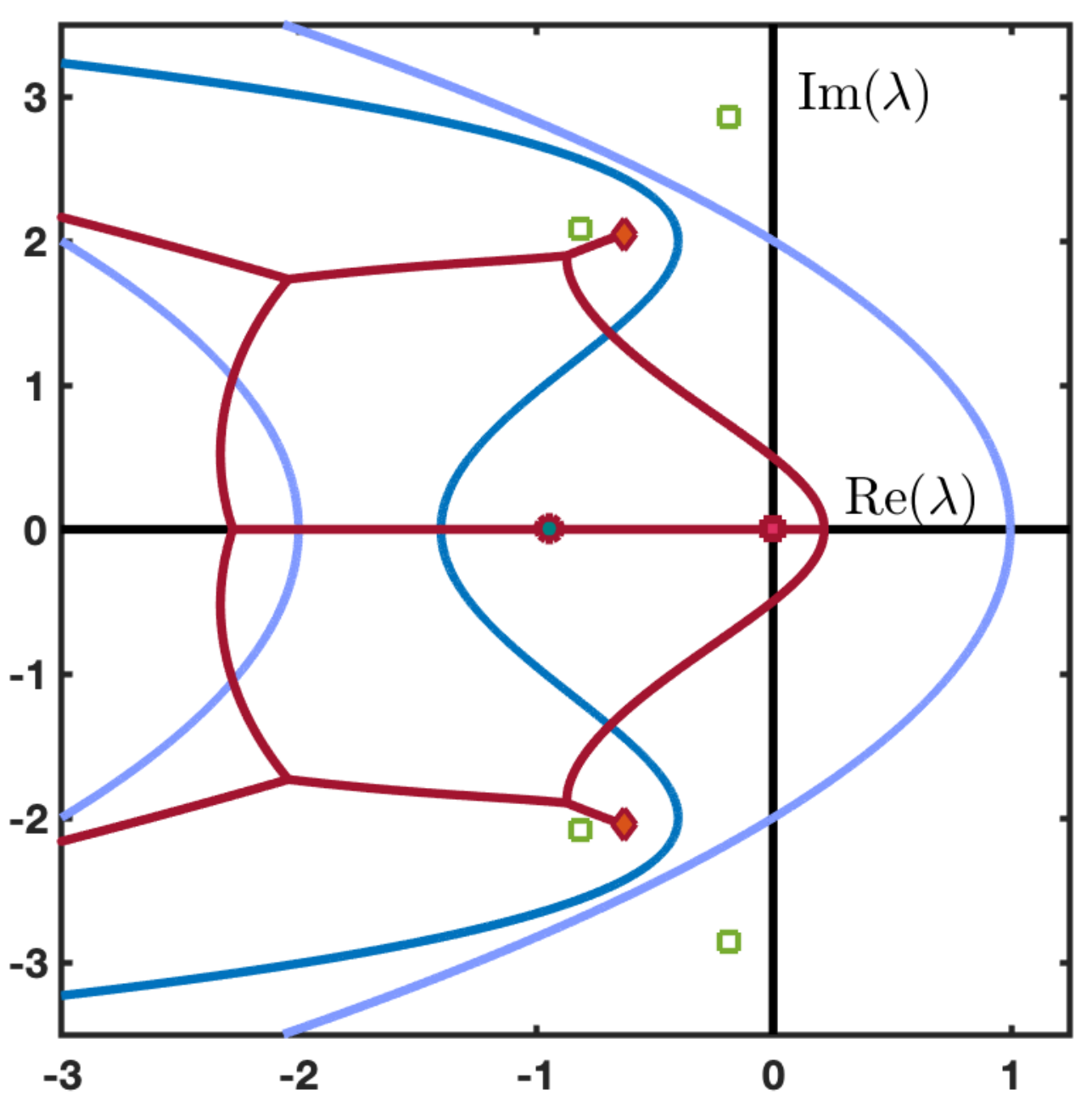}}
\subfigure[$(d,\alpha,\mu)=(0.5,2,-1)$.]{\includegraphics[width=0.245\textwidth]{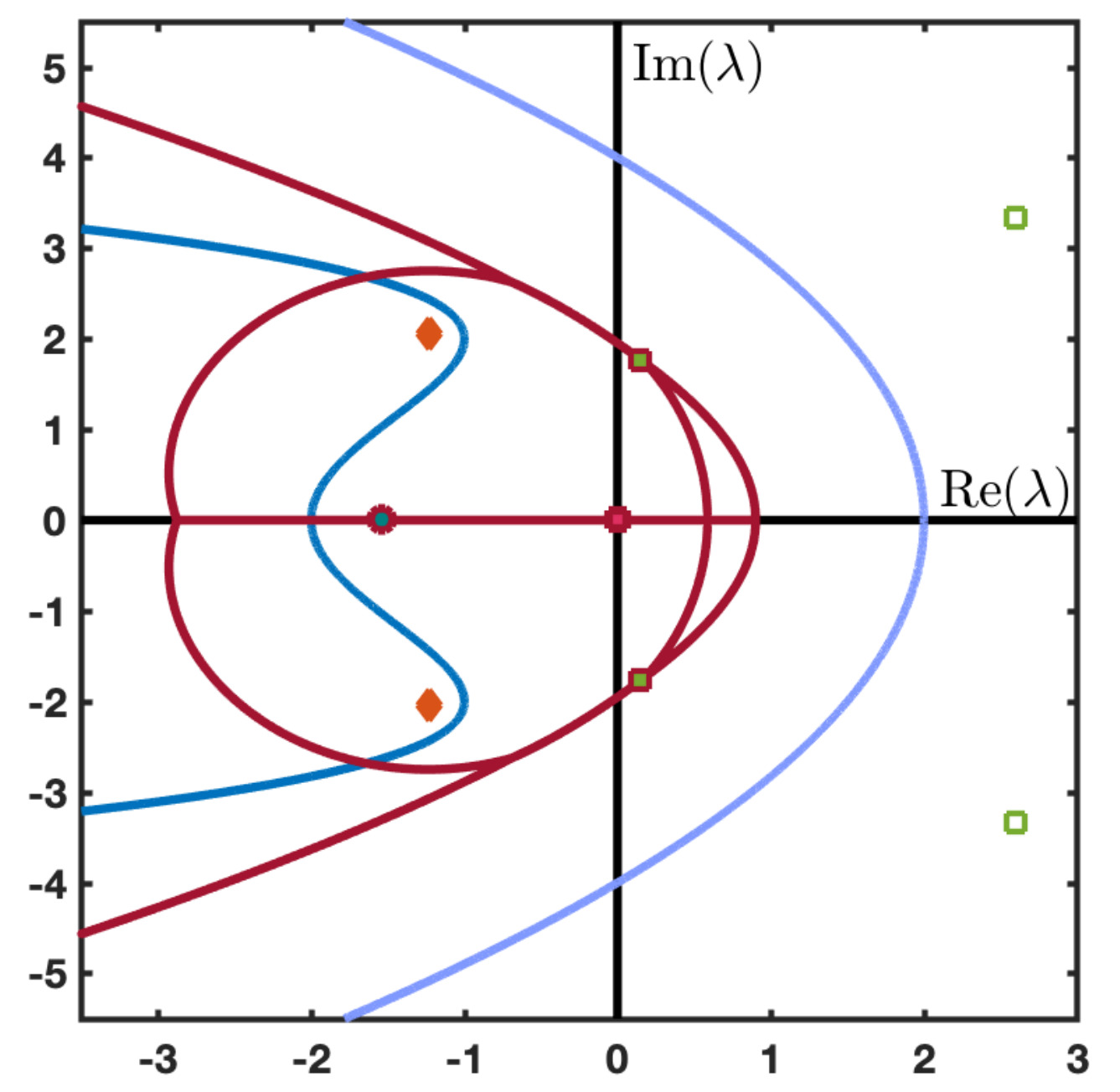}}
\caption{Spectra of $\calL^+$ at $s=s_*=2\sqrt{\alpha d}$ for various parameter values: essential spectra of $\calL^\pm_u$ (light blue) and $\calL^\pm_v$ (dark blue); absolute spectrum of $\calL^+$ (dark red) from numerical continuation; double root  $\lambda^\textnormal{dr}_1$ (pink dot),  $\lambda^\textnormal{dr}_{2,3}$, (orange diamonds), and $\lambda^\textnormal{dr}_4$ (dark cyan dot);  resonance poles $\lambda^{\textnormal{rp}_+}_\pm,\lambda^{\textnormal{rp}_-}_\pm$ (green squares, filled if pinched. In cases (a) and (b), the resonance poles in the absolute spectrum are stable, the other two are not pinched. In case (c) resonance poles do not belong to the absolute spectrum and are not pinched; in case (d), two resonance poles unstable, pinched and elements of $\Sabs(\calL^+)$, while the other two are not pinched.}
\label{fig:AbsSpec1}
\end{figure}

\begin{prop}[Lower bound on $\mathcal{R}_\mathrm{abs}$]\label{p:rabs}
 The absolute spectrum is unstable when $\mu_0^\mathrm{abs}<\mu<0$, where we recall the formula for  $\mu_0^\mathrm{abs}$ that was defined in~\eqref{e:muremabs},
 \[ \mu_0^\mathrm{abs}(\alpha,d)= \frac{d^2}{4} -\frac{4\alpha}{d} -\frac{4 \alpha^2}{d^2}.\]
\end{prop}
\begin{rmk}[Boundary between $\mathcal{R}_\mathrm{rem}$ and $\mathcal{R}_\mathrm{abs}$]
 We computed the onset of absolute instability through actual numerical computation of the absolute spectra and illustrated results in Figure \ref{fig:RegionsRj}. We found that the onset actually agrees with $\mu_0^\mathrm{abs}$, such that the lower bound in Proposition~\ref{p:rabs} is actually sharp. In fact the real triple point appears to coincide with the rightmost point of the absolute spectrum. Out of the three branches of the absolute spectrum is real and to the left of the triple point, the two others emerge vertically, with a negative tangency; see for instance Figure \ref{fig:AbsSpec1}, (a). 
\end{rmk}

\begin{prop}[Boundary of $\mathcal{R}_\mathrm{pw}$]\label{p:rpw}
 There exists a resonance pole in the unstable complex half plane precisely when $\alpha- d-d^2/2>0$ and $\mu^\mathrm{pw}(\alpha,d)<\mu<0$, where 
 \bqq
 \mu^\mathrm{pw}(\alpha,d)=\alpha - \frac{\alpha^2}{4 d^2} - \frac{d^2 (2 + d)^4}{64 \alpha^2} + 
 \frac{1}{8} (4 - 4 d - d^2). \label{eq:mupw}
 \eqq
 Moreover, when $\alpha-d-d^2/2>0$ it holds that $\mu^\mathrm{abs}_0(\alpha,d)<\mu^\mathrm{pw}(\alpha,d)$ so that $\mathcal{R}_{\mathrm{abs}}\neq \emptyset$.
\end{prop}

The remainder of this section is occupied by the proofs of these statements.

\begin{Proof}[ of Proposition \ref{p:rrem}.]
 We claim that for $\mu=\mu^\mathrm{rem}(\alpha,d)$, the absolute spectrum on the imaginary axis consists of precisely a simple double root at the origin. Continuity of the absolute spectrum \cite{rademacher07} together with the fact that the simple double root at the origin remains at the origin for all nearby parameter values then implies that the absolute spectrum possesses non-negative real part for nearby values of $\mu$, thus proving the claim. Recall that the absolute spectrum is to the left of the weighted essential spectrum, which at $\mu=\mu^\mathrm{rem}$ is to the left of the imaginary axis except for $\lambda=0$ and $\lambda^\mathrm{rem}=\sigma_v(\sqrt{3\eta_*^2+1};\eta_*))$. Since in the $\eta_*$-weight, the essential spectrum of the $v$-component is strictly to the left of the origin and stable, the double root from the $u$-equation locally gives rise to a simple curve of absolute spectrum on the negative real line, also for nearby values of $\mu$. It remains to verify that $\lambda^\mathrm{rem}$ does not belong to the absolute spectrum  for $\mu=\mu^\mathrm{rem}$. At this value of $\lambda$, the roots $\nu_j$ from the $v$-component satisfy $\Re\nu_1> \Re\nu_2=\eta_*>\Re\nu_3>\Re\nu_4,$ as a quick calculation shows. On the other hand, $\Re\nu_1>\eta_*>\Re\nu_2$ for the two roots of the $u$-component, hence excluding absolute spectrum.  
 \end{Proof}

\begin{Proof}[ of Proposition \ref{p:rabs}.]
 We label roots of the dispersion relation as $\Re\nu_+^0(\lambda)>\Re\nu_-^0(\lambda)$ for the $u$-component and $\rho_{1,2,3,4}(\lambda)$ for the $v$-component, ordered by increasing real part $\Re\rho_1(\lambda)\leq \Re\rho_2(\lambda)\leq \Re\rho_3(\lambda)\leq \Re\rho_4(\lambda)$. The proof proceeds in two main steps. 

 \textbf{Step 1: $0\in\Sabs(\calL^+)$.} We claim that for any $d, \alpha>0$, and $s=s_*$ if $\mu=\mu^\textnormal{abs}_0(\alpha,d)<0$ then 
$ \mathrm{Re}(\nu_\pm^0(0))=\mathrm{Re}(\rho_{1,2}(0))$. To see this, recall that  $\nu_\pm^0(0)=-\sqrt{\frac{\alpha}{d}}$ and $\eta_*=-\sqrt{\frac{\alpha}{d}}$.  We claim that  $D_v(0,\eta_*+\mbi k)=0$ for some $k>0$.  Since $0$ is to the right of $\Sess(\mathcal{L}_v^+)$ and $\eta_*<0$ this implies that $\eta_*\pm \mbi k$ are precisely the roots $\nu_{1,2}(0)$ and thereby yields our claim. To find $k$, we write explicitly 
\[ D_v(0,\eta_*+\mbi k)=-\eta_*^4 +6\eta_*^2 k^2-k^4 -2\eta_*^2+2k^2-1+\mu +s_*\eta_*+\mbi \left( -4\eta_*^3 k +4\eta_* k^3-4\eta_* k+s_*k\right)=0.  \] 
The imaginary part of this equation gives  $k^2=\frac{1}{4\eta_*} \left( 4\eta_*^3+4\eta_*-s_*\right)$, which when substituted into the real part, with the explicit expression for $\eta_*$, gives 
\[ -\frac{\alpha^2}{d^2}+(6\frac{\alpha}{d}+2)\left(\frac{\alpha}{d}+1+\frac{d}{2}\right)-\left(\frac{\alpha}{d}+1+\frac{d}{2}\right)^2-1+\mu-2\alpha =0,\]
which is true precisely when $\mu=\mu_0^{\textnormal{abs}}(\alpha,d)$, as claimed.

 We note that 
 \begin{equation}\label{eq:TransitionMuValue1}
\mu^\textnormal{abs}_0=\frac{d^2}{4}-\frac{4\alpha}{d}-\frac{4\alpha^2}{d^2}<0 \quad \Leftrightarrow \alpha>-d/2+d/4\sqrt{d^2+4}>0.
\end{equation}
 
\textbf{Step 2:  $\mu^\textnormal{abs}_0(\alpha,d)<\mu<0$ implies $\Sabs(\mathcal{L}^+)$ is unstable.}
The idea is to find a value $\lambda^{\textnormal{tr}}>0$ so that  $\mathrm{Re}(\nu_+^0(\lambda^{\textnormal{tr}}))=\mathrm{Re}(\rho_{1,2}(\lambda^{\textnormal{tr}}))$.  Such a triple point is automatically an element of $\Sabs(\mathcal{L}^+)$. We will locate this triple point through intersections of the weighted essential spectra for  $u$ and $v$ components, which in turn are explicit parameterized curves $\sigma_{u/v}$, respectively, 
\begin{align*}
 \sigma_u(k;\eta)&= -dk^2+2\mbi d\eta k+d\eta^2+2\mbi \sqrt{\alpha d}k + 2\sqrt{\alpha d}\eta + \alpha,\\
\sigma_v(k;\eta)&=-k^4+4\mbi\eta k^3+(2+6\eta^2)k^2+\mbi(2\sqrt{\alpha d}-4\eta-4\eta^3)k-\eta^4-2\eta^2+2\sqrt{\alpha d}\eta-1+\mu.
\end{align*}
Note that $\sigma_u(0,\eta)=(\sqrt{d} \eta+\sqrt{\alpha})^2$ such that $\nu_+^0((\sqrt{d} \eta+\sqrt{\alpha})^2)=\eta$. We therefore consider $\eta_*<\eta<0$.  Moreover,
\[\Im(\sigma_v(k;\eta))=0\quad \Longleftrightarrow \quad k_v^0=0, k_v^\pm=\pm\sqrt{\eta^2+1-\frac{\sqrt{\alpha d}}{2\eta}}.\]
We disregard the zero solution, since this occurs for $\lambda$ in the left half of the complex plane. 
We now obtain that there exists a triple point, if there exists a (real) weight $\eta_*<\eta<0$ (recall $s=s_*>0$) such that
$$\Psi(\eta)\coloneqq\sigma_v(k_v^+;\eta)-\sigma_u(0;\eta)=4\eta^4+(4-d)\eta^2-2\sqrt{\alpha d}\eta+\mu-\alpha-\frac{\alpha d}{4\eta^2}=0.$$
We see that  $\Psi(\eta)$ is continuous for all $\eta<0$ with $\lim_{\eta\to 0^-} \Psi(\eta)=-\infty$.  We also see that 
\[ \Psi(\eta_*)=-\mu^\textnormal{abs}_0(\alpha,d)+\mu\] 
and is therefore positive for $\mu^\textnormal{abs}_0(\alpha,d)<\mu<0$.  Finally, since $\Psi'(\eta)<0$ for all $\eta_*<\eta<0$ we obtain the existence of a unique $\eta^{\textnormal{tr}}$ for which $\Psi(\eta^{\textnormal{tr}})=0$ which in turn establishes the existence of $\lambda^{\textnormal{tr}}=\sigma_u(0,\eta^{\textnormal{tr}})>0$, which concludes the proof. 
\end{Proof}

\begin{Proof}[ of Proposition \ref{p:rpw}.]
One can find all double roots explicitly, and we list the results below after setting $s=s_*$:
\begin{align}
\lambda_1^\textnormal{dr}&=0,\label{e:udr}\\
 \lambda_{2,3}^\textnormal{dr}&=\mu-\frac{1}{3}\pm\frac{(\mbi\sqrt{3}\mp 1)(27\alpha d-2)}{3C_3}+\frac{1}{24}(1\pm\mbi\sqrt{3})C_3,\label{eq:vbranch}\\
\lambda_4^\textnormal{dr}&=\mu-\frac{1}{3}+\frac{2(27\alpha d-2)}{3C_3}-\frac{1}{12}C_3,\\
\lambda^{\textnormal{rp}_-}_\pm&= \alpha-d-\frac{d^2}{2}\pm dC_4-2\sqrt{\alpha d(\pm C_4-1-d/2)},\\
\lambda^{\textnormal{rp}_+}_\pm&= \alpha-d-\frac{d^2}{2}\pm dC_4+2\sqrt{\alpha d(\pm C_4-1-d/2)}.
\end{align}
where
\[
C_3\coloneqq \left(64+54\alpha d(40-27\alpha d)+6\sqrt{3}\sqrt{\alpha d(16+27\alpha d)^3} \right)^{1/3}>0,\qquad C_4\coloneqq \sqrt{\mu+d-\alpha+d^2/4}.\]
Double roots to the $u$- and the $v$-equation, \eqref{e:udr} and \eqref{eq:vbranch}, respectively, are always stable with $\mu<0$. This is clear for $\lambda_1^\textnormal{dr}$ and for $\lambda^\textnormal{dr}_{2,3}$ (see\ Lemma~\ref{l:shabs}). The roots $\lambda_{2,3,4}^\mathrm{dr}$ are induced by the $v$-equation, alone. Since pinched double roots necessarily lie to the left of the essential spectrum, these double root are either located in the stable complex half plane or are not pinched. Since $\lambda_{2,3}^\mathrm{dr}$ are complex conjugate and involve all 4 roots of the $v$-equation with equal real part, they are in fact pinched and lie to the left of the essential spectrum.

The double roots $\lambda^{\textnormal{rp}_\pm}_\pm$ are in fact resonance poles, involving a root from the SH and one from the KPP subsystem. We first focus on $\lambda^{\textnormal{rp}_-}_\pm$. We distinguish two cases depending on the sign of  $\mu-\alpha+d+d^2/4$.

First, assume that  $\mu-\alpha+d+d^2/4>0$. Then 
\bqs
0<C_4 = \sqrt{(1+d/2)^2+\mu-\alpha-1}<1+d/2,
\eqs
such that
\bqs
\Re\left(\lambda^{\textnormal{rp}_-}_- \right)\leq \Re\left(\lambda^{\textnormal{rp}_-}_+ \right)=\mu-(C_4-d/2)^2<0.
\eqs

Next, suppose that $\mu-\alpha+d+d^2/4<0$, then 
$C_4 = \mbi \sqrt{-(\mu-\alpha+d+d^2/4)}$,
and thus
\begin{align}
\Re\left(\lambda^{\textnormal{rp}_-}_\pm \right) &= \alpha-d-\frac{d^2}{2}-2\sqrt{\alpha d}\,\Re\left( \sqrt{-1-d/2\pm \mbi \sqrt{-(\mu-\alpha+d+d^2/4)}}\right) \nonumber \\
&=\alpha-d-\frac{d^2}{2}-\sqrt{2\alpha d}\sqrt{-1-d/2+\sqrt{1-\mu +\alpha}}, \label{eq:rplambda}
\end{align}
after a short calculation. We claim that $\left(\nu^{\textnormal{rp}_-}_\pm,\lambda^{\textnormal{rp}_-}_\pm\right)$ is pinched whenever $\Re\left(\lambda^{\textnormal{rp}_-}_\pm \right)\geq0$, where, explicitly,
\bqs
\nu^{\textnormal{rp}_-}_\pm=-\sqrt{-1-d/2\pm\sqrt{C_4}}=-\sqrt{-1-d/2\pm\mbi\sqrt{-(\mu-\alpha+d+d^2/4)}}.
\eqs
Clearly, $\lambda^{\textnormal{rp}_-}_\pm$ is to the right of the essential spectrum in the $v$-component and 
\bqs
\Re\left(\nu^{\textnormal{rp}_-}_\pm \right)=-\frac{\sqrt{2}}{2}\sqrt{-1-d/2+\sqrt{1-\mu +\alpha}}<0.
\eqs
As a consequence, the double root $\left(\nu^{\textnormal{rp}_-}_\pm,\lambda^{\textnormal{rp}_-}_\pm\right)$ is pinched in case we can prove that
$\nu^{\textnormal{rp}_-}_\pm=\nu^0_+(\lambda^{\textnormal{rp}_-}_\pm).
$ For this, it suffices to verify
\bqs
\Re\left(\nu^{\textnormal{rp}_-}_\pm \right)>\eta_*=-\sqrt{\frac{\alpha}{d}},
\eqs
or, after some simplifications, 
\bqs
- \mu < \frac{4\alpha^2}{d^2}+\frac{4\alpha}{d}+\frac{d^2}{4}+\alpha+d=-\mu_0^\textnormal{abs}  +\frac{d^2}{2}+\alpha+d.
\eqs
This last inequality is a consequence of $\mu_0^\textnormal{abs}<\mu<0$ and therefore  $\left(\lambda^{\textnormal{rp}_-}_\pm \right)$ is pinched when located in the right half plane.  To determine conditions leading to instability of this resonance pole, we consult (\ref{eq:rplambda}) and note that $\alpha-d-d^2/2$ must be positive.  If this condition holds, then we may compute that the resonance pole becomes unstable at  $\mu^{\mathrm{pw}}(\alpha,d)$ given in (\ref{eq:mupw}).  It is a short calculation to verify that $\alpha-d-d^2/2>0$ implies that $\mu_0^{\mathrm{abs}}(\alpha,d)<\mu^{\mathrm{pw}}(\alpha,d)$.

Next, we turn to $\lambda^{\textnormal{rp}_+}_\pm$,% when $ \mu_0^\mathrm{abs}<\mu<0$?
starting with the case  $\mu-\alpha+d+d^2/4>0$. Again, $C_4>0$ and 
\bqs
\Re\left(\lambda^{\textnormal{rp}_+}_- \right)\leq \Re\left(\lambda^{\textnormal{rp}_+}_+ \right)=\mu-(C_4-d/2)^2<0.
\eqs

In case $\mu-\alpha+d+d^2/4<0$, we find
\bqs
\Re\left(\lambda^{\textnormal{rp}_+}_\pm \right) =\alpha-d-\frac{d^2}{2}+\sqrt{2\alpha d}\sqrt{-1-d/2+\sqrt{1-\mu +\alpha}}.
\eqs
Note that $\left(\nu^{\textnormal{rp}_+}_\pm,\lambda^{\textnormal{rp}_+}_\pm\right)$ is not pinched whenever $\Re\left(\lambda^{\textnormal{rp}_+}_\pm \right)\geq0$, where, explicitly, 
\bqs
\nu^{\textnormal{rp}_+}_\pm=\sqrt{-1-d/2\pm\sqrt{C_4}}=\sqrt{-1-d/2\pm\mbi\sqrt{-(\mu-\alpha+d+d^2/4)}},\qquad 
\Re\left(\nu^{\textnormal{rp}_+}_\pm \right)>0.
\eqs
Since $\Re\left(\lambda^{\textnormal{rp}_+}_\pm \right)\geq0$, then $\lambda^{\textnormal{rp}_+}_\pm$ is to the right of the $v$-essential spectrum and $\Re\left(\nu^{\textnormal{rp}_+}_\pm \right)>0$ ensures that it cannot be pinched.\end{Proof}

\section{Numerical simulations and the appearance of faster invasion modes}\label{sec:Num}

In Theorem~\ref{thmMain}, we have established the stability of the critical Fisher-KPP front in the presence of inhomogeneous coupling to a secondary equation (the linearized Swift-Hohenberg equation).  This result holds for parameters in $\mathcal{R}_{\mathrm{rem}}$ and $\mathcal{R}_{\mathrm{abs}}$ in the case when the system has a remnant instability but lacks unstable resonance poles (with the caveat of an additional condition excluding possible resonances).  The main takeaway of this result is that the inability to stabilize essential spectrum using exponential weights does not preclude stability of the traveling front.

In this final section, we investigate front stability and speed selection in numerical simulations, distinguishing in particular between regions $\mathcal{R}_{\mathrm{rem}}$ and $\mathcal{R}_{\mathrm{abs}}$. For parameters in $\mathcal{R}_{\mathrm{rem}}$, the absolute spectrum is stable and we observe stability of the front in numerical simulations. For parameters in $\mathcal{R}_\mathrm{abs}$, we observe stability of the front over large periods of time, thereby corroborating our analytical result.  Eventually, however, the front appears to be overtaken by an incoherent invasion mode traveling with a faster, approximately constant, average speed.  Below, we explain this eventual acceleration and offer a prediction faster ensuing speed.  In summary, the faster invasion speed is mediated by a resonance caused by small couplings between modes due to numerical  round off errors. The strongest such resonance enables propagation at a speed given by marginal stability of the absolute spectrum of the zero state rather than marginal stability of resonance poles as identified in our result. Numerical evidence supports our prediction of invasion at this faster speed which we will refer to as the {\em absolute spreading speed},
\begin{equation}\label{e:sabs}
 s_\mathrm{abs}\coloneqq\sup\{s|\, \Sigma_\mathrm{abs} \text{ is unstable}\}.
\end{equation}

\paragraph{Numerical simulations ---   confirmation of stability and absolute spreading speeds.}

\begin{figure}[ht!]
\centering
\subfigure{\includegraphics[width=0.215\textwidth]{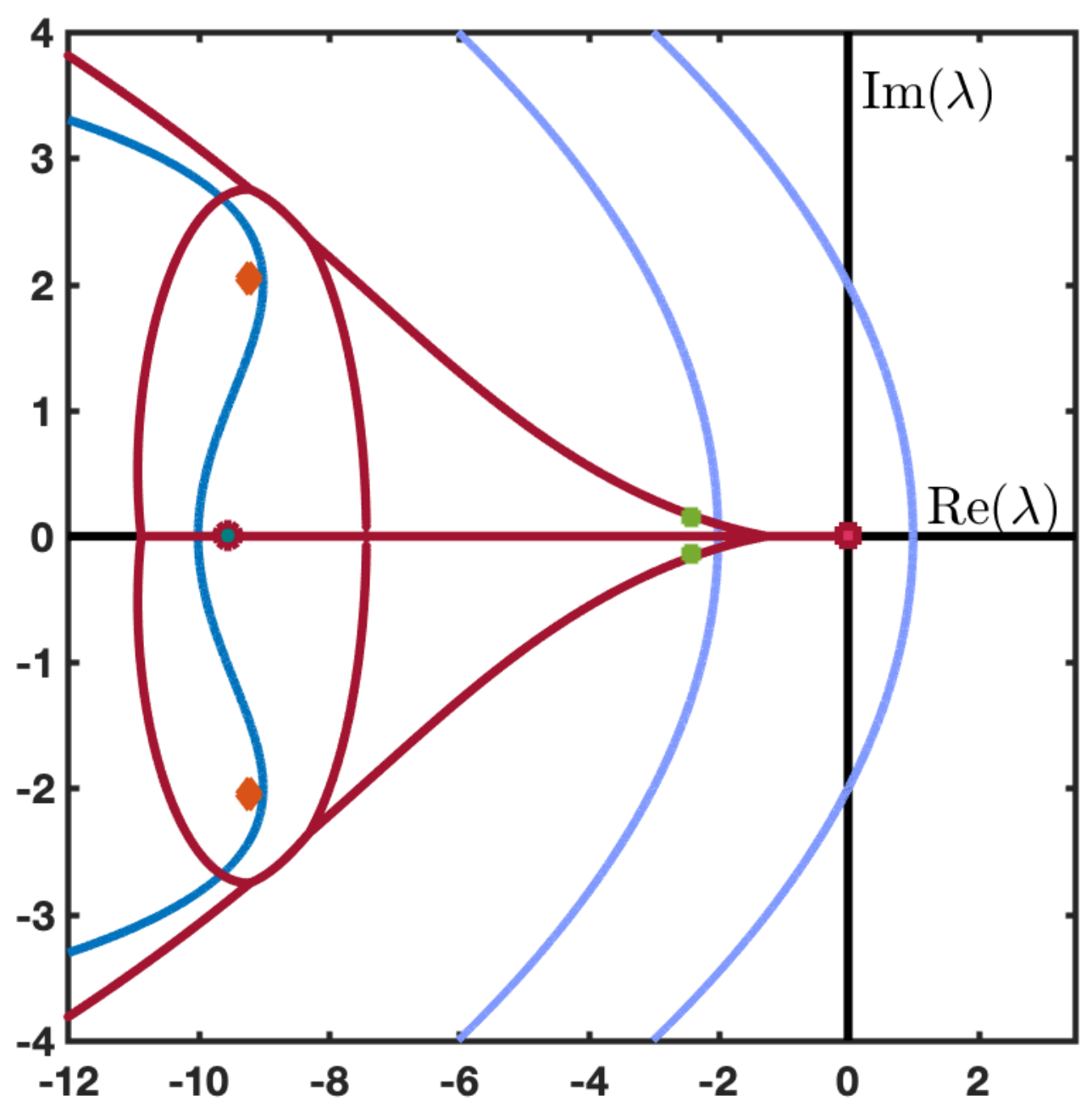}}
\qquad
\subfigure{\includegraphics[width=0.222\textwidth]{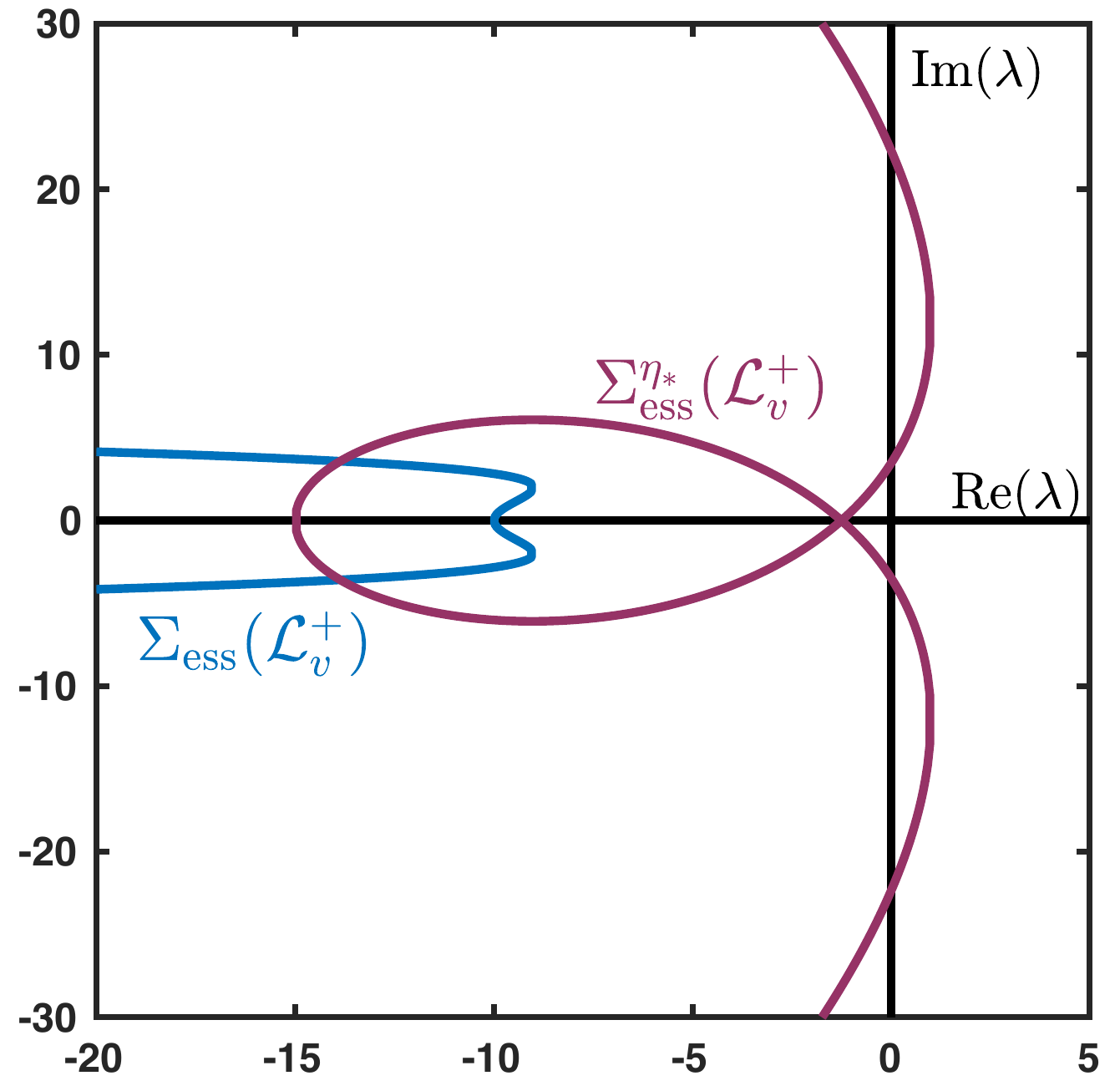}}\qquad
\subfigure{\includegraphics[width=0.25\textwidth]{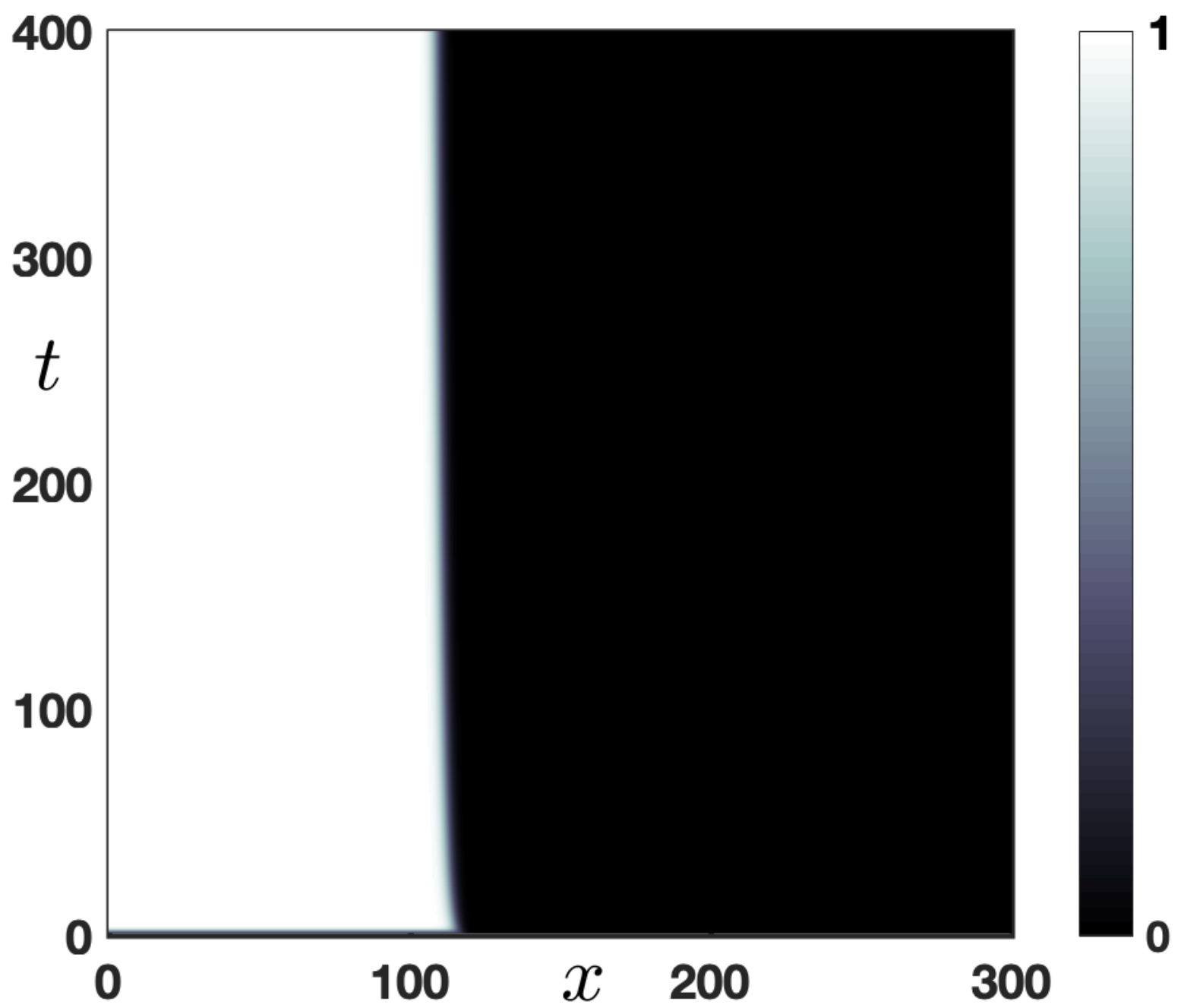}}
\caption{Essential and absolute spectra for $d=\alpha=1$, $\mu=-9$, showing stability of the absolute spectrum (left, see Figure \ref{fig:AbsSpec1} for explanation); instability of the  weighted $v$-essential spectrum with weight $\eta_*$, showing a remnant instability  (center) and space time plot (right), demonstrating front stability in  $\mathcal{R}_{\mathrm{rem}}$}
\label{fig:notremnant}
\end{figure}

\begin{figure}[h!]
\centering
\includegraphics[width=0.16\textwidth]{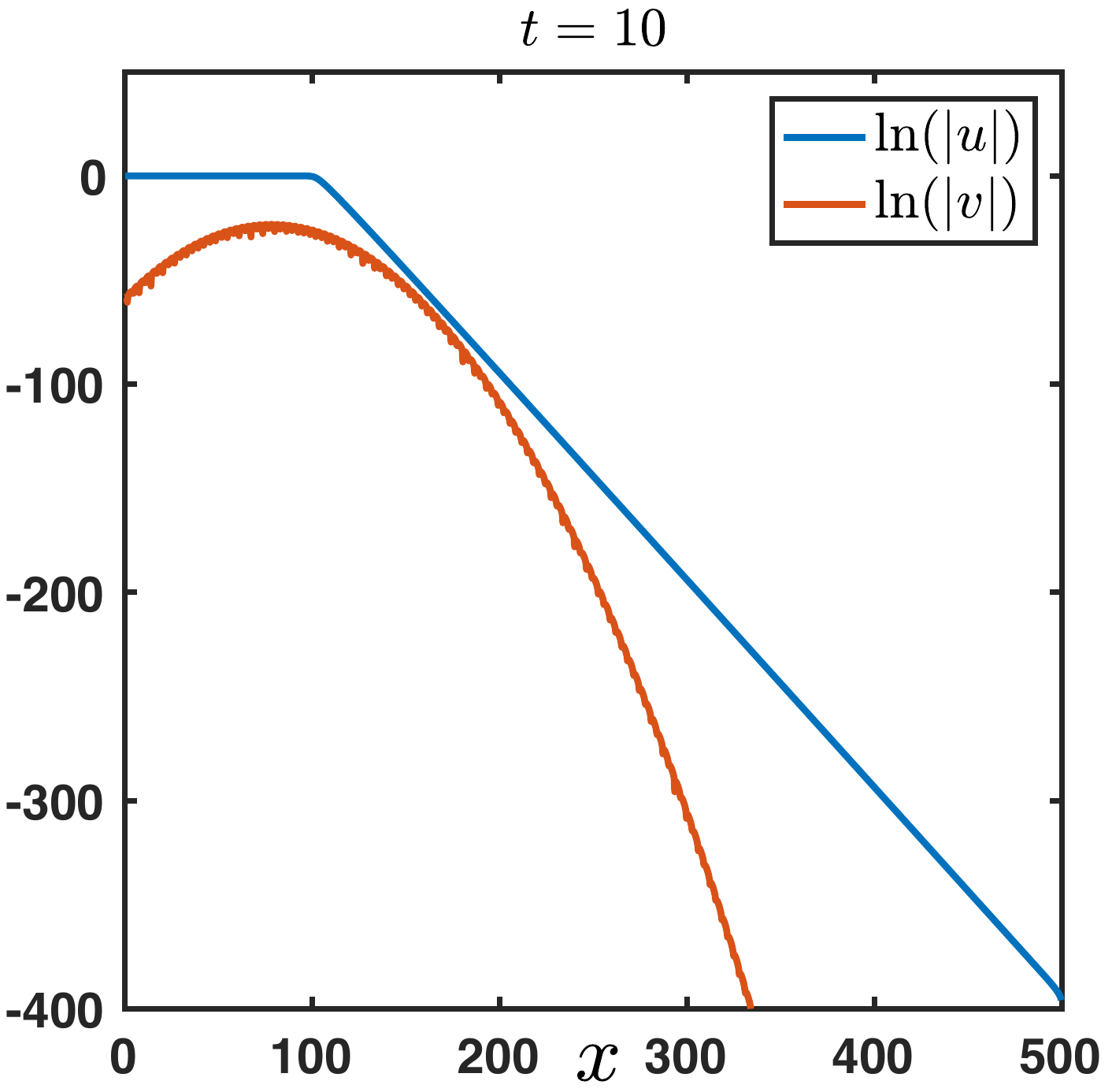}
\includegraphics[width=0.16\textwidth]{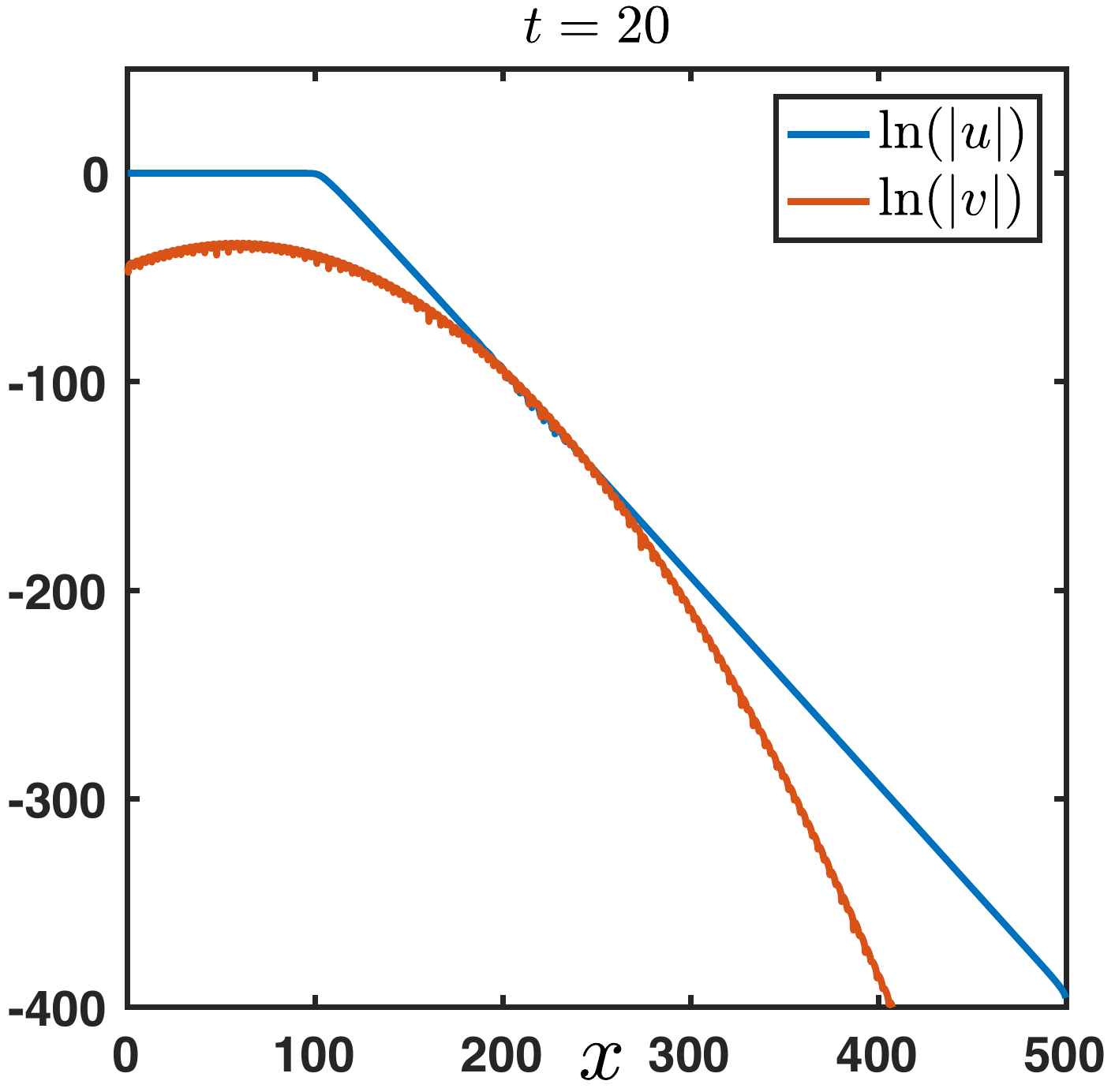}
\includegraphics[width=0.16\textwidth]{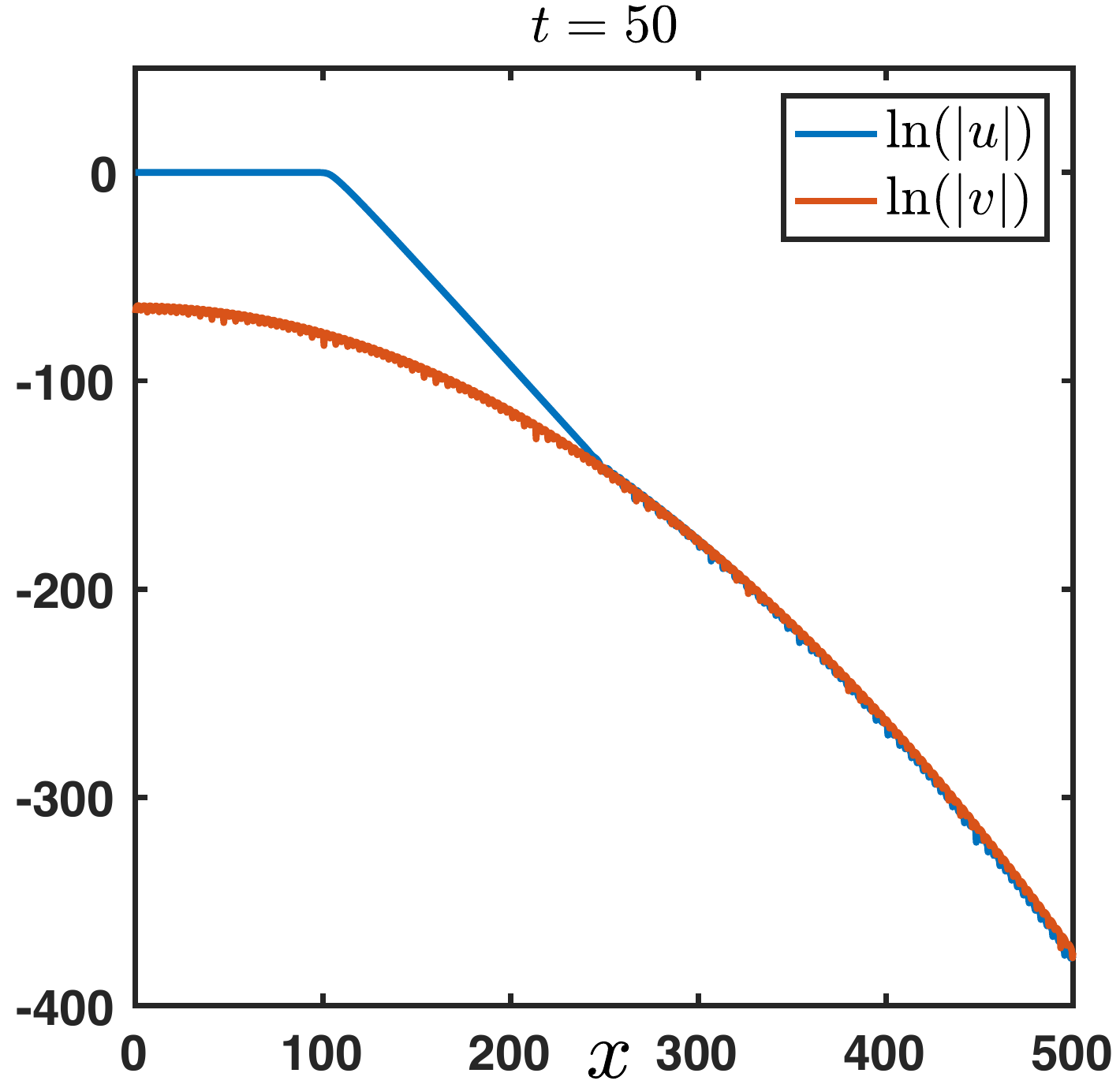}
\includegraphics[width=0.16\textwidth]{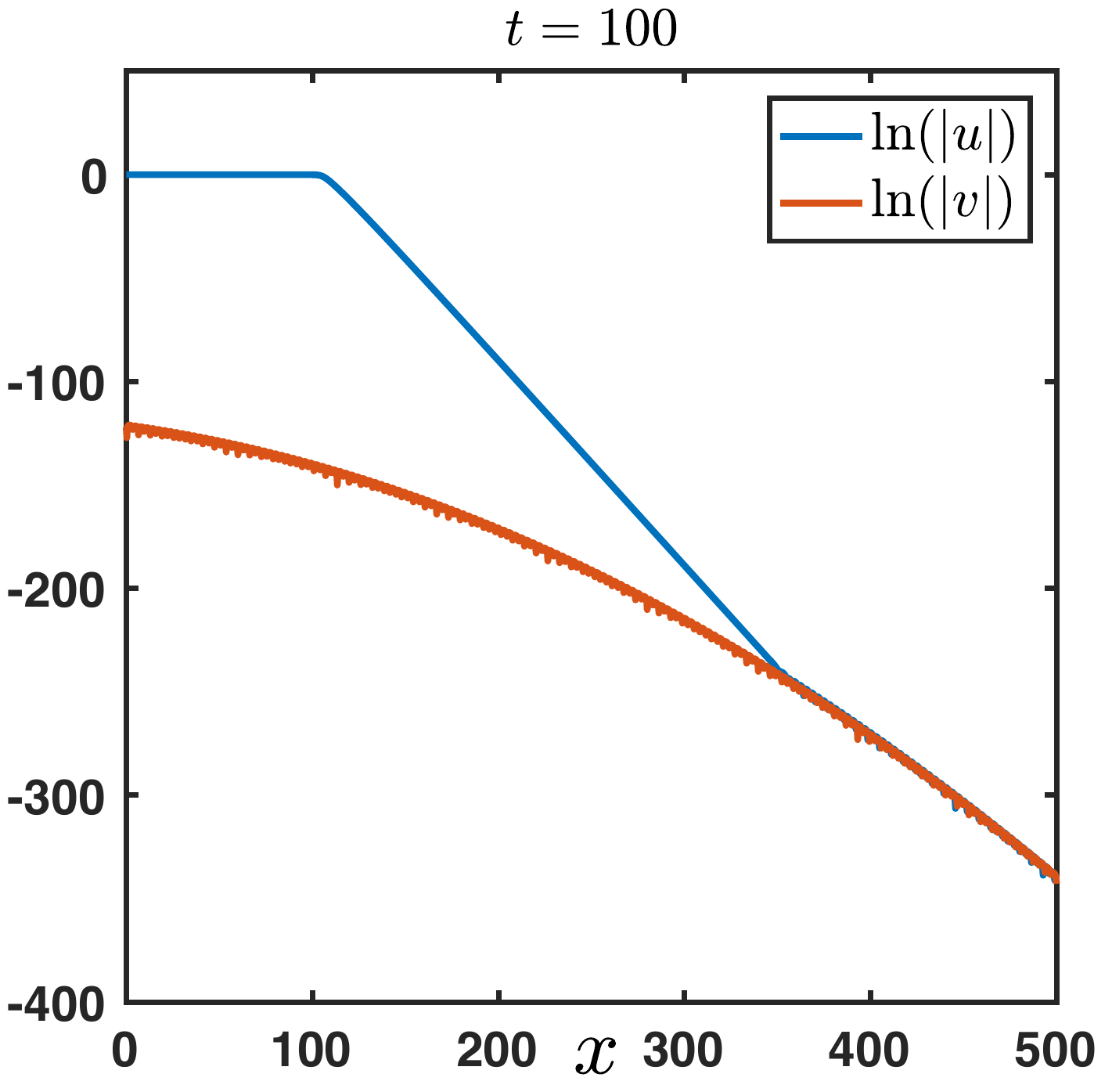}
\includegraphics[width=0.16\textwidth]{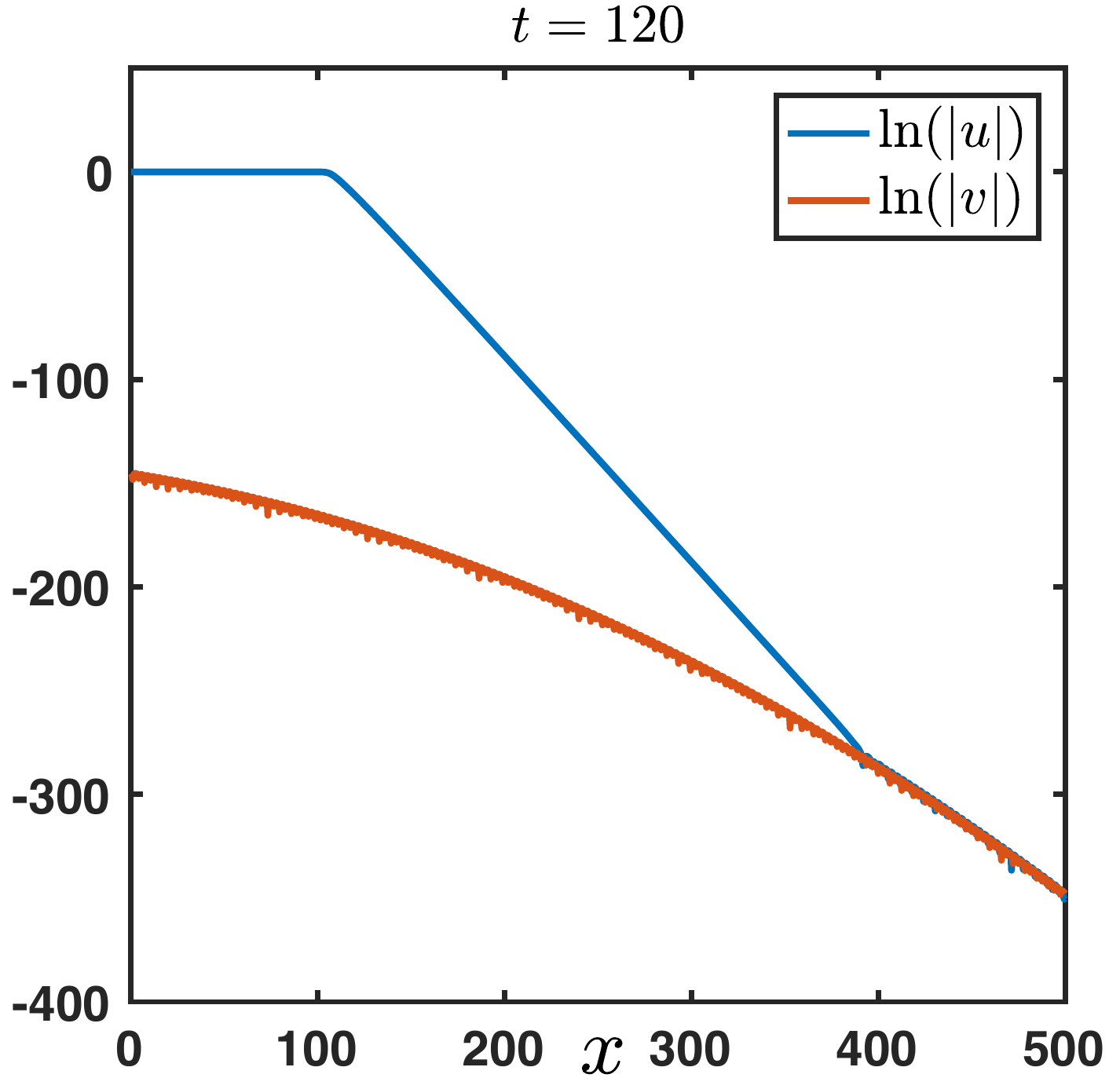}
\includegraphics[width=0.16\textwidth]{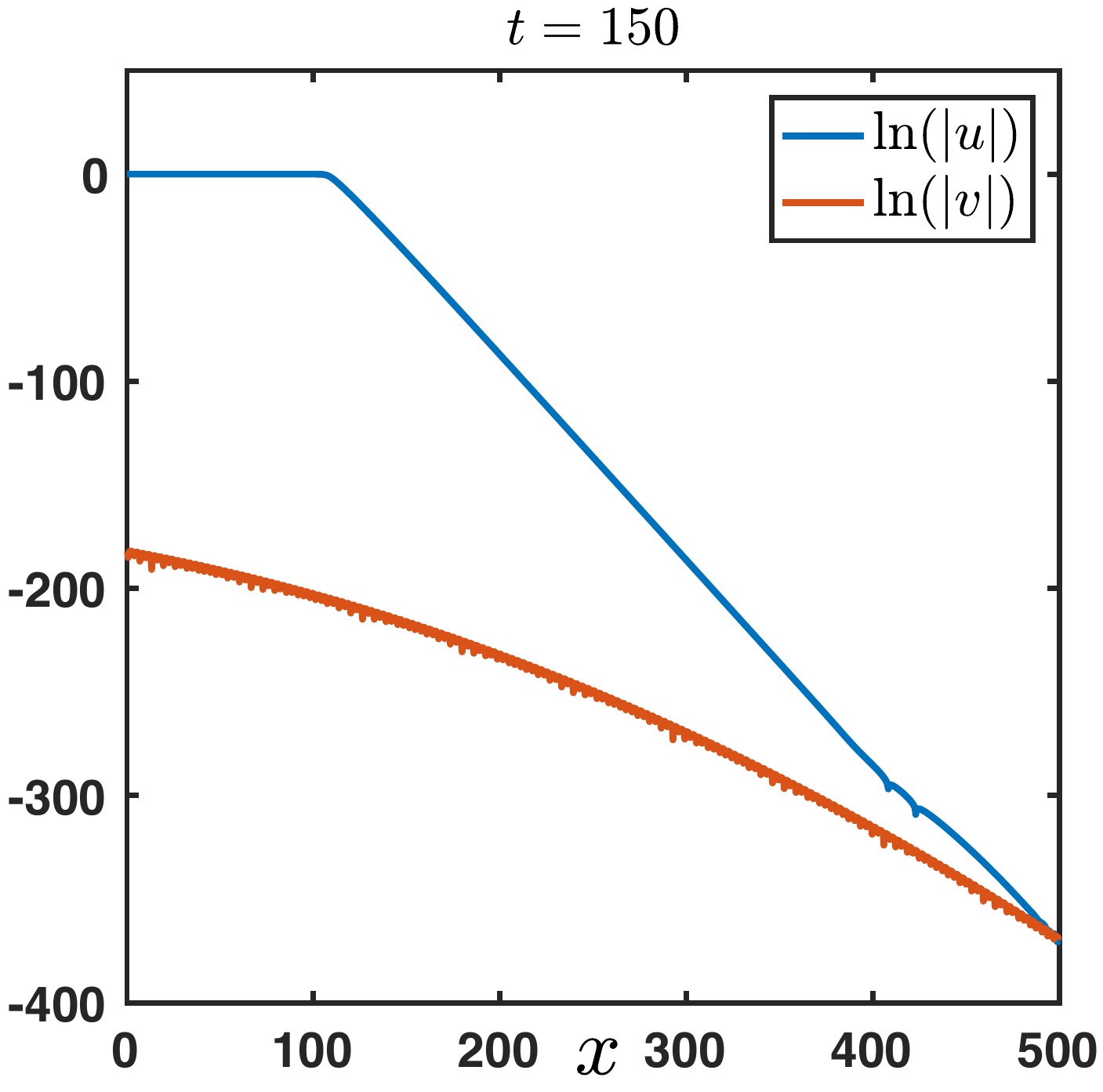}
\caption{Snapshots of the logarithms of $u$- and $v$-components at $t\in\left\{10,20,50,100,120,150 \right\}$ with $(d,\alpha,\mu)=(1,1,-1)$ and $\beta=10^{-4}$. The envelope of the $v$-equation follows a parabola, $v(t,x)=\frac{\me^{\mu t}}{4\sqrt{\pi t}} \cos(x+s_*t) \me^{-\frac{(x+s_*t)^2}{16t}}$, and induces weak decay in the $u$-equation, which however is not supported in the absence of the $v$-component.}
\label{fig:Logplotsolutions}
\end{figure}

\begin{figure}[ht!]
\centering
\subfigure{\includegraphics[width=0.2\textwidth]{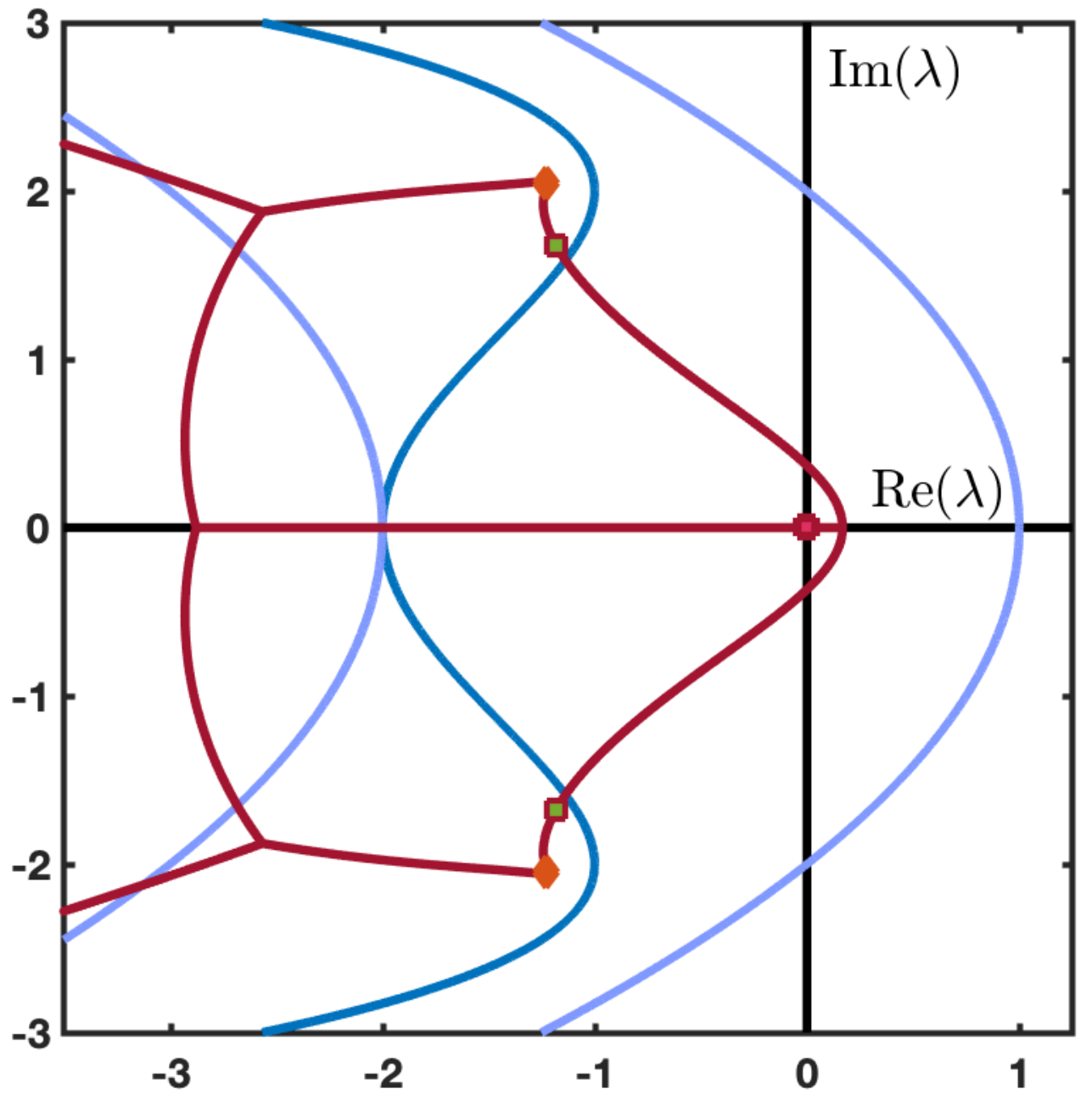}}\qquad
\subfigure{\includegraphics[width=0.2\textwidth]{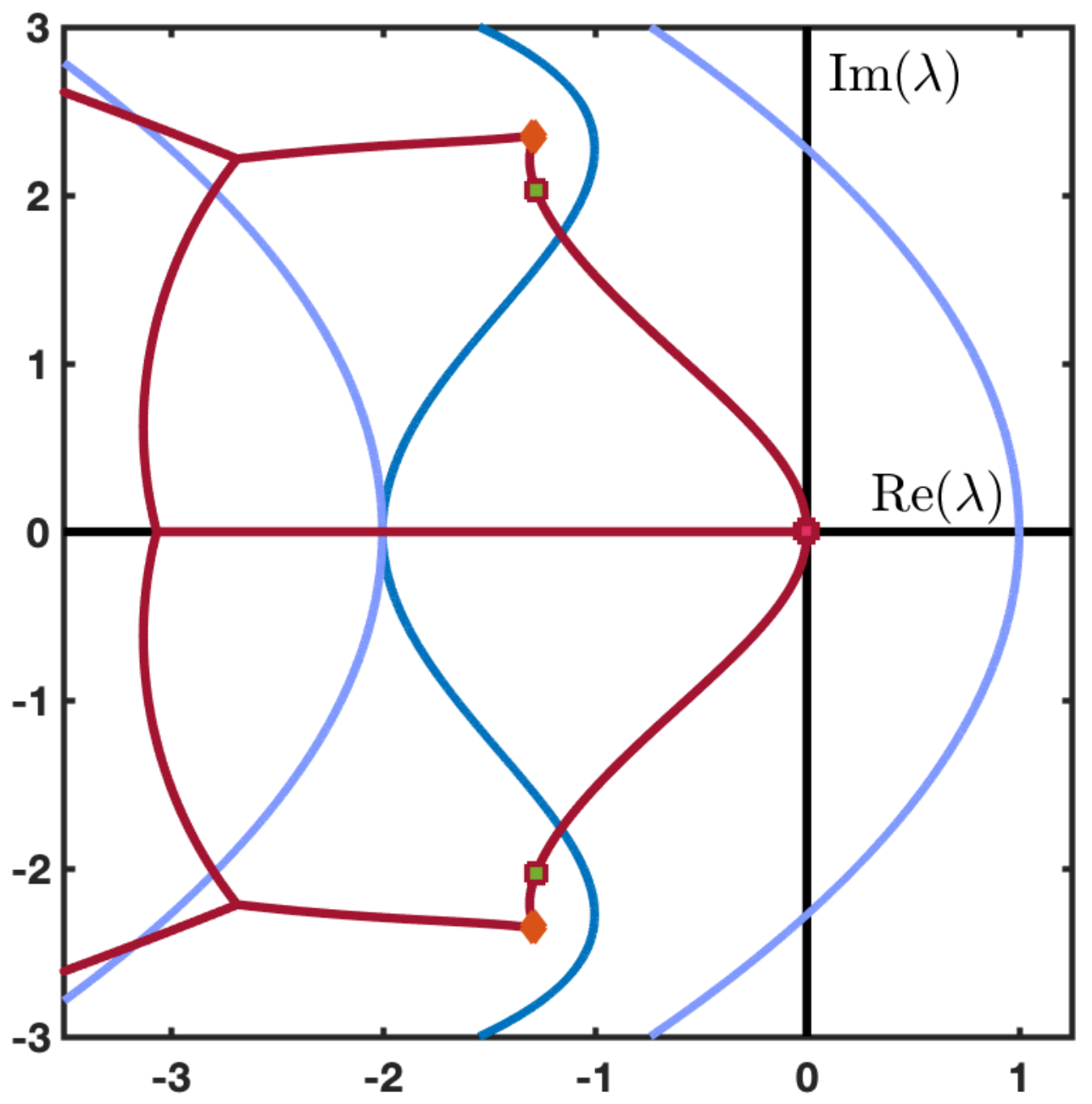}}\qquad
\subfigure{\includegraphics[width=0.243\textwidth]{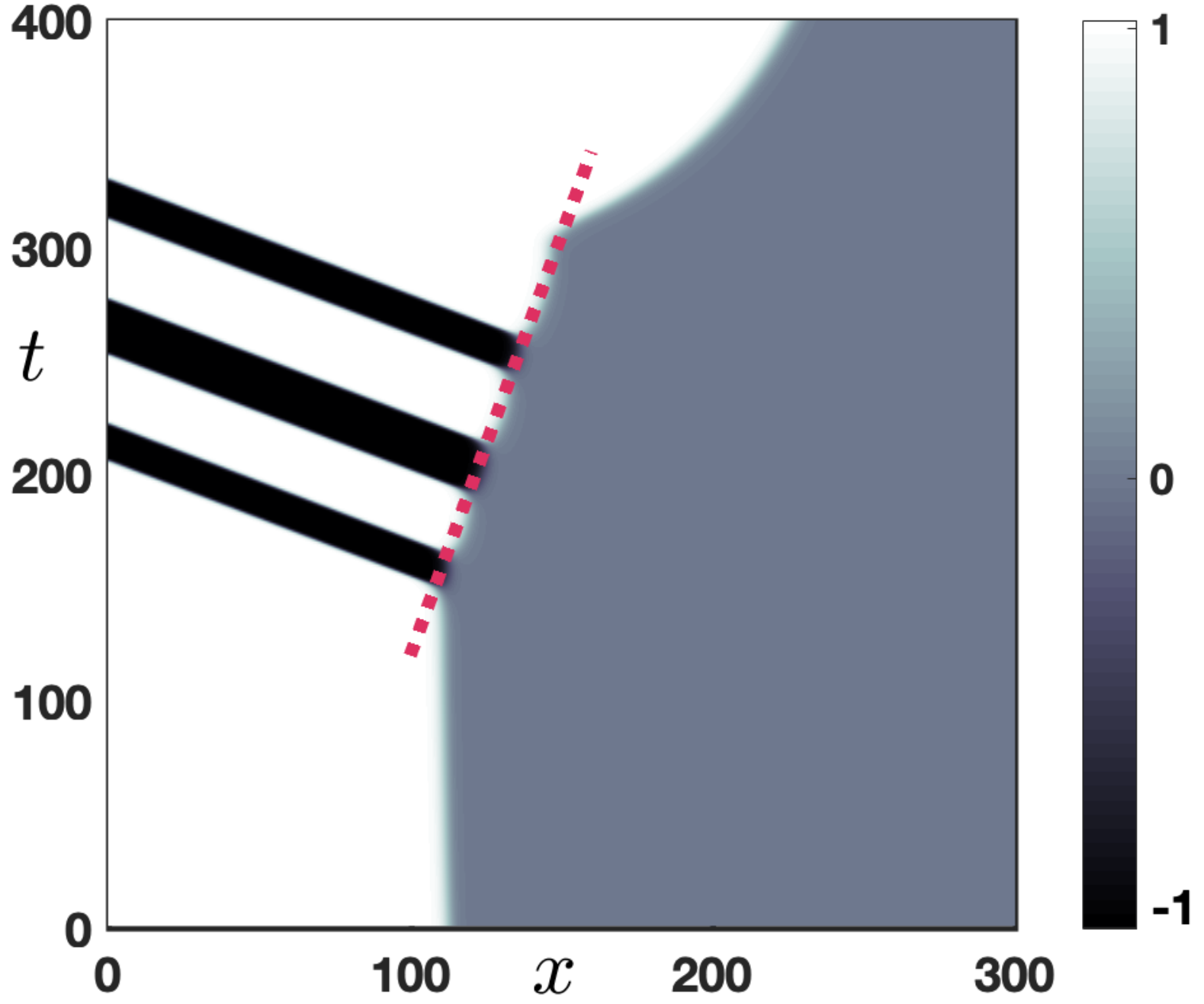}}\\
\subfigure{\includegraphics[width=0.2\textwidth]{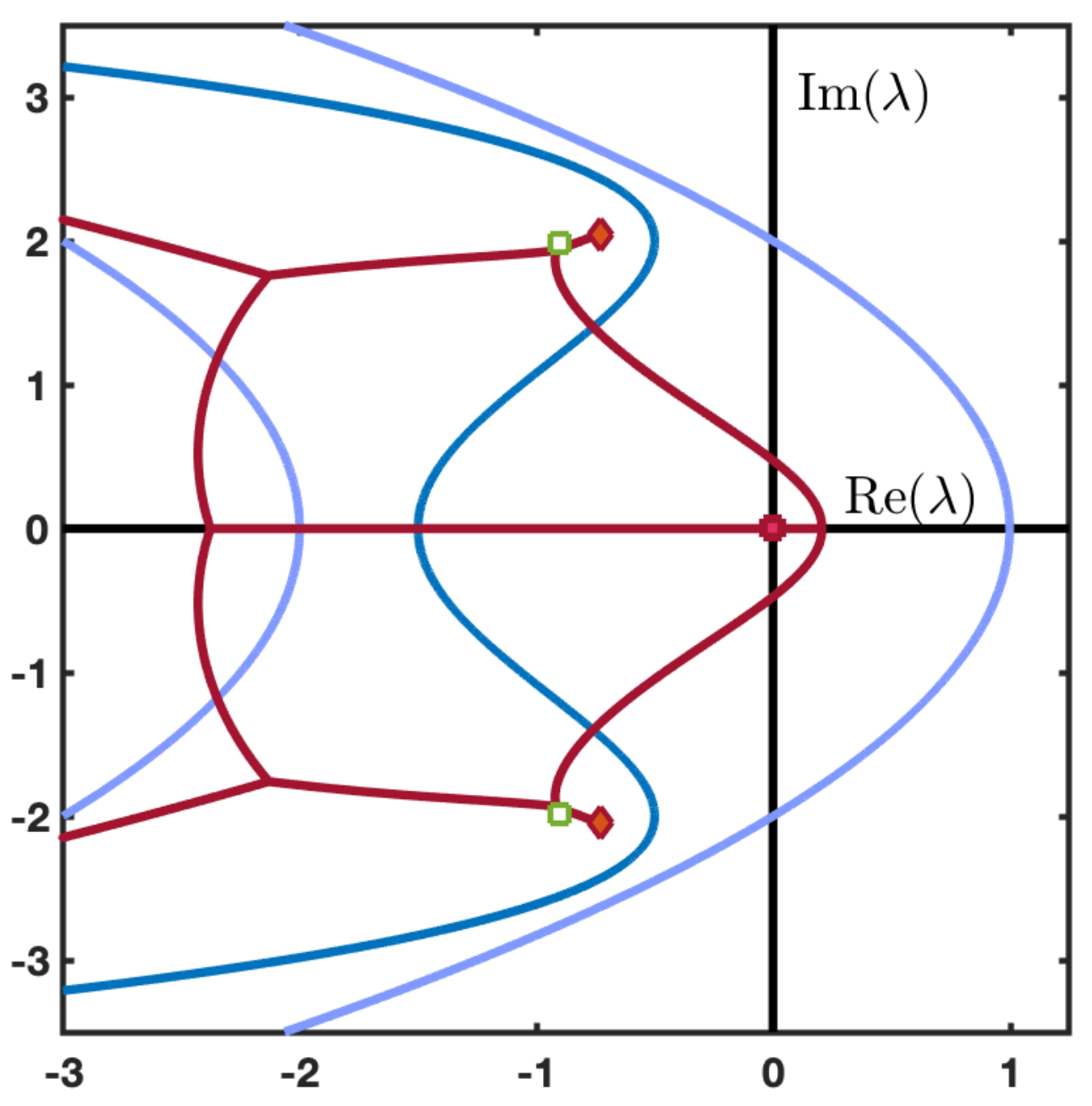}}\qquad
\subfigure{\includegraphics[width=0.2\textwidth]{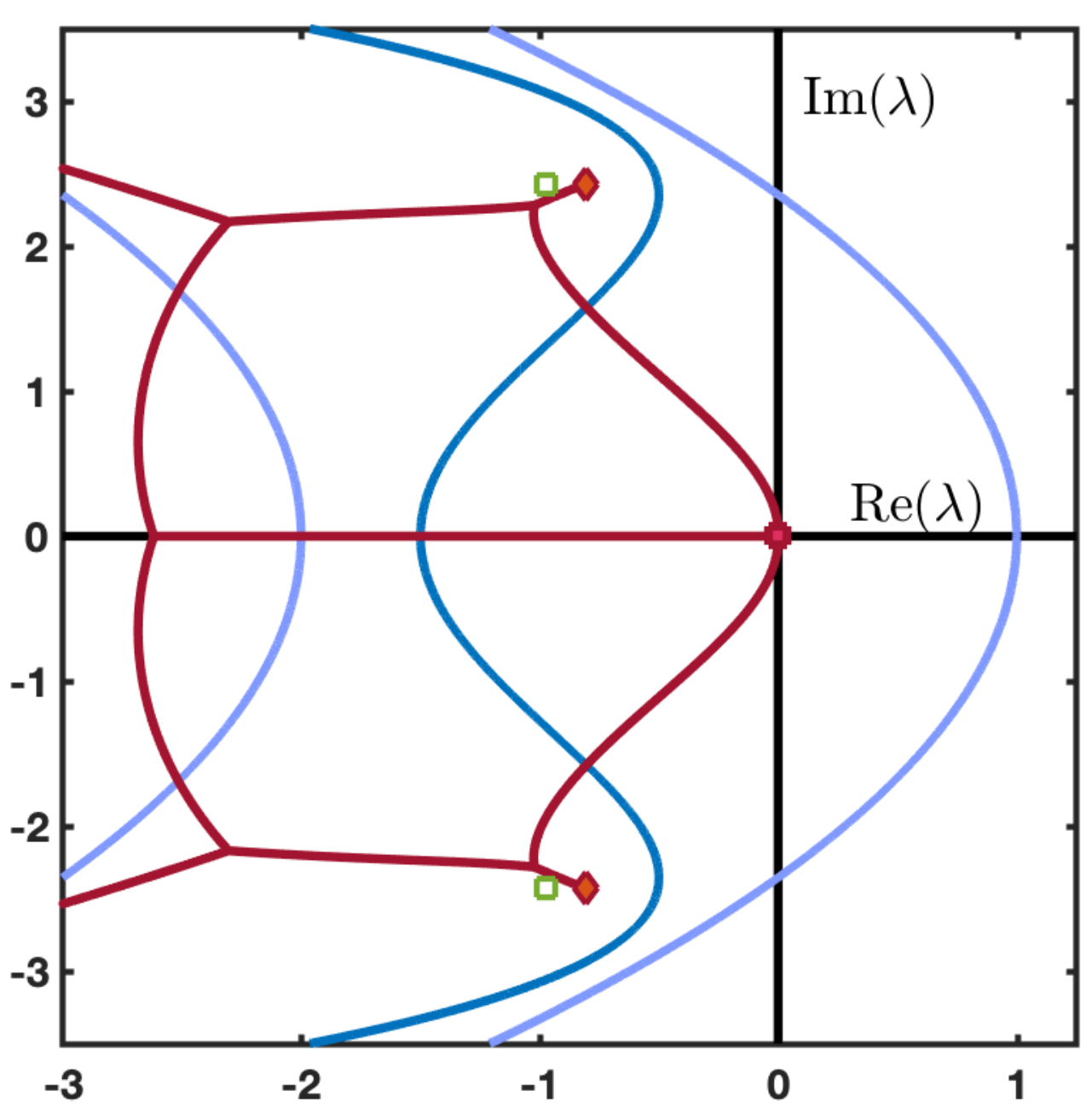}}\qquad
\subfigure{\includegraphics[width=0.243\textwidth]{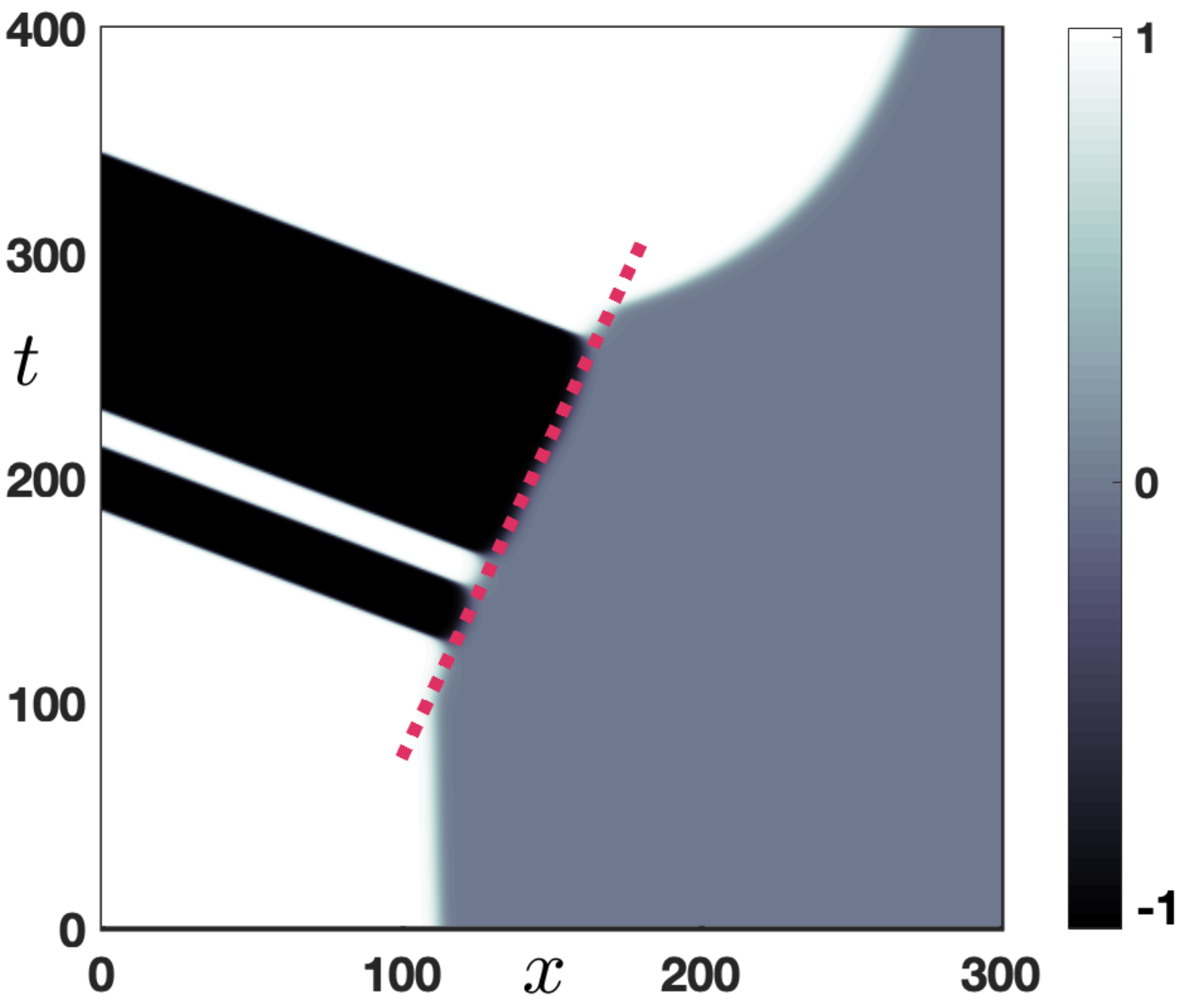}}\\
\subfigure{\includegraphics[width=0.2\textwidth]{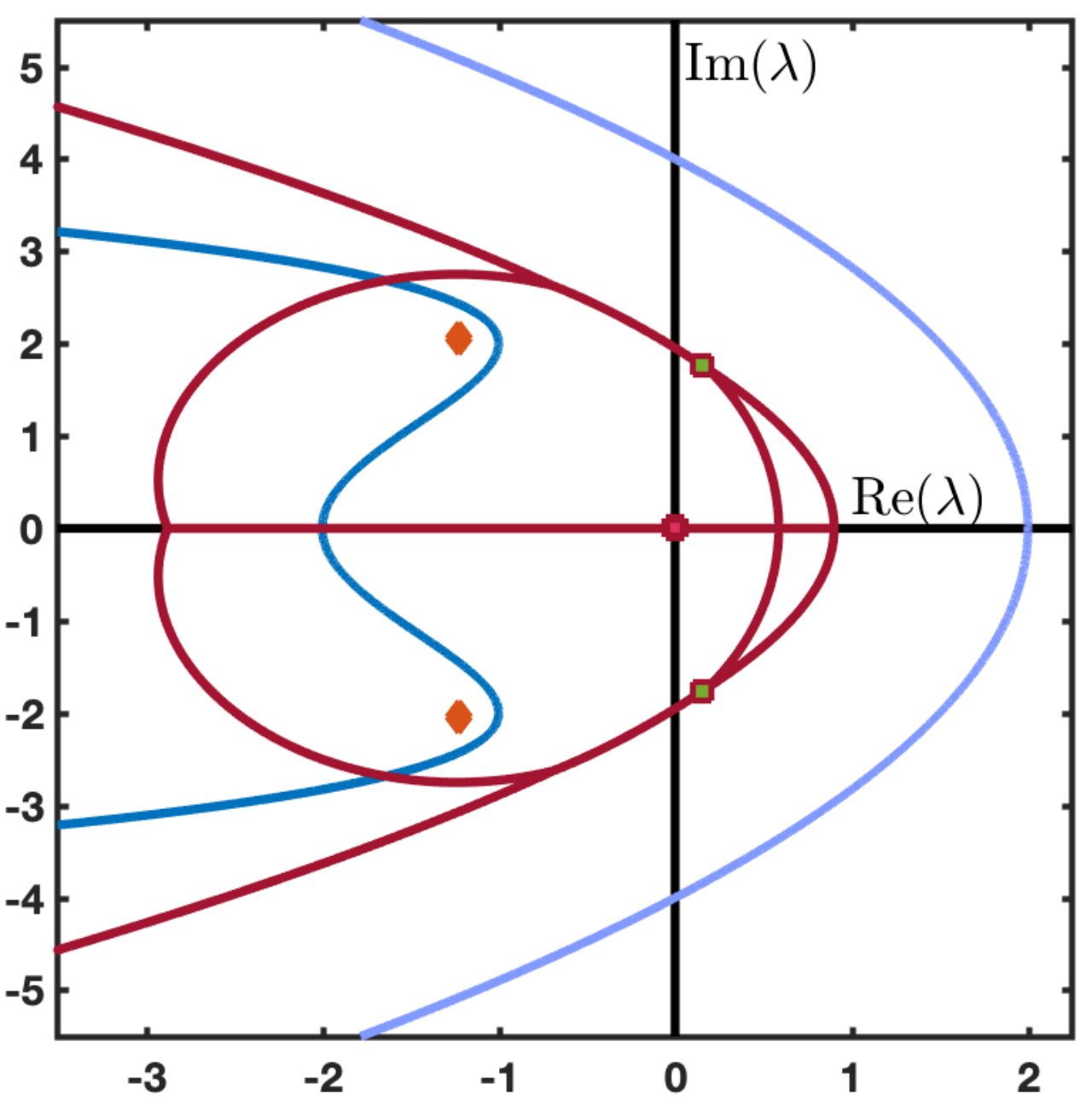}}\qquad
\subfigure{\includegraphics[width=0.2\textwidth]{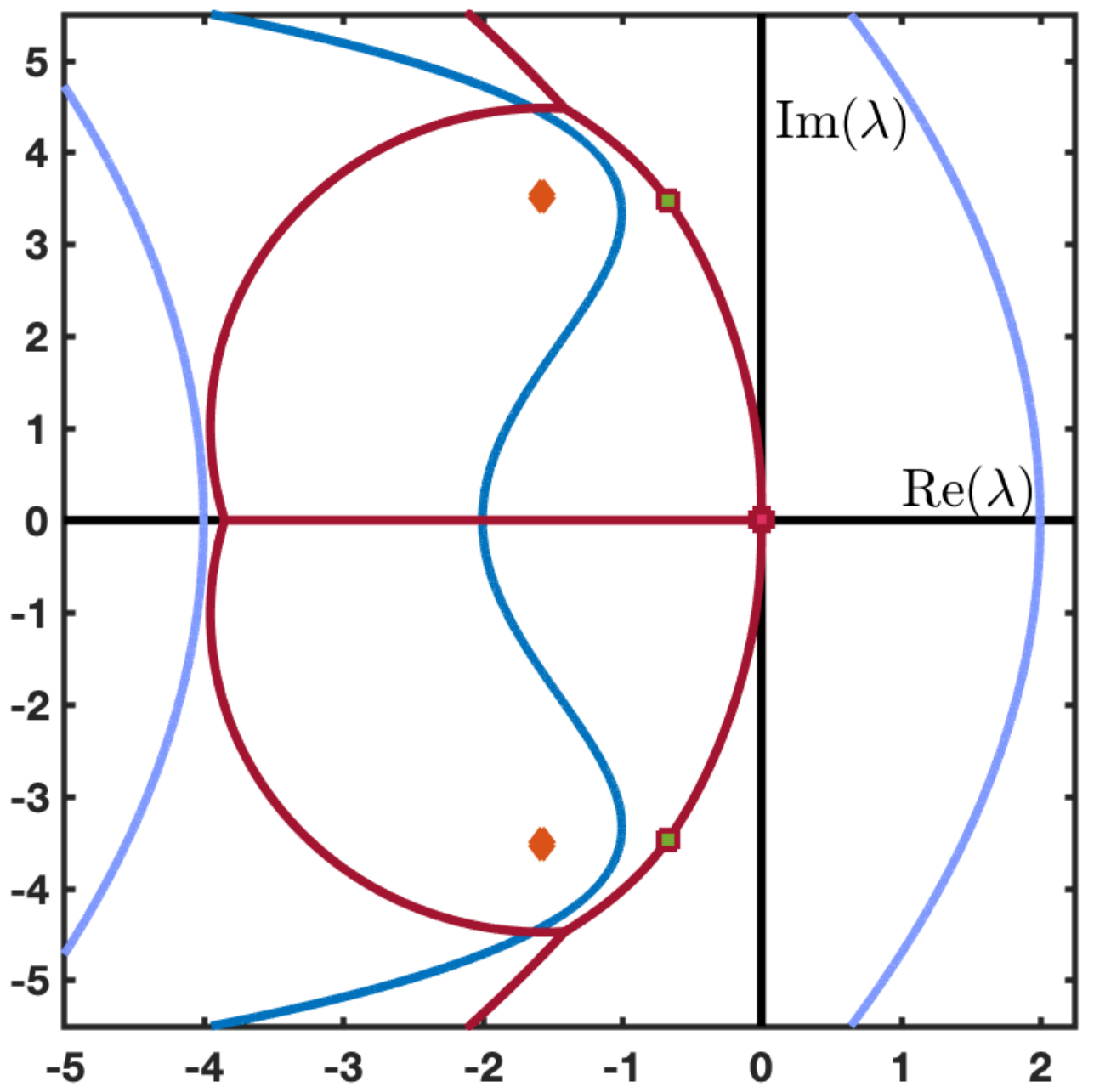}}\qquad
\subfigure{\includegraphics[width=0.243\textwidth]{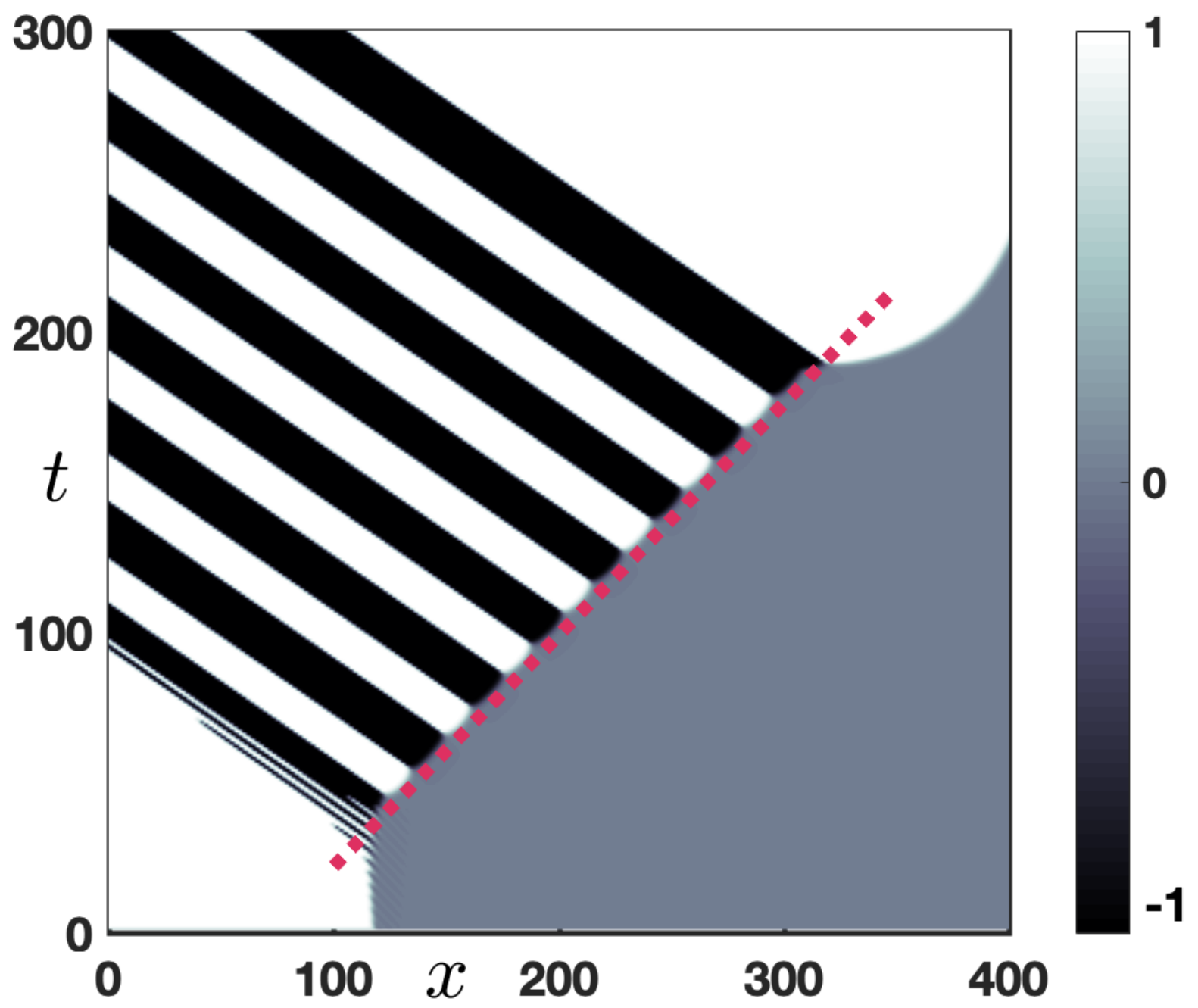}}
\caption{Spectra and simulations in $\mathcal{R}_\mathrm{abs}$: Essential spectra of $\calL^\pm_u$ (light blue), $\calL^\pm_v$ (dark blue), and absolute spectra of $\Sabs(\calL^+)$ (dark red) computed via continuation, computed with $s=s_*$ (left) and $s=s_\mathrm{abs}$ (middle; see \eqref{e:sabs}), for parameters  $(d,\alpha,\mu)=(1,1,-1)$ (top row), $(d,\alpha,\mu)=(1,1,-1/2)$ (middle row), and  $(d,\alpha,\mu)=(1/2,2,-1)$ (bottom row). The right column shows space-time plots of the $u$-component with best fits for the accelerated measured speed $s_\mathrm{m}=2.27,\, 2.35,\, 3.3$ (top to bottom), which compare well with the theoretically computed $s_\mathrm{abs}=2.2762,\,2.3547,\,3.3382.$}
\label{fig:NumericalSimulations}
\end{figure}

We carried out numerical simulations of \eqref{eq:main} in a frame moving with speed $s=s_*=\sqrt{\alpha d}$ based on a semi-implicit finite difference scheme, with localized initial conditions in both $u$- and $v$-components, varying parameters $d$, $\alpha$, and $\mu$. Figure \ref{fig:notremnant} shows spectra and a space-time plot of the $u$-component for a parameter choice in $\mathcal{R}_\mathrm{rem}$, where the front, stationary in the comoving frame, is stable as predicted. 
The subtle stabilization mechanism is visible in a log plot of the solution, shown in Figure \ref{fig:Logplotsolutions}. Clearly the $v$-component induces a weak, oscillatory  decay in the $u$-component, which travels at a larger speed to the right and decreases in overall amplitude. A $u$-front with such \emph{monotone} weaker decay would travel faster in the KPP equation. 
For parameter values in $\mathcal{R}_\mathrm{abs}$, our main result predicts stability, but simulations demonstrate instability after a long stable initial transient. Figure \ref{fig:NumericalSimulations} shows simulations in this parameter regime, demonstrating the acceleration and comparing a fit of the observed accelerated speed to the absolute spreading speed from \eqref{e:sabs}.

\begin{figure}[h!]
\centering
\subfigure[$(d,\alpha,\mu)=(1,1,-9)$.]{\includegraphics[width=0.243\textwidth]{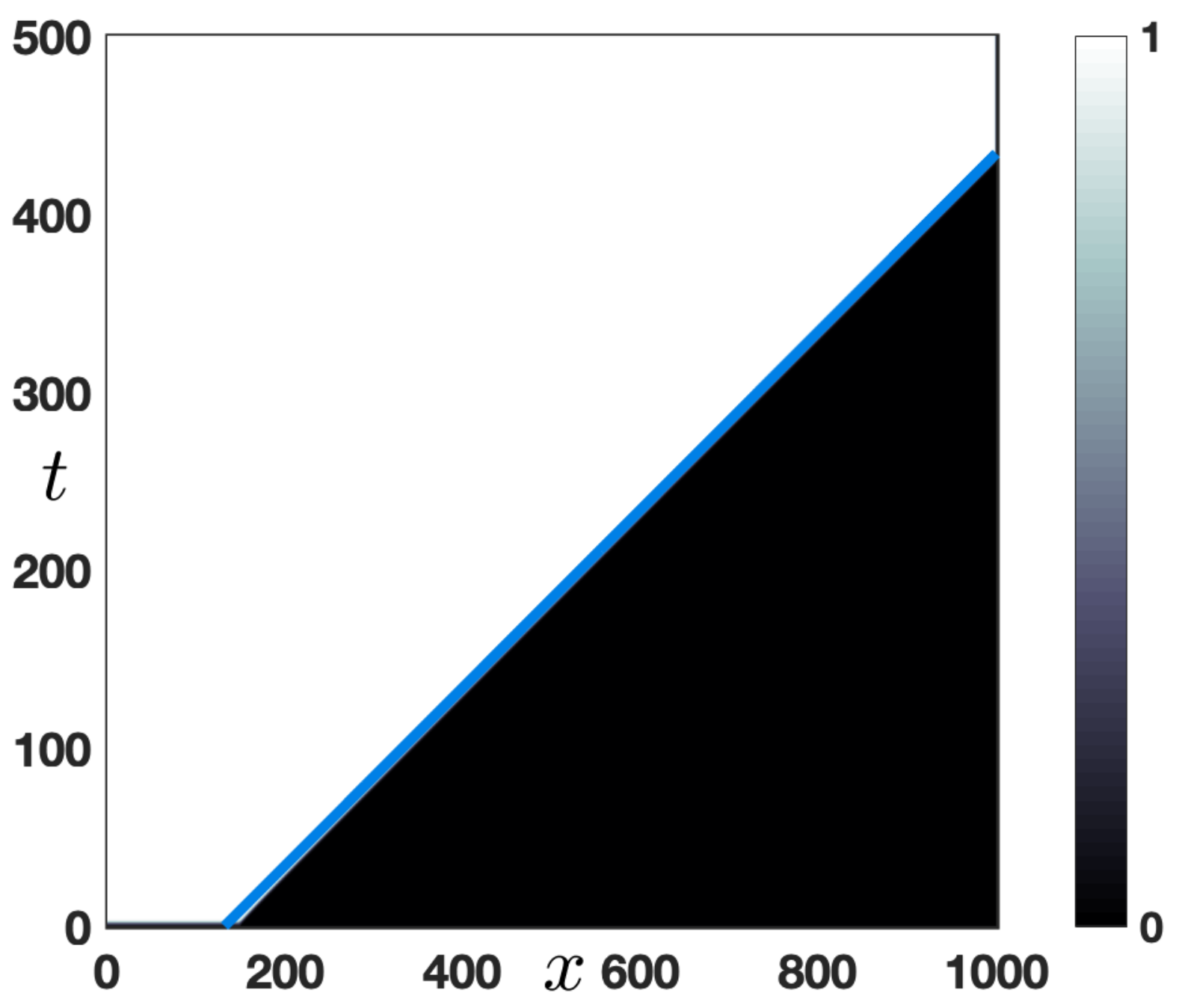}}\qquad
\subfigure[$(d,\alpha,\mu)=(1,1,-1/2)$.]
{\includegraphics[width=0.243\textwidth]{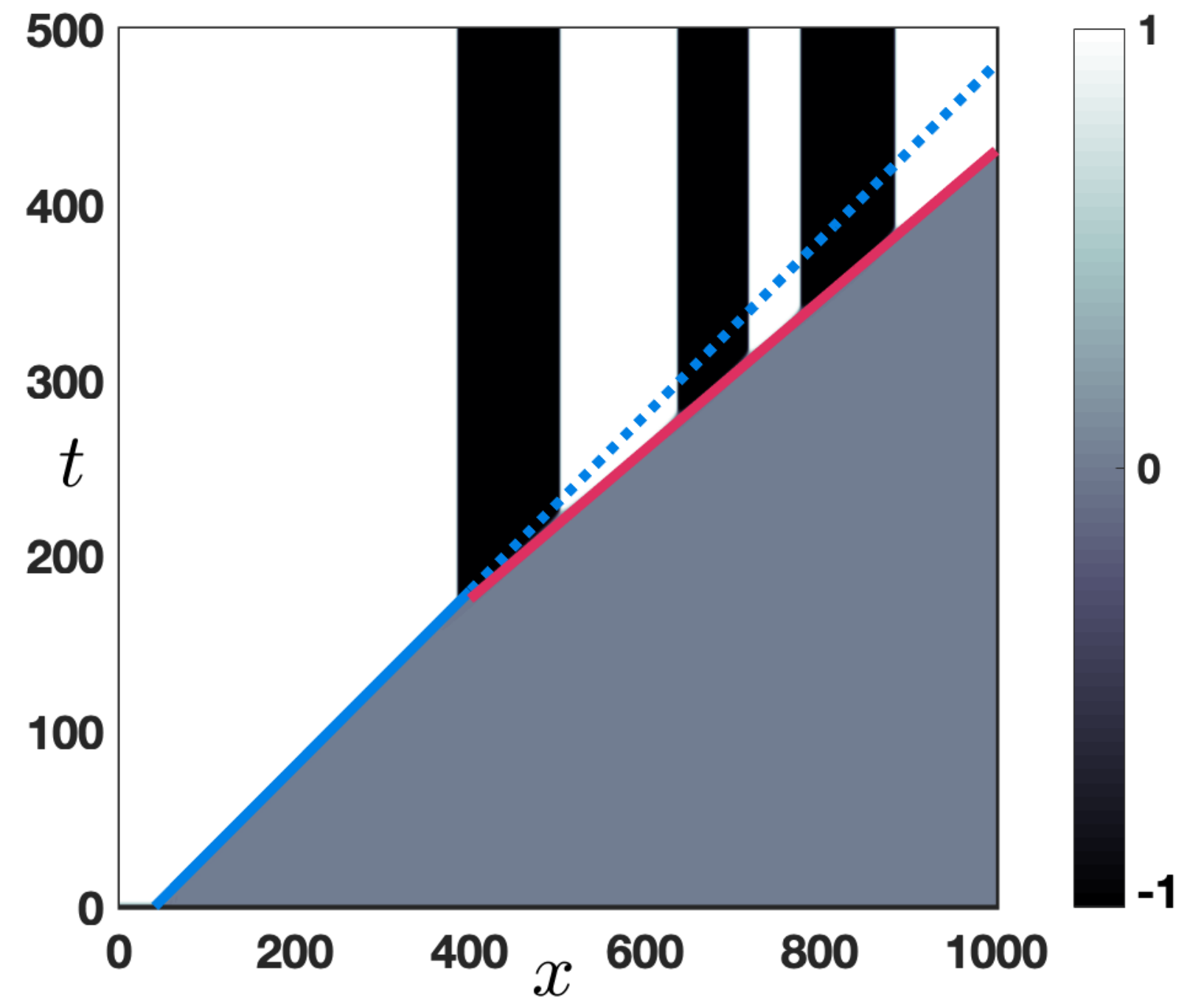}}\qquad
\subfigure[$(d,\alpha,\mu)=(1/2,2,-1)$.]{\includegraphics[width=0.24\textwidth]{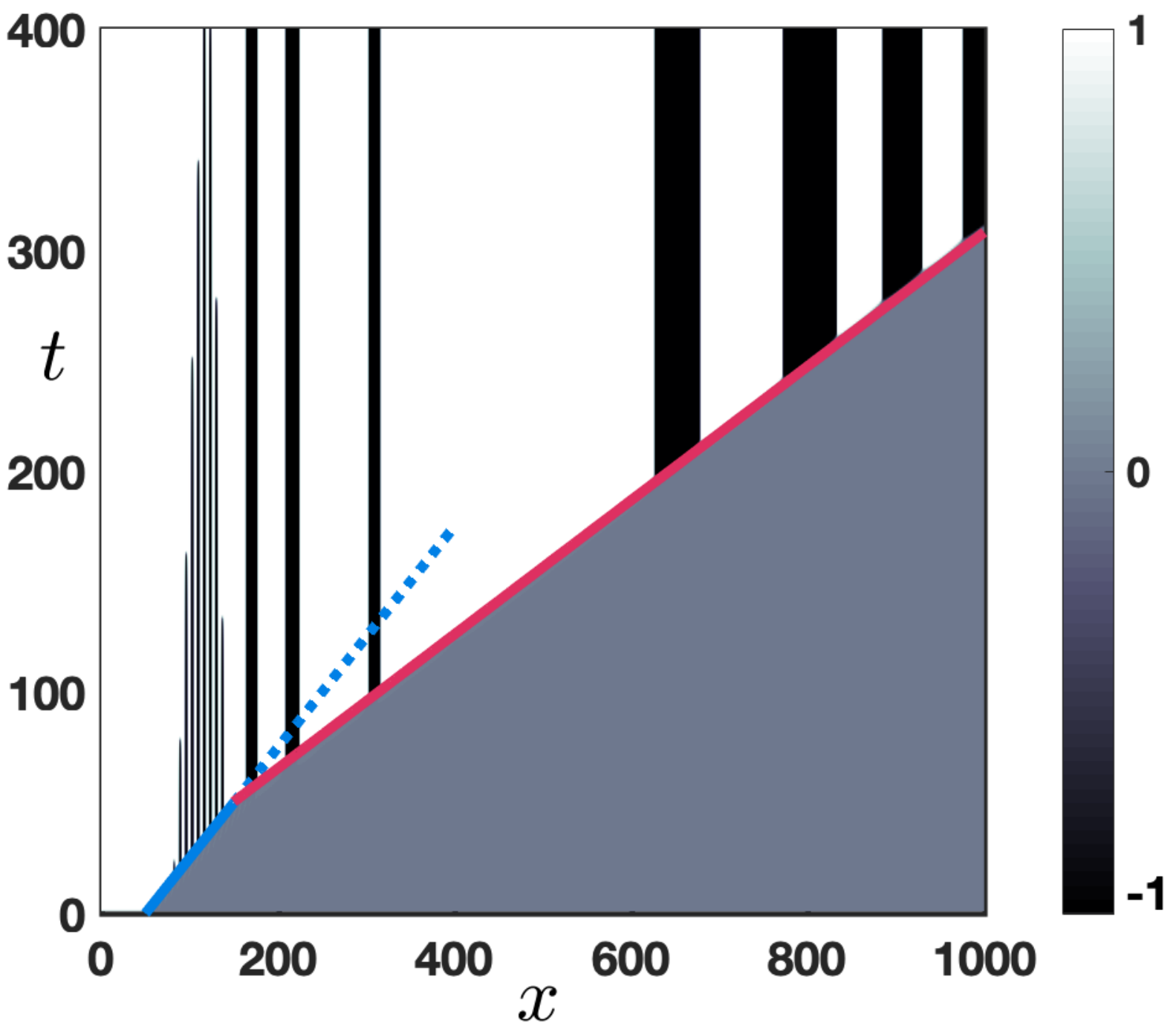}}
\caption{Space-time plots of the $u$-component in a steady frame $s=0$; parameters as in Figure~\ref{fig:notremnant}  and Figure~\ref{fig:NumericalSimulations}, middle and bottom row, that is, showing $\mathcal{R}_\mathrm{rem}$, $\mathcal{R}_\mathrm{abs}$, and  $\mathcal{R}_\mathrm{pw}$ (left to right). Line fits are initial speed $s_*$ (blue) and, after an initial transient, $s_\mathrm{abs}$ (magenta, center and right).}
\label{fig:NumSTPs0}
\end{figure}

In the steady laboratory frame, $s=0$, we observe the transient stability followed by acceleration; see Figure~\ref{fig:NumSTPs0}. 
 
In the remainder we offer an explanation of the eventual faster invasion speed $s_\mathrm{abs}$.  Our explanation is inspired by recent studies of resonant (or anomalous) invasion speeds  \cite{faye17,holzer14,holzer16,holzerscheel14}. We first study a linear, scalar equation toy model with (complex) exponential forcing and then derive  predictions for the original system (\ref{eq:main}).  
 
\paragraph{A motivating example.}
We consider 
\bqq u_t=du_{xx}+su_x+\alpha u +\beta (1+\sigma(x)) \me^{\nu x}, \label{eq:motivating} \eqq
for $\nu\in\mathbb{C}$ with $\mathrm{Re}(\nu)<0$ and where for simplicity,  we assume that $\sigma(x)$ is $2\pi$ periodic with Fourier Series expansion 
\[ \sigma(x)=\sum_{\ell\in\mathbb{Z}} c_\ell \me^{\mbi \ell x}. \]
The solution to this equation following a Laplace transform takes the form,
\[ u_\lambda(x)= c(\lambda)\me^{\nu_0^-(\lambda)x} +\frac{\beta}{D_u^0(\lambda,\nu)}\me^{\nu x}+\sum_{\ell\in\mathbb{Z}} \frac{\beta c_\ell }{D_u^0(\lambda,\nu+\mbi \ell)} \me^{(\nu+\mbi \ell)x}, \quad x>0. \]
Resonance poles arise at values of $\lambda$ for which $\nu_0^+(\lambda)=\nu+\mbi \ell$, for any $\ell\in\mathbb{Z}$. In fact, poles where $\nu_0^-(\lambda)=\nu+\mbi \ell$,  are not {\em relevant} in the sense of  \cite{holzer14,holzerscheel14}.  The most unstable modes are those with pure exponential decay, so the value of $l$ for which the imaginary part of $\nu+\mbi \ell$ is smallest in magnitude has the singularity with the largest real part and hence the dominant temporal growth rate.  Taking the inverse Laplace transform and using the theory of residues, the inhomogeneous terms contribute as
\[ \sum_{\ell\in\mathbb{Z}} \beta c_\ell \mathrm{Res}\left(\frac{1}{D_u^0(\lambda,\nu+\mbi \ell)},\lambda_\ell\right) \me^{\lambda_\ell t}\me^{(\nu+\mbi \ell)x}, \qquad \lambda_\ell=d(\nu+\mbi \ell)^2+s(\nu+\mbi \ell)+\alpha.  \]
Each term in the sum spreads at its envelope speed
\[ s_\textnormal{env}(\ell)= -\frac{\mathrm{Re}(\lambda_\ell)}{\mathrm{Re}(\nu)}, \]
where we recall that $\mathrm{Re}(\nu)<0$.  We expect the fastest mode to be observed which in this case is the one with largest temporal growth rate $\lambda_\ell$.  It is important to note that this prediction is independent of the magnitude of $c_\ell$, provided that it is non-zero.  To observe this, write
\bqq \me^{\ln(|c_\ell|)+\mathrm{Re}(\nu)x+\mathrm{Re}(\lambda_\ell) t }= \me^{\mathrm{Re}(\nu) \left(x+\frac{\ln(|c_\ell|)}{\mathrm{Re}(\nu)} +\frac{\mathrm{Re}(\lambda_\ell)}{\mathrm{Re}(\nu)} t\right) } = \me^{\mathrm{Re}(\nu) \left(x +\frac{\mathrm{Re}(\lambda_l)}{\mathrm{Re}(\nu)}\left( t+\frac{\ln(|c_\ell|)}{\mathrm{Re}(\lambda_\ell)} \right) \right) }, \label{eq:envelopecalc} \eqq
from which we predict a $\mathcal{O}\left(-\frac{\ln(|c_\ell|)}{\mathrm{Re}(\lambda_\ell)} \right)$ transient before the mode appears.  

\paragraph{The absolute spreading speed.}
Now consider the case where it is no longer a fixed mode $\me^{\nu x}$ that forms the inhomogeneous forcing term in (\ref{eq:motivating}) but rather a secondary partial differential equation as is the case in \eqref{eq:main}.  Then, after applying the Laplace transform, the forcing term is comprised of complex exponentials of the form $ \me^{\nu_v(\lambda)x}$ where, crucially, the mode $\nu_v$ now depends on $\lambda$.  Repeating the calcuations above, we now find that singularities arise whenever,
\[ D^0_\lambda(\lambda,\nu_v(\lambda)+\mbi \ell)=0. \]
We thus require values of $\lambda$ for which $\nu_0^+(\lambda)=\nu_v(\lambda)+\mbi \ell$. In other words, we would seek values of $\lambda$ for which two modes (one from the homogeneous $u$ equation and one from the coupled $v$ equation) have the same real part but whose imaginary parts may differ.  From the discussion in Section~\ref{sec:absoluteformain}, we note that this is equivalent to asking for $\lambda\in\Sigma_\textnormal{abs}(\cL^+)$.   Based upon this discussion, we would expect that any invasion front with unstable absolute spectrum  would be unstable when inhomogeneous coupling is introduced.  A prediction for the speed of the invasion front would be the speed at which there are marginally stable singularities, but no unstable singularities.  In other words, we expect the invasion speed to be the absolute spreading speed, defined in \eqref{e:sabs}; see also \cite{goh11}. The analysis presented thus far presents the basis of our Conjecture \ref{prop:sabs} stated in the introduction.

Of course, this faster invasion speed is predicated on the existence of the inhomgeneity $\sigma(x)$ which couples modes with the same exponential decay rate; a coupling which is notably absent in (\ref{eq:main}).  Yet, we did observe this faster speed in numerical simulations!  We attribute this to numerical effects, in particular round-off errors,  that effectively couple these modes in numerical simulations.

\begin{figure}[ht!]
\centering
\subfigure{\includegraphics[width=0.3\textwidth]{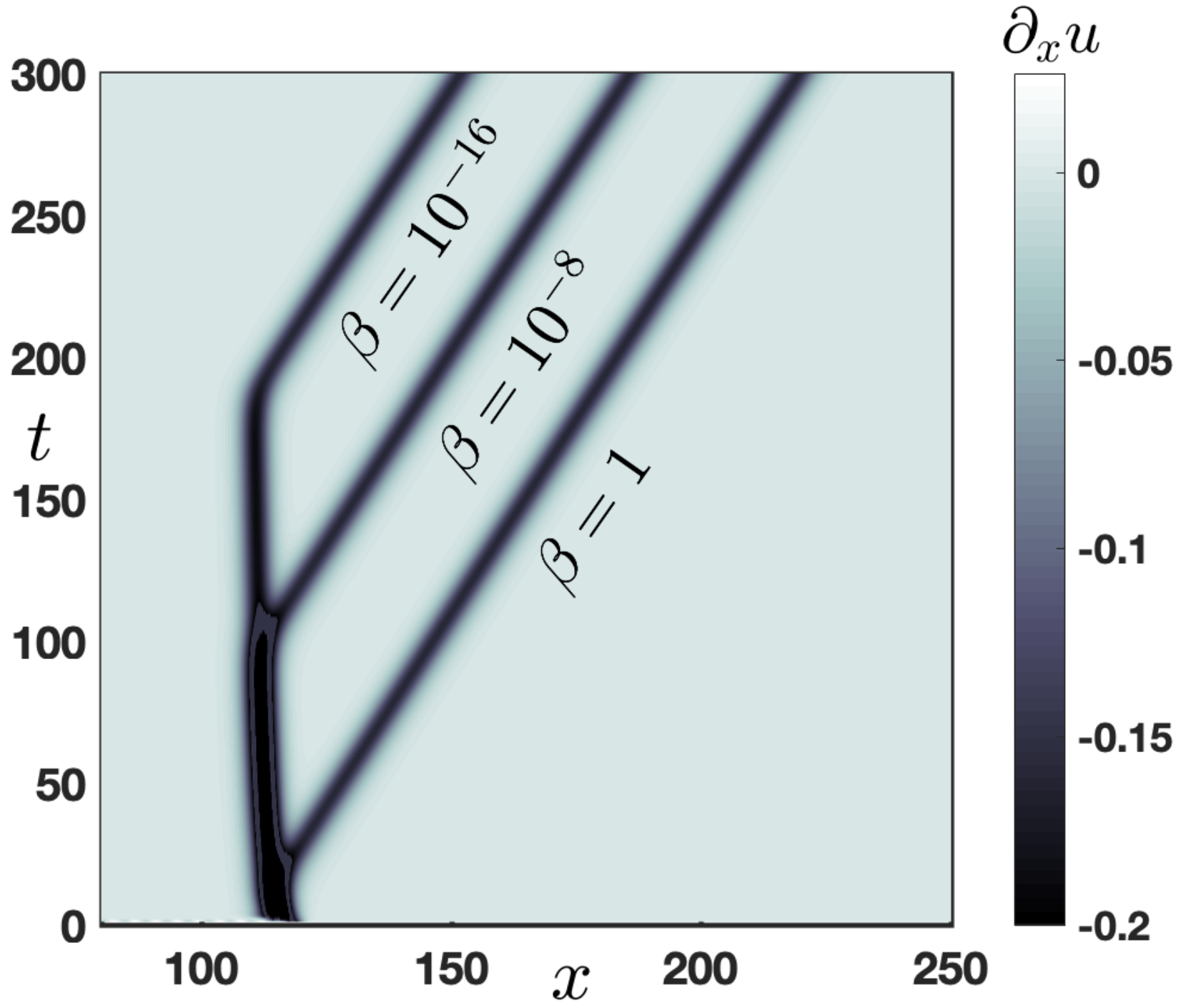}}
\hspace{0.8cm}
\subfigure{\includegraphics[width=0.3\textwidth]{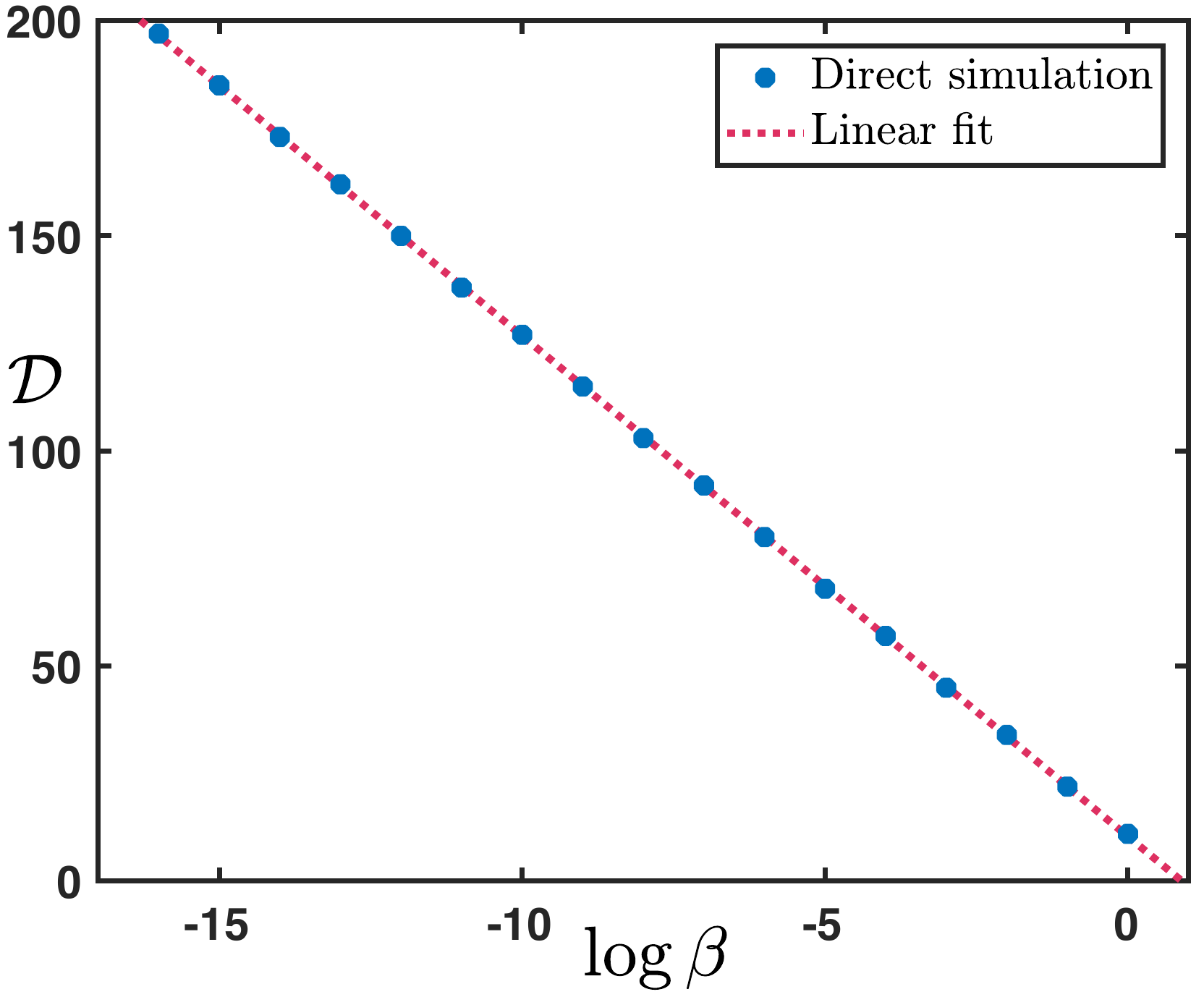}}
\caption{Left: Superimposed grayscale space-time plots of $\partial_x u$ to indicate the front location for  \eqref{eq:mainwithcos} with $s=s_*$ for values $\beta\in\left\{1,10^{-8},10^{-16}\right\}$ while other parameters are fixed to $(d,\alpha,\mu)=(1,1,-1/2)$. In that case $\ell_*=1.4872$ and $\nu=-0.5456$. The right panel shows the delay $\mathcal{D}$ of appearance of the faster invasion mode as a function of $\log_{10}\beta$. The measured slope is found to be $-11.6348$ which compares well with the prediction  $-\frac{\ln(10)}{\mathrm{Re}(\lambda_0)}=-11.1516$ from \eqref{eq:envelopecalc}.}
\label{fig:NumSimbeta}
\end{figure}

We will not embark on a full study of the numerical aspects of this phenomenon, here, but point to some evidence corroborating the relevance of the above calculations, beyond the fact that numerically observed accelerated speeds and theoretical predictions of $s_\mathrm{abs}$ agree well; see  Figure~\ref{fig:NumericalSimulations}.  We simulated
\bqq
\begin{cases}
\partial_t u&= d \partial_{xx} u+s\partial_x u+f(u)+\beta \cos(\ell_* x) v, \\
\partial_t v&= -(\partial_{xx} +1)^2v+s\partial_x v+\mu v,
\end{cases} 
\label{eq:mainwithcos}
\eqq
where the key difference to (\ref{eq:main}) is the $\cos(\ell_* x)$-term in the coupling.  We choose $\ell_*$ as follows.  Let $\lambda_{max}$ be the $\lambda$ value corresponding to the most unstable part of the absolute spectrum.  Then there exist a mode of the $v$ component (assume it is $\nu_2(\lambda)$) such that $\nu_0^+(\lambda_{max})=\nu_2(\lambda_{max})+\mbi \ell_*$.  Thus, the coupling term in (\ref{eq:mainwithcos}) should produce an unstable singularity which is (marginally) stabilized in a frame of reference moving with speed $s_\textnormal{abs}$. Results of simulations are presented in Figure~\ref{fig:NumSimbeta} for a variety of $\beta$ values.  For $\beta=1$, the faster invasion mode appears immediately, but for smaller values, the appearance of the faster invasion mode is delayed in a fashion consistent with formal calculation in (\ref{eq:envelopecalc}).

\begin{appendix}

\section{Proof of Lemma~\ref{lem:G12infinity}}\label{sec:G12infinityproof}
We proof Lemma~\ref{lem:G12infinity}.  We first consider the case $y>0$.  Continuity of $\bG_\lambda^{12,\infty}(x,y)$ at both $x=0$ and $x=y$ implies
\begin{align*}
b_2+b_3+b_4 &= \beta \sum_{j=3}^4 \left(\frac{1}{D_u^0(\lambda,\nu_j(\lambda))}-\frac{1}{D_u^1(\lambda,\nu_j(\lambda))} \right) c_j(\lambda) \me^{-\nu_j (\lambda)y},\\
b_1\me^{\nu_-^0(\lambda)y}+b_2\me^{\nu_-^0(\lambda)y}+b_3\me^{\nu_+^0(\lambda)y}&=\beta \sum_{j=1}^4 \frac{c_j(\lambda)}{D_u^0(\lambda,\nu_j(\lambda))}.
\end{align*}
Continuity of $\partial_x \bG_\lambda^{12,\infty}(x,y)$ at both $x=0$ and $x=y$ then implies 
\begin{align*}
b_2\nu_-^0(\lambda)+b_3\nu_+^0(\lambda)+b_4\nu_+^1(\lambda) &= \beta \sum_{j=3}^4 \left(\frac{1}{D_u^0(\lambda,\nu_j(\lambda))}-\frac{1}{D_u^1(\lambda,\nu_j(\lambda))} \right) \\
& \qquad\qquad \cdot c_j(\lambda)\nu_j(\lambda) \me^{-\nu_j (\lambda)y},\\
b_1\nu_-^0(\lambda)\me^{\nu_-^0(\lambda)y}+b_2\nu_-^0(\lambda)\me^{\nu_-^0(\lambda)y}+b_3\nu_+^0(\lambda)\me^{\nu_+^0(\lambda)y}&=\beta \sum_{j=1}^4 \frac{c_j(\lambda)\nu_j(\lambda)}{D_u^0(\lambda,\nu_j(\lambda))}.
\end{align*}
We thus arrive at 
\bqs
\underbrace{\left( \begin{matrix}\me^{\nu_-^0(\lambda)y} & \me^{\nu_-^0(\lambda)y} & \me^{\nu_+^0(\lambda)y} & 0 \\
\nu_-^0(\lambda)\me^{\nu_-^0(\lambda)y} & \nu_-^0(\lambda)\me^{\nu_-^0(\lambda)y} & \nu_+^0(\lambda)\me^{\nu_+^0(\lambda)y} & 0 \\
0 & 1 & 1 & 1 \\
0 & \nu_-^0(\lambda)&\nu_+^0(\lambda)& \nu_+^1(\lambda) \end{matrix}\right)}_{\coloneqq\mathcal{A}_1(\lambda)}\left(\begin{matrix}b_1 \\ b_2 \\ b_3 \\ b_4 \end{matrix} \right)=\beta \mathcal{V},
\eqs
where the right-hand side is given by
\bqs
\mathcal{V}=\sum_{j=1}^4 \frac{c_j(\lambda)}{D_u^0(\lambda,\nu_j(\lambda))} \left(\begin{matrix} 1 \\ \nu_j(\lambda) \\ 0 \\ 0\end{matrix}\right) + \sum_{j=3}^4 \left(\frac{1}{D_u^0(\lambda,\nu_j(\lambda))}-\frac{1}{D_u^1(\lambda,\nu_j(\lambda))} \right) c_j(\lambda) \me^{-\nu_j (\lambda)y}\left(\begin{matrix} 0 \\ 0 \\ 1 \\ \nu_j(\lambda) \end{matrix}\right).
\eqs
The inverse of the matrix $\mathcal{A}_1(\lambda)$ reads
\bqs
\mathcal{A}^{-1}_1(\lambda)=\left( \begin{matrix}
\frac{\nu_-^0(\lambda) p_1(\lambda) \me^{-\nu_+^0(\lambda)y} -\nu_+^0(\lambda) \me^{-\nu_-^0(\lambda)y}}{\nu_-^0(\lambda)-\nu_+^0(\lambda)} & \frac{- p_1(\lambda) \me^{-\nu_+^0(\lambda)y} +\me^{-\nu_-^0(\lambda)y}}{\nu_-^0(\lambda)-\nu_+^0(\lambda)} & -\frac{\nu_+^1(\lambda)}{\nu_+^1(\lambda)-\nu_-^0(\lambda)} & \frac{1}{\nu_+^1(\lambda)-\nu_-^0(\lambda)} \\
-\frac{\nu_-^0(\lambda) p_1(\lambda) \me^{-\nu_+^0(\lambda)y}}{\nu_-^0(\lambda)-\nu_+^0(\lambda)}  &  \frac{p_1(\lambda) \me^{-\nu_+^0(\lambda)y}}{\nu_-^0(\lambda)-\nu_+^0(\lambda)}   & \frac{\nu_+^1(\lambda)}{\nu_+^1(\lambda)-\nu_-^0(\lambda)} & -\frac{1}{\nu_+^1(\lambda)-\nu_-^0(\lambda)} \\
 \frac{\nu_-^0(\lambda) \me^{-\nu_+^0(\lambda)y}}{\nu_-^0(\lambda)-\nu_+^0(\lambda)}  &  -\frac{\me^{-\nu_+^0(\lambda)y}}{\nu_-^0(\lambda)-\nu_+^0(\lambda)} & 0 & 0 \\
 - \frac{\nu_-^0(\lambda) \me^{-\nu_+^0(\lambda)y}}{\nu_-^0(\lambda)-\nu_+^1(\lambda)}   & \frac{\me^{-\nu_+^0(\lambda)y}}{\nu_-^0(\lambda)-\nu_+^1(\lambda)} & -\frac{\nu_-^0(\lambda)}{\nu_+^1(\lambda)-\nu_-^0(\lambda)} & \frac{1}{\nu_+^1(\lambda)-\nu_-^0(\lambda)}
\end{matrix}\right)
\eqs
where we denoted
\bqs
p_1(\lambda)=\frac{\nu_+^0(\lambda)-\nu_+^1(\lambda)}{\nu_-^0(\lambda)-\nu_+^1(\lambda)}.
\eqs
As a consequence, we get the following expressions for the $b_j(\lambda,y)$:
\begin{align*}
\frac{b_1(\lambda,y)}{\beta}&=\frac{p_1(\lambda)}{\nu_-^0(\lambda)-\nu_+^0(\lambda)}\sum_{j=1}^4 \frac{c_j(\lambda)(\nu_-^0(\lambda)-\nu_j(\lambda))}{D_u^0(\lambda,\nu_j(\lambda))} \me^{-\nu_+^0(\lambda)y} \nonumber  \\
&~~~-\frac{1}{\nu_-^0(\lambda)-\nu_+^0(\lambda)}\sum_{j=1}^4 \frac{c_j(\lambda)(\nu_+^0(\lambda)-\nu_j(\lambda))}{D_u^0(\lambda,\nu_j(\lambda))} \me^{-\nu_-^0(\lambda)y} \nonumber\\
&~~~-\frac{1}{\nu_+^1(\lambda)-\nu_-^0(\lambda)}\sum_{j=3}^4 \left(\frac{1}{D_u^0(\lambda,\nu_j(\lambda))}-\frac{1}{D_u^1(\lambda,\nu_j(\lambda))} \right) c_j(\lambda)(\nu_+^1(\lambda)-\nu_j(\lambda)) \me^{-\nu_j (\lambda)y} \nonumber 
\end{align*}
\begin{align*}
\frac{b_2(\lambda,y)}{\beta}&=-\frac{p_1(\lambda)}{\nu_-^0(\lambda)-\nu_+^0(\lambda)}\sum_{j=1}^4 \frac{c_j(\lambda)(\nu_-^0(\lambda)-\nu_j(\lambda))}{D_u^0(\lambda,\nu_j(\lambda))} \me^{-\nu_+^0(\lambda)y} \nonumber\\
&~~~+\frac{1}{\nu_+^1(\lambda)-\nu_-^0(\lambda)}\sum_{j=3}^4 \left(\frac{1}{D_u^0(\lambda,\nu_j(\lambda))}-\frac{1}{D_u^1(\lambda,\nu_j(\lambda))} \right) c_j(\lambda)(\nu_+^1(\lambda)-\nu_j(\lambda)) \me^{-\nu_j (\lambda)y}  \nonumber 
\end{align*}
\begin{align*}
\frac{b_3(\lambda,y)}{\beta}&=\frac{1}{\nu_-^0(\lambda)-\nu_+^0(\lambda)}\sum_{j=1}^4 \frac{c_j(\lambda)(\nu_-^0(\lambda)-\nu_j(\lambda))}{D_u^0(\lambda,\nu_j(\lambda))} \me^{-\nu_+^0(\lambda)y}  \nonumber \\
\frac{b_4(\lambda,y)}{\beta}&=-\frac{1}{\nu_-^0(\lambda)-\nu_+^1(\lambda)}\sum_{j=1}^4 \frac{c_j(\lambda)(\nu_-^0(\lambda)-\nu_j(\lambda))}{D_u^0(\lambda,\nu_j(\lambda))} \me^{-\nu_+^0(\lambda)y} \nonumber \\
&~~~-\frac{1}{\nu_-^0(\lambda)-\nu_+^1(\lambda)}\sum_{j=3}^4 \left(\frac{1}{D_u^0(\lambda,\nu_j(\lambda))}-\frac{1}{D_u^1(\lambda,\nu_j(\lambda))} \right) c_j(\lambda)(\nu_-^0(\lambda)-\nu_j(\lambda)) \me^{-\nu_j (\lambda)y}. \label{eq:bs}
\end{align*}
For $y<0$, we arrive at a system
\bqs
\underbrace{\left( \begin{matrix}
1 & 1 & 1 & 0 \\
\nu_-^0(\lambda)&\nu_-^1(\lambda)& \nu_+^1(\lambda) & 0 \\
0 & \me^{\nu_-^1(\lambda)y} & \me^{\nu_+^1(\lambda)y} & \me^{\nu_+^1(\lambda)y}  \\
0 & \nu_-^1(\lambda)\me^{\nu_-^1(\lambda)y} & \nu_+^1(\lambda)\me^{\nu_+^1(\lambda)y} & \nu_+^1(\lambda)\me^{\nu_+^1(\lambda)y} 
\end{matrix}\right)}_{\coloneqq\mathcal{A}_2(\lambda)}\left(\begin{matrix}h_1 \\ h_2 \\ h_3 \\ h_4 \end{matrix} \right)=\beta  \mathcal{W},
\eqs
where
\[
\mathcal{W}\coloneqq\sum_{j=1}^4 \frac{c_j(\lambda)}{D_u^1(\lambda,\nu_j(\lambda))} \left(\begin{matrix}  0 \\ 0 \\ 1 \\ \nu_j(\lambda) \end{matrix}\right)+ \sum_{j=1}^2 \left(\frac{1}{D_u^0(\lambda,\nu_j(\lambda))}-\frac{1}{D_u^1(\lambda,\nu_j(\lambda))} \right) c_j(\lambda) \me^{-\nu_j (\lambda)y}\left(\begin{matrix} 1 \\ \nu_j(\lambda) \\ 0 \\ 0\end{matrix}\right).
\]
The inverse of the matrix $\mathcal{A}_2(\lambda)$ reads
\bqs
\mathcal{A}_2^{-1}(\lambda)=\left(
\begin{matrix}
-\frac{\nu_+^1(\lambda)}{\nu_-^0(\lambda)-\nu_+^1(\lambda)} & \frac{1}{\nu_-^0(\lambda)-\nu_+^1(\lambda)} & \frac{\nu_+^1(\lambda)\me^{-\nu_-^1(\lambda)y}}{\nu_-^0(\lambda)-\nu_+^1(\lambda)} & -\frac{\me^{-\nu_-^1(\lambda)y}}{\nu_-^0(\lambda)-\nu_+^1(\lambda)} \\
0  & 0 & -\frac{\nu_+^1(\lambda)\me^{-\nu_-^1(\lambda)y}}{\nu_-^1(\lambda)-\nu_+^1(\lambda)} & \frac{\me^{-\nu_-^1(\lambda)y}}{\nu_-^1(\lambda)-\nu_+^1(\lambda)} \\
\frac{\nu_-^0(\lambda)}{\nu_-^0(\lambda)-\nu_+^1(\lambda)}   & -\frac{1}{\nu_-^0(\lambda)-\nu_+^1(\lambda)} & \frac{p_2(\lambda)\nu_+^1(\lambda)\me^{-\nu_-^1(\lambda)y}}{\nu_-^1(\lambda)-\nu_+^1(\lambda)} & -\frac{p_2(\lambda)\me^{-\nu_-^1(\lambda)y}}{\nu_-^1(\lambda)-\nu_+^1(\lambda)} \\
-\frac{\nu_-^0(\lambda)}{\nu_-^0(\lambda)-\nu_+^1(\lambda)}    & \frac{1}{\nu_-^0(\lambda)-\nu_+^1(\lambda)} & \frac{-p_2(\lambda)\nu_+^1(\lambda)\me^{-\nu_-^1(\lambda)y}+\nu_-^1(\lambda)\me^{-\nu_+^1(\lambda)y}}{\nu_-^1(\lambda)-\nu_+^1(\lambda)} & \frac{p_2(\lambda)\me^{-\nu_-^1(\lambda)y}-\me^{-\nu_+^1(\lambda)y}}{\nu_-^1(\lambda)-\nu_+^1(\lambda)}
\end{matrix}
\right)
\eqs
where we denoted
\bqs
p_2(\lambda)\coloneqq\frac{\nu_-^0(\lambda)-\nu_-^1(\lambda)}{\nu_-^0(\lambda)-\nu_+^1(\lambda)} .
\eqs
As a consequence, we get the following expressions for the $h_j(\lambda,y)$: 
\begin{align*}
\frac{h_1(\lambda,y)}{\beta}&=\frac{1}{\nu_-^0(\lambda)-\nu_+^1(\lambda)}\sum_{j=1}^4 \frac{c_j(\lambda)(\nu_+^1(\lambda)-\nu_j(\lambda))}{D_u^1(\lambda,\nu_j(\lambda))}\me^{-\nu_-^1(\lambda)y} \nonumber \\
&~~~-\frac{1}{\nu_-^0(\lambda)-\nu_+^1(\lambda)}\sum_{j=1}^2 \left(\frac{1}{D_u^0(\lambda,\nu_j(\lambda))}-\frac{1}{D_u^1(\lambda,\nu_j(\lambda))} \right) c_j(\lambda)(\nu_+^1(\lambda)-\nu_j(\lambda)) \me^{-\nu_j (\lambda)y},\nonumber \\
\frac{h_2(\lambda,y)}{\beta}&=-\frac{1}{\nu_-^1(\lambda)-\nu_+^1(\lambda)}\sum_{j=1}^4 \frac{c_j(\lambda)(\nu_+^1(\lambda)-\nu_j(\lambda))}{D_u^1(\lambda,\nu_j(\lambda))}\me^{-\nu_-^1(\lambda)y},\nonumber 
\end{align*}
\begin {align*}
\frac{h_3(\lambda,y)}{\beta}&=\frac{p_2(\lambda)}{\nu_-^1(\lambda)-\nu_+^1(\lambda)}\sum_{j=1}^4 \frac{c_j(\lambda)(\nu_+^1(\lambda)-\nu_j(\lambda))}{D_u^1(\lambda,\nu_j(\lambda))}\me^{-\nu_-^1(\lambda)y}\nonumber \\
&~~~+\frac{1}{\nu_-^0(\lambda)-\nu_+^1(\lambda)}\sum_{j=1}^2 \left(\frac{1}{D_u^0(\lambda,\nu_j(\lambda))}-\frac{1}{D_u^1(\lambda,\nu_j(\lambda))} \right) c_j(\lambda)(\nu_-^0(\lambda)-\nu_j(\lambda)) \me^{-\nu_j (\lambda)y},\nonumber\\
\frac{h_4(\lambda,y)}{\beta}&=-\frac{p_2(\lambda)}{\nu_-^1(\lambda)-\nu_+^1(\lambda)}\sum_{j=1}^4 \frac{c_j(\lambda)(\nu_+^1(\lambda)-\nu_j(\lambda))}{D_u^1(\lambda,\nu_j(\lambda))}\me^{-\nu_-^1(\lambda)y}\nonumber\\
&~~~+\frac{1}{\nu_-^1(\lambda)-\nu_+^1(\lambda)}\sum_{j=1}^4 \frac{c_j(\lambda)(\nu_-^1(\lambda)-\nu_j(\lambda))}{D_u^1(\lambda,\nu_j(\lambda))}\me^{-\nu_+^1(\lambda)y}\nonumber\\
&~~~-\frac{1}{\nu_-^0(\lambda)-\nu_+^1(\lambda)}\sum_{j=1}^2 \left(\frac{1}{D_u^0(\lambda,\nu_j(\lambda))}-\frac{1}{D_u^1(\lambda,\nu_j(\lambda))} \right) c_j(\lambda)(\nu_-^0(\lambda)-\nu_j(\lambda)) \me^{-\nu_j (\lambda)y}
\end{align*}

Let us define the following quantities:
\begin{align*}
\mathbb{D}^0_\pm(\lambda)&\coloneqq\sum_{j=1}^4\frac{c_j(\lambda)(\nu_\pm^0(\lambda)-\nu_j(\lambda))}{D_u^0(\lambda,\nu_j(\lambda))}, \quad \mathbb{D}^1_\pm(\lambda)\coloneqq\sum_{j=1}^4\frac{c_j(\lambda)(\nu_\pm^1(\lambda)-\nu_j(\lambda))}{D_u^1(\lambda,\nu_j(\lambda))},\\
\mathbb{B}_j^{1,+}(\lambda)&\coloneqq\left(\frac{1}{D_u^0(\lambda,\nu_j(\lambda))}-\frac{1}{D_u^1(\lambda,\nu_j(\lambda))} \right) c_j(\lambda)(\nu_+^1(\lambda)-\nu_j(\lambda)),\\
\mathbb{B}_j^{0,-}(\lambda)&\coloneqq\left(\frac{1}{D_u^0(\lambda,\nu_j(\lambda))}-\frac{1}{D_u^1(\lambda,\nu_j(\lambda))} \right) c_j(\lambda)(\nu_-^0(\lambda)-\nu_j(\lambda)).
\end{align*}
Then we have the condensed expressions for:
\begin{flalign}
\frac{h_1(\lambda,y)}{\beta}&=\frac{\mathbb{D}^1_+(\lambda)}{\nu_-^0(\lambda)-\nu_+^1(\lambda)}\me^{-\nu_-^1(\lambda)y}-\frac{1}{\nu_-^0(\lambda)-\nu_+^1(\lambda)}\sum_{j=1}^2 \mathbb{B}_j^{1,+}(\lambda) \me^{-\nu_j (\lambda)y}, &&\nonumber\\
\frac{h_2(\lambda,y)}{\beta}&=-\frac{\mathbb{D}^1_+(\lambda)}{\nu_-^1(\lambda)-\nu_+^1(\lambda)}\me^{-\nu_-^1(\lambda)y},&&\nonumber 
\end{flalign}
\begin{flalign}
\frac{h_3(\lambda,y)}{\beta}&=\frac{p_2(\lambda)}{\nu_-^1(\lambda)-\nu_+^1(\lambda)}\mathbb{D}^1_+(\lambda)\me^{-\nu_-^1(\lambda)y}+\frac{1}{\nu_-^0(\lambda)-\nu_+^1(\lambda)}\sum_{j=1}^2 \mathbb{B}_j^{0,-}(\lambda) \me^{-\nu_j (\lambda)y},\nonumber \\
\frac{h_4(\lambda,y)}{\beta}&=-\frac{p_2(\lambda)}{\nu_-^1(\lambda)-\nu_+^1(\lambda)}\mathbb{D}^1_+(\lambda)\me^{-\nu_-^1(\lambda)y}+\frac{\mathbb{D}^1_-(\lambda)}{\nu_-^1(\lambda)-\nu_+^1(\lambda)}\me^{-\nu_+^1(\lambda)y}&&\nonumber\\
&\quad -\frac{1}{\nu_-^0(\lambda)-\nu_+^1(\lambda)}\sum_{j=1}^2 \mathbb{B}_j^{0,-}(\lambda)   \me^{-\nu_j (\lambda)y}, &&\label{eq:hs}
\end{flalign}
together with
\begin{flalign}
\frac{b_1(\lambda,y)}{\beta}&=\frac{p_1(\lambda)}{\nu_-^0(\lambda)-\nu_+^0(\lambda)}\mathbb{D}_-^0(\lambda) \me^{-\nu_+^0(\lambda)y} \nonumber&&\\
&\quad -\frac{1}{\nu_-^0(\lambda)-\nu_+^0(\lambda)}\mathbb{D}_+^0(\lambda)\me^{-\nu_-^0(\lambda)y}-\frac{1}{\nu_+^1(\lambda)-\nu_-^0(\lambda)}\sum_{j=3}^4 \mathbb{B}_j^{1,+}(\lambda) \me^{-\nu_j (\lambda)y}, \nonumber 
\end{flalign}
\begin {flalign}
\frac{b_2(\lambda,y)}{\beta}&=-\frac{p_1(\lambda)}{\nu_-^0(\lambda)-\nu_+^0(\lambda)}\mathbb{D}_-^0(\lambda) \me^{-\nu_+^0(\lambda)y}+\frac{1}{\nu_+^1(\lambda)-\nu_-^0(\lambda)}\sum_{j=3}^4 \mathbb{B}_j^{1,+}(\lambda) \me^{-\nu_j (\lambda)y},\nonumber \\
\frac{b_3(\lambda,y)}{\beta}&=\frac{1}{\nu_-^0(\lambda)-\nu_+^0(\lambda)}\mathbb{D}_-^0(\lambda) \me^{-\nu_+^0(\lambda)y},\nonumber&& \\
\frac{b_4(\lambda,y)}{\beta}&=-\frac{1}{\nu_-^0(\lambda)-\nu_+^1(\lambda)}\mathbb{D}_-^0(\lambda) \me^{-\nu_+^0(\lambda)y}-\frac{1}{\nu_-^0(\lambda)-\nu_+^1(\lambda)}\sum_{j=3}^4  \mathbb{B}_j^{0,-}(\lambda)\me^{-\nu_j (\lambda)y}. &&\label{eq:bs}
\end{flalign}

\begin{rmk} \label{rmk:nosingularity}
The quantities $b_1(\lambda,y)$,  $b_2(\lambda,y)$, and  $b_3(\lambda,y)$ all contain $\sqrt{\lambda}$ singularities.  Note that this singularity in  $b_1(\lambda,y)$ is removable after consulting the formulas for $\mathbb{D}_\pm^0(\lambda)$.   Likewise, the terms  $b_2(\lambda,y)$, and  $b_3(\lambda,y)$ appear in tandem in the expression for $\bG_\lambda^{12}(x,y)$ where again the singularity is removable due to cancellation. 
\end{rmk}

\section{Proof of Lemma~\ref{lem:G12full}} \label{sec:G12full}
In this section, we present the proof of Lemma~\ref{lem:G12full}.  Rather than present all six cases, we show the details for two representative cases.  Recall the definition
\[ h(x)=-\left(f'(Q_*(x))-f_\infty(x)\right).\]
We will also use throughout these estimates the fact that 
\[  \W_\lambda(\tau)=\W_\lambda(0)\me^{-\frac{s}{d}\tau}, \quad \tau\in\R. \]

\paragraph{Case $y<0<x$.} For this arrangement the pointwise Green's function has the following expression
\begin{align*} \bG_\lambda^{12}(x,y) & =h_1(\lambda,y)\me^{\nu_-^0(\lambda)x}-\beta \sum_{j=1}^2\frac{c_j(\lambda)}{D_u^0(\lambda,\nu_j(\lambda))}\me^{\nu_j(\lambda)(x-y)}\\
&~~~+\varphi^+(x) \int_{-\infty}^x  
\frac{\varphi^-(\tau)}{\W_\lambda(\tau)}h(\tau) \bG_\lambda^{12,\infty}(\tau,y)\md \tau +\varphi^-(x) \int_x^\infty  
\frac{\varphi^+(\tau)}{\W_\lambda(\tau)}h(\tau) \bG_\lambda^{12,\infty}(\tau,y)\md \tau.
\end{align*}
For the leading order terms the expression for $h_1(\lambda,y)$ shows that the first term can be absorbed into $\boldH_\lambda(x,y)$ (after factoring out the singularity at $\lambda_v^{\textnormal{bp}}$) and $\bI_{\lambda,j}(x,y)$ (after factoring out an additional singularity at $\lambda^{\textnormal{rp}}$) while the inhomogeneous terms contribute to $\boldJ_{\lambda,j}(x,y)$.  We then begin our treatment of the integral terms with the second integral where we use the specific form of $\bG_\lambda^{12,\infty}(\tau,y)$ for $y<0<\tau$.
Using the decomposition of $\varphi^-(x)$ valid for $x>0$ and focusing on the terms proportional to $h_1(\lambda,y)$ the integral then assumes the form,
\begin{align*} \frac{C(\lambda)}{\W_\lambda(0)}h_1(\lambda,y)\me^{\nu_+^0(\lambda)x}(1+\kappa_+(x,\lambda))\int_x^\infty \me^{(\nu_-^0(\lambda)-\nu_+^0(\lambda))\tau}(1+\theta_+(\tau,\lambda)) h(\tau) \md\tau \\
 + \frac{D(\lambda)}{\W_\lambda(0)}h_1(\lambda,y)\me^{\nu_-^0(\lambda)x}(1+\theta_+(x,\lambda))\int_x^\infty(1+\theta_+(\tau,\lambda)) h(\tau) \md\tau. \end{align*}
Convergence of the integrals then follows and we see that the first integral is $\mathcal{O}(\me^{(\nu_-^0(\lambda)-\nu_+^0(\lambda)-\vartheta)x})$ and after again consulting the formula for $h_1(\lambda,y)$ we see that this term again contributes to $\boldH_\lambda(x,y)$.  The remaining terms for this integral are of the form
\begin{align*} -\beta \sum_{j=1}^2 \frac{C(\lambda)(1+\kappa_+(x,\lambda))}{\W_\lambda(0)D_u^0(\lambda,\nu_j(\lambda))}c_j(\lambda,y)\me^{\nu_+^0(\lambda)x}\me^{-\nu_j(\lambda)y}\int_x^\infty \me^{(-\nu_+^0(\lambda)+\nu_j(\lambda))\tau}(1+\theta_+(\tau,\lambda)) h(\tau) \md\tau \\
-\beta \sum_{j=1}^2 \frac{D(\lambda)(1+\theta_+(x,\lambda))}{\W_\lambda(0)D_u^0(\lambda,\nu_j(\lambda))}c_j(\lambda,y)\me^{\nu_-^0(\lambda)x}\me^{-\nu_j(\lambda)y}\int_x^\infty \me^{(-\nu_+^0(\lambda)+\nu_j(\lambda))\tau}(1+\theta_+(\tau,\lambda)) h(\tau) \md\tau. 
 \end{align*}
The integral terms converge due to the gap lemma condition and are $\mathcal{O}(\me^{(-\nu_+^0(\lambda)+\nu_j(\lambda)-\vartheta)x})$ and consequently these terms can be absorbed into $\boldJ_{\lambda,j}(x,y)$.

We now work out the first integral which must be broken into three pieces:
\begin{align*}
\int_{-\infty}^x \frac{\varphi^-(\tau)}{\W_\lambda(\tau)} h(\tau)\bG_\lambda^{12,\infty}(\tau,y)\md \tau&=\int_{-\infty}^y \frac{\varphi^-(\tau)}{\W_\lambda(\tau)} h(\tau)\bG_\lambda^{12,\infty}(\tau,y)\md \tau+\int_y^0 \frac{\varphi^-(\tau)}{\W_\lambda(\tau)} h(\tau)\bG_\lambda^{12,\infty}(\tau,y)\md \tau\\
&~~+\int_0^x \frac{\varphi^-(\tau)}{\W_\lambda(\tau)} h(\tau)\bG_\lambda^{12,\infty}(\tau,y)\md \tau.
\end{align*}
For the first of these integrals we have to consider the term
\[ \frac{h_4(\lambda,y)}{\W_\lambda(0)}\int_{-\infty}^y \me^{(\nu_+^1(\lambda)-\nu_-^1(\lambda))\tau}(1+\theta_-(\tau,\lambda))h(\tau) \md \tau, \]
from which we observe that the integral is bounded and is $\mathcal{O}( \me^{(\nu_+^1(\lambda)-\nu_-^1(\lambda)+\vartheta)y})$ and we see that this contributions can be distributed among $\boldH_\lambda(x,y)$ and $\bI_{\lambda,j}(x,y)$.  The other terms from $\bG_\lambda^{12,\infty}(x,y)$ are those proportional to $\me^{\nu_j(\lambda)}$ for $j=3,4$.  For these the integral becomes 
\[ \beta \sum_{j=3}^4 \frac{c_j(\lambda)}{\W_\lambda(0)D_u^1(\lambda,\nu_j(\lambda))}\me^{-\nu_j(\lambda)y} \int_{-\infty}^y \me^{(\nu_j(\lambda)-\nu_-^1(\lambda))\tau}(1+\theta_-(\tau,\lambda))h(\tau) \md \tau, \]
for which convergence of the integral is automatic and these terms are incorporated into $\boldH_\lambda(x,y)$.

We now evaluate the second integral for which we need
\bqs
\bG_\lambda^{12,\infty}(\tau,y)=-h_2(\lambda,y)\me^{\nu_-^1(\lambda)\tau}-h_3(\lambda,y)\me^{\nu_+^1(\lambda)\tau}-\beta\sum_{j=1}^2\frac{c_j(\lambda)}{D_u^1(\lambda,\nu_j(\lambda))}\me^{\nu_j (\lambda)(\tau-y)}, \quad y<\tau <0.
\eqs

We must estimate the integral,
\[
\int_y^0 \me^{-\nu_-^1(\lambda) \tau }(1+\theta_-(\tau,\lambda))h(\tau)\bG_\lambda^{12,\infty}(\tau,y)\md \tau, \]
from which we note that the integral is bounded uniformly in $y$ and therefore this term can be incorporated into $\boldH_\lambda(x,y)$ and $\bI_\lambda(x,y)$ after consulting the formula for $h_2(\lambda,y)$ and $h_3(\lambda,y)$.  A similar decomposition occurs for the inhomogeneous terms proportional to $\me^{\nu_j(\lambda)(\tau-y)}$ .

This bring us to the final integral with $0<\tau<x$ for which we must use the form of $\bG_\lambda^{12,\infty}$ for $y<0<\tau$.  For $\varphi^-(\tau)$ we must use its representation valid for $\tau>0$ and we then need to estimate the integrals
\begin{align}
\frac{C(\lambda)}{\W_\lambda(0)}\int_0^x \me^{-\nu_-^0(\lambda) \tau }(1+\kappa_+(\tau,\lambda))h(\tau)\bG_\lambda^{12,\infty}(\tau,y)\md \tau \nonumber  \\
+\frac{D(\lambda)}{\W_\lambda(0)}\int_0^x \me^{-\nu_+^0(\lambda) \tau }(1+\theta_+(\tau,\lambda))h(\tau)\bG_\lambda^{12,\infty}(\tau,y)\md \tau. \end{align}
On this interval the Green's function $\bG_\lambda^{12,\infty}(x,y)$ has one term of the form $h_1(\lambda,y)\me^{\nu_-^0(\lambda)\tau}$ for which we see that the integrals are uniformly bounded in $x$ and $y$ and the contributions can be divided among $\boldH_\lambda(x,y)$ and $\bI_{\lambda,j}(x,y)$ following the dependence of $h_1(\lambda,y)$ on $y$.  The final terms in $\bG_\lambda^{12,\infty}(x,y)$ are those proportional to $\me^{\nu_j(\lambda)(\tau-y)}$.  In this case, the integrals have contributions of $\mathcal{O}(\me^{-\nu_-^0(\lambda)x }\me^{\nu_j(\lambda)(x-y)})$ and $\mathcal{O}(\me^{-\nu_j(\lambda)y})$.  Therefore, these terms contribute to $\bI_{\lambda,j}(x,y)$ and $\boldJ_{\lambda,j}(x,y)$.  This concludes the analysis in this case.

\paragraph{Case $0<y<x$.} As before, the pointwise Green's function for this arrangement can be expressed as a sum of homogeneous and particular solutions for the asymptotic system plus a correction term expressed by a variation of constants formula,
\begin{align*} \bG_\lambda^{12}(x,y) & =b_1(\lambda,y)\me^{\nu_-^0(\lambda)x}-\beta \sum_{j=1}^2\frac{c_j(\lambda)}{D_u^0(\lambda,\nu_j(\lambda))}\me^{\nu_j(\lambda)(x-y)}\\
&~~~+
\varphi^+(x) \int_{-\infty}^x  
\frac{\varphi^-(\tau)}{\W_\lambda(\tau)}h(\tau) \bG_\lambda^{12,\infty}(\tau,y)\md \tau +\varphi^-(x) \int_x^\infty  
\frac{\varphi^+(\tau)}{\W_\lambda(\tau)}h(\tau) \bG_\lambda^{12,\infty}(\tau,y)\md \tau.
\end{align*}
To obtain bounds, the two integrals must once again be decomposed into four integrals depending on the relative positions of $\tau$, $y$ and zero.  The first such integral is 
\[ \varphi^+(x)\int_{-\infty}^0 \frac{\varphi^-(\tau)}{\W_\lambda(\tau)}h(\tau) \bG_\lambda^{12,\infty}(\tau,y)\md \tau.\]
Here we must use the expression for the asymptotic Green's function valid for $\tau<0<y$, 
\[  \bG_\lambda^{12,\infty}(\tau,y)=b_4(\lambda,y)\me^{\nu_+^1(\lambda)\tau}+\frac{\beta}{D_u^1(\lambda,\nu_3(\lambda))}c_3(\lambda) \me^{\nu_3(\lambda)(\tau-y)}+\frac{\beta}{D_u^1(\lambda,\nu_4(\lambda))}c_4(\lambda) \me^{\nu_4 (\lambda)(\tau-y)}. \]
Convergence of the integral is guaranteed since
\[ \mathrm{Re}\left(\nu+\vartheta-\nu_-^1(\lambda)\right)>0,\quad \text{for} \ \nu=\nu_+^1(\lambda), \nu_3(\lambda), \text{and} \ \nu_4(\lambda)  \]
and the integral inherits the domain of analyticity of the integrand.  Focusing first on the term involving $b_4(\lambda,y)$ we can write this as 
\[ \me^{\nu_-^0(\lambda)(x-y)} \left( \frac{b_4(\lambda,y)}{\W_\lambda(0)} \me^{\nu_-^0(\lambda)y} (1+\theta^+(x,\lambda))\int_{-\infty}^0 \me^{-\nu_-^1(\lambda) \tau} (1+\theta^-(\tau,\lambda)) h(\tau) \me^{\nu_+^1(\lambda) \tau } \md\tau\right), \]
Note that $b_4(\lambda,y)$ has singularities when $\nu_2(\lambda)=\nu_3(\lambda)$ (i.e. $\lambda=\lambda_v^{\textnormal{bp}}$) or $\nu_+^0(\lambda)=\nu_2(\lambda)$ (i.e. $\lambda=\lambda^{\textnormal{rp}}$).  Factoring out these singularities and consulting the expression for  $b_4(\lambda,y)$ we observe that the remaining terms in the parenthesis are uniformly bounded in $x$ and $y$ and contribute to $\boldH_\lambda(x,y)$.

Next consider the terms involving $\nu_{j}(\lambda)$ with $j=3,4$.  In each case we can repeat the argument above and obtain
\[ \me^{\nu_-^0(\lambda)(x-y)} \left( \frac{\beta (1+\theta^+(x,\lambda))c_j(\lambda)}{\W_\lambda(0)D_u^1(\lambda,\nu_j(\lambda))}  \me^{\nu_-^0(\lambda)y}\me^{-\nu_j(\lambda)y}\int_{-\infty}^0 \me^{-\nu_-^1(\lambda) \tau} (1+\theta^-(\tau,\lambda)) h(\tau) \me^{\nu_j(\lambda) \tau } \md\tau\right), \]
Convergence of the integral follows due to the spectral gap for the asymptotic system at $-\infty$.  Factoring out the singularity at the branch point where $\nu_2(\lambda)=\nu_3(\lambda)$ we obtain a similar bound.  

The second integral is 
\[ \varphi^+(x)\int_0^y \frac{\varphi^-(\tau)}{\W_\lambda(\tau)}h(\tau) \bG_\lambda^{12,\infty}(\tau,y)\md \tau.\]
To analyze this integral we must replace $\varphi^-(\tau)$ by its representation for $\tau>0$ and use the expression for the asymptotic Green's function valid for $0<\tau<y$,
\[ \bG_\lambda^{12,\infty}(\tau,y)=-b_2(\lambda,y)\me^{\nu_-^0(\lambda)\tau}-b_3(\lambda,y)\me^{\nu_+^0(\lambda)\tau}+\beta \sum_{j=3}^4\frac{c_j(\lambda)}{D_u^0(\lambda,\nu_1(\lambda))}\me^{\nu_j(\lambda)(x-y)}.\]
Recall Remark~\ref{rmk:nosingularity} which says that the singularity at $\lambda=0$ due to the expressions for $b_2(\lambda,y)$ and $b_3(\lambda,y)$ is removable in the expression for $\bG_\lambda^{12}$.  To obtain estimates, we begin with the terms involving  $b_3(\lambda,y)$ for which we can factor
\begin{align*}  -\me^{\nu_-^0(\lambda)(x-y)} \left( \frac{C(\lambda)}{\W_\lambda(0)} b_3(\lambda,y) \me^{\nu_-^0(\lambda)y}(1+\theta^+(x,\lambda))\int_0^y \me^{(\nu_+^0(\lambda)-\nu_-^0(\lambda))\tau} (1+\kappa^+(\tau,\lambda)) h(\tau) \md\tau \right) \\ 
-\me^{\nu_-^0(\lambda)(x-y)} \left( \frac{D(\lambda)}{\W_\lambda(0)} b_3(\lambda,y) \me^{\nu_-^0(\lambda)y}(1+\theta^+(x,\lambda))\int_0^y (1+\theta^+(\tau,\lambda)) h(\tau) \md\tau \right)
\end{align*}
whereas for $b_2(\lambda,y)$ we have
\begin{align*} -\me^{\nu_-^0(\lambda)(x-y)} \left( \frac{C(\lambda)}{\W_\lambda(0)} b_2(\lambda,y) \me^{\nu_-^0(\lambda)y}(1+\theta^+(x,\lambda))\int_0^y  (1+\kappa^+(\tau,\lambda)) h(\tau) \md\tau \right), \\
-\me^{\nu_-^0(\lambda)(x-y)} \left( \frac{D(\lambda)}{\W_\lambda(0)} b_2(\lambda,y) \me^{\nu_-^0(\lambda)y}(1+\theta^+(x,\lambda))\int_0^y \me^{(\nu_-^0(\lambda)-\nu_+^0(\lambda))\tau}  (1+\theta^+(\tau,\lambda)) h(\tau) \md\tau \right).
\end{align*}
Inspection of the integrands in both cases reveals that the integrals are uniformly bounded in $y$ and by consulting the formula for $b_2(\lambda)$ we obtain that these terms can be factored into $\boldH_\lambda(x,y)$.  
The remaining terms involving $\nu_j(\lambda)$ for $j=3,4$ have a similar form
\begin{align*} \me^{\nu_-^0(\lambda)(x-y)} \left( \frac{\beta C(\lambda)(1+\theta^+(x,\lambda))c_j(\lambda)}{\W_\lambda(0)D_u^0(\lambda,\nu_j(\lambda))}  \me^{\nu_-^0(\lambda)y} \me^{-\nu_j(\lambda)y} \int_0^y \me^{(\nu_j(\lambda)-\nu_-^0(\lambda))\tau)} (1+\kappa^+(\tau,\lambda)) h(\tau) \md\tau \right), \\
\me^{\nu_-^0(\lambda)(x-y)} \left( \frac{\beta D(\lambda) (1+\theta^+(x,\lambda))c_j(\lambda)}{\W_\lambda(0)D_u^0(\lambda,\nu_j(\lambda))} \me^{\nu_-^0(\lambda)y}\me^{-\nu_j(\lambda)y}\int_0^y \me^{(\nu_j(\lambda)-\nu_+^0(\lambda))\tau}  (1+\theta^+(\tau,\lambda)) h(\tau) \md\tau \right).
\end{align*}
Since $\mathrm{Re}(\nu_{3,4}(\lambda)>0$, upon moving the exponentials involving $y$ into the integral we observe once again that the terms in the parenthesis are bounded and can be absorbed into $\boldH_\lambda(x,y)$.

The third integral is 
\[ \varphi^+(x)\int_y^x \frac{\varphi^-(\tau)}{\W_\lambda(\tau)}h(\tau) \bG_\lambda^{12,\infty}(\tau,y)\md \tau.\]
In this case, we have $0<y<\tau$ and the asymptotic Green's function takes the form
\[ \bG_\lambda^{12,\infty}(\tau,y)=b_1(\lambda,y)\me^{\nu_-^0(\lambda)\tau}-\frac{\beta}{D_u^0(\lambda,\nu_1(\lambda))}c_1(\lambda) \me^{\nu_1(\lambda)(\tau-y)}-\frac{\beta}{D_u^0(\lambda,\nu_2(\lambda))}c_2(\lambda) \me^{\nu_2 (\lambda)(\tau-y)}. \]
For the terms involving $b_1(\lambda,y)$ we have
\begin{align*} \me^{\nu_-^0(\lambda)(x-y)} \left( \frac{C(\lambda)}{\W_\lambda(0)} b_1(\lambda,y) \me^{\nu_-^0(\lambda)y}(1+\theta^+(x,\lambda))\int_y^x  (1+\kappa^+(\tau,\lambda)) h(\tau) \md\tau \right), \\
+\me^{\nu_-^0(\lambda)(x-y)} \left( \frac{D(\lambda)}{\W_\lambda(0)} b_1(\lambda,y) \me^{\nu_-^0(\lambda)y}(1+\theta^+(x,\lambda))\int_y^x \me^{(\nu_-^0(\lambda)-\nu_+^0(\lambda))\tau}  (1+\theta^+(\tau,\lambda)) h(\tau) \md\tau \right).
\end{align*}

The analysis here resembles the previous case and these terms can be incorporated into $\boldH_\lambda(x,y)$.  On the other hand, the terms involving $\nu_{1,2}(\lambda)$ require a bit more effort.  To begin we write
\begin{align*}
\me^{\nu_-^0(\lambda)x}(1+\theta^+(x))\frac{\beta C(\lambda) c_j(\lambda)}{\W_\lambda(0)D_u^0(\lambda,\nu_j(\lambda))}\int_y^x \me^{-\nu_-^0(\lambda) \tau} (1+\kappa^+(\tau))h(\tau) \me^{\nu_j(\lambda)(\tau-y)} \md \tau ,\\
+\me^{\nu_-^0(\lambda)x}(1+\theta^+(x))\frac{\beta D(\lambda) c_j(\lambda)}{\W_\lambda(0)D_u^0(\lambda,\nu_j(\lambda))}\int_y^x \me^{-\nu_+^0(\lambda) \tau} (1+\theta^+(\tau))h(\tau) \me^{\nu_j(\lambda)(\tau-y)} \md \tau.
\end{align*}
For the first integral, we use that $\mathrm{Re}(\nu_{1,2}(\lambda))<0$ while $\mathrm{Re}(-\nu_-^0(\lambda)-\vartheta)>0$ and see that these terms can be incorporated into $\boldH_\lambda(x,y)$.  For the second integral, the gap condition on the eigenvalues implies that the integral can be bounded by the exponentials evaluated at $\tau=y$ from which we see that these terms can be absorbed into $\boldH_\lambda(x,y)$.  

The final integral is
\[ \varphi^-(x)\int_x^\infty \frac{\varphi^+(\tau)}{\W_\lambda(\tau)}h(\tau) \bG_\lambda^{12,\infty}(\tau,y)\md \tau.\]
The same expression for the asymptotic Green's function holds here since $0<y<\tau$.   The terms involving $b_1(\lambda,y)$ are rather straightforward to handle, so we focus instead on those involving $\nu_j(\lambda)$. Once again decomposing $\varphi^-(x)$ into an expression valid for $x>0$, we find terms including
\[ C(\lambda)\me^{\nu_+^0(\lambda)x}(1+\kappa^+(x)) \frac{-\beta c_j(\lambda)}{\W_\lambda(0)D_u^0(\lambda,\nu_j(\lambda))} \int_x^\infty \me^{-\nu_+^0(\lambda)\tau}(1+\theta^+(\tau))h(\tau)\me^{\nu_j(\lambda)(\tau-y)} \md \tau. \]
The integral converges provided due to the gamma lemma condition.  After re-arranging we obtain
\[ \me^{\nu_j(\lambda)(x-y)}\left( \frac{-\beta C(\lambda) c_j(\lambda)}{\W_\lambda(0)D_u^0(\lambda,\nu_j(\lambda))} (1+\kappa^+(x))  \me^{(\nu_+^0(\lambda)-\nu_j(\lambda)) x} \int_x^\infty
\me^{-\nu_+^0(\lambda)\tau}(1+\theta^+(\tau))h(\tau)\me^{\nu_j(\lambda)\tau}\right). \]
The terms in the parenthesis are bounded by $Ce^{-\vartheta x}$ and analytic aside from singularities occurring when $\nu_+^0(\lambda)=\nu_2(\lambda)$ or $\nu_2(\lambda)=\nu_3(\lambda)$.  They contribute to $\boldJ_{\lambda,j}(x,y)$.

\section{Absolute spectrum of the Swift-Hohenberg equation}\label{sec:SwiftHohenbergAbsSpec}
We prove Lemma~\ref{l:shabs} abd provide explicit formulas for the absolute spectrum of the Swift-Hohenberg equation. 
We utilize $\Sess^\eta(\calL^+_v)$ introduced in Section~\ref{sec:ExponentialWeights} and therefore recap the weighted spectral curve
\[\sigma_v(k;\eta)=-k^4+4\mbi \eta k^3 +(2+6\eta^2)k^2 + \mbi(s-4\eta^3-4\eta)k-(1+\eta^2)^2+s\eta+\mu,\]
with $\eta<0$. Also, recall the ordering~\eqref{eq:rootsvorderingabs} as well as the definition~\eqref{e:mi}, so that we can write
\[\Sabs(\calL^+_v)\coloneqq\lbrace \lambda\in\C: D_v(\lambda,\rho(\lambda))=0, \Re(\rho_2(\lambda))=\Re(\rho_3(\lambda)) \rbrace.\]

The proof of Lemma~\ref{l:shabs} is divided into three parts: (i) computation of the real triple point and showing that all real points to the left (on the real line) of it are elements of $\Sabs(\calL^+_v)$. (ii) computation of the two simple and pinched double roots. (iii) deriving an explicit formula for two curves that connect the triple point and the double roots and which are monotone with respect to the real part.

\begin{Proof}[ of Lemma~\ref{l:shabs}.]

\textbf{(i)} First, we obtain that a self-intersection of (exactly) three segments of $\sigma_v$ for different wavenumbers $k\in \R$ can only occur on the real axis. Second, since we look for a spectral point, later denoted by $\lambda^\textnormal{tr}_v$, for which solving $D(\lambda^\textnormal{tr}_v,\rho)=0$ leads to solutions that satisfy $\Re(\rho_j(\lambda^\textnormal{tr}_v))=\Re(\rho_{j+1}(\lambda^\textnormal{tr}_v))=\Re(\rho_{j+2}(\lambda^\textnormal{tr}_v))$ for an index $j\in \lbrace 1,2\rbrace$, we check that $\sigma_v(k;\eta)$ intersects the real axis for non-zero wavenumbers at
\[\overline{\sigma}_v(\eta)\coloneqq\sigma_v(\pm\sqrt{\eta^2+1-s/(4\eta)};\eta)=4\eta^4+4\eta^2+\mu-\frac{s^2}{16\eta^2},\]
where $\eta^2+1-s/(4\eta)>0$ holds for any $s>0$ and $\eta<0$. Moreover, the intersection point $\overline{\sigma}_v$ tends to $-\infty$ for $\eta\to 0$ and to $+\infty$ for $\eta\to -\infty$.

Since two segments of weighted essential spectrum intersect for different (non-zero) wavenumbers at $\overline{\sigma}_v$ for any given $\eta<0$, it holds that the real part of (at least) two spatial roots $\rho_j$ of $D_v(\overline{\sigma}_v,\rho)$ are equal to $\eta$, i.e.\ it exists an index $j\in\lbrace 1,\dots,3\rbrace$ such that $\Re(\rho_j(\overline{\sigma}_v))=\eta=\Re(\rho_{j+1}(\overline{\sigma}_v))$.
Next, we find that $\Im(\sigma_v(0;\eta))=0$ and obtain
\[\sigma_v^0(\eta)\coloneqq\sigma_v(0;\eta)=-\eta^4-2\eta^2+s\eta-1+\mu,\]
where $\sigma_v^0<0$ for $\eta<0$ and any given $s>0$ as well as $\mu<0$. 

Since $\overline{\sigma}_v$ tends from $-\infty$ to $+\infty$ for $\eta$ from zero to $-\infty$ and at the same time $\sigma_v^0<0$ holds, we expect that there exists a weight $\eta<0$ such that $\overline{\sigma}_v=\sigma^0_v$. Thus, we compute that $\overline{\sigma}_v=\sigma^0_v$ whenever $\eta^\textnormal{tr}_v$  is the (unique) real solution to $20\eta^3+4\eta+s=0$. We have $\eta^\textnormal{tr}_v<0$ holds for any $s>0$.

The key point is that at $\lambda^\textnormal{tr}_v\coloneqq\overline{\sigma}_v(\eta^\textnormal{tr}_v)=\sigma^0_v(\eta^\textnormal{tr}_v)$, the spatial roots of $D_v(\lambda^\textnormal{tr}_v, \rho)$ satisfy $\Re(\rho_j(\lambda^\textnormal{tr}_v))=\Re(\rho_{j+1}(\lambda^\textnormal{tr}_v))=\Re(\rho_{j+2}(\lambda^\textnormal{tr}_v))$ for an index $j \in \lbrace 1,2\rbrace$ and hence $\lambda^\textnormal{tr}_v$ is a triple point in $\Sabs(\calL^+_v)$ (recall $i_\infty^v=2$), where we refer to~\cite{rademacher2006geometric, rademacher07} for more information on singularities of the absolute spectrum.

The explicit formula for the triple point is given by~\eqref{eq:TriplePointSwiftHohenberg} via $\sigma^0_v(\eta^\textnormal{tr}_v)$ and we note that $\lambda^\textnormal{tr}_v<-1$ for any $s>0$ and $\mu<0$. Additionally, since $\lambda^\textnormal{tr}_v$ lies to the left of the rightmost real point of $\Sess(\calL^+_v)$, it actually holds that
\[\Re(\rho_1(\lambda^\textnormal{tr}_v))=\Re(\rho_2(\lambda^\textnormal{tr}_v))=\Re(\rho_3(\lambda^\textnormal{tr}_v))<0<\Re(\rho_4(\lambda^\textnormal{tr}_v)).\]

Next, we use~\cite[Proposition 2.3.1]{fiedler03} which states that the (spatial) Morse index increases or decreases by one upon crossing $\sigma_v$ in the complex plane from left to right depending on its orientation. The orientations locally around $\overline{\sigma}_v$ for all $s>0$ and $\eta^\textnormal{tr}_v\le\eta<0$ are
\begin{equation}\label{eq:OrientationSwiftHohenbergTriplePoint}
\begin{split}
\Re(\sigma'_v(-\sqrt{\eta^2+1-s/(4\eta)};\eta)&>0,\\
\Re(\sigma'_v(+\sqrt{\eta^2+1-s/(4\eta)};\eta)&<0,\\
\Im(\sigma'_v(\pm\sqrt{\eta^2+1-s/(4\eta)};\eta)&<0,
\end{split}
\end{equation}
where prime denotes the derivative with respect to the wavenumber $k$. From this, we obtain that $\overline{\sigma}_v\in \Sabs(\calL^+_v)$ for $\eta^\textnormal{tr}_v<\eta<0$ since in this case it holds that
\[\Re(\rho_1(\overline{\sigma}_v))\le\Re(\rho_2(\overline{\sigma}_v))=\Re(\rho_3(\overline{\sigma}_v))<0<\Re(\rho_4(\overline{\sigma}_v))\]
and thus the only (purely) real points in $\Sabs(\calL^+_v)$ are elements of the interval $(-\infty, \lambda^\textnormal{tr}_v]$, with $\lambda^\textnormal{tr}_v<-1$.

\medskip
\textbf{(ii)} Utilizing the resultant function, we compute that two of the the three (candidates) for double roots in $\Sabs(\calL^+_v)$ are given by~\eqref{eq:DoubleRootsSwiftHohenberg}. Note that the third candidate for a double root is purely real and always larger than $\lambda^\textnormal{tr}_v$. Based on step (i), it is not an element of $\Sabs(\calL^+_v)$ and thus we neglect it in what follows. We further find $\lambda^\textnormal{tr}_v<\Re(\lambda^\textnormal{dr}_v)$ for all $s>0$ and $\mu<0$ and since $\Im(\lambda^\textnormal{tr}_v)=0$ but $\Im(\lambda^\textnormal{dr}_v)\neq 0$, we expect the emergence of the two complex conjugated branches from $\lambda^\textnormal{tr}_v$.

Next, we show that $\lambda^\textnormal{dr}_v$ are indeed elements of $\Sabs(\calL^+_v)$ and that they correspond to the rightmost points. Since $\Re(\sigma'_v)=0$ for $k^r_\pm(\eta)\coloneqq\pm\sqrt{3\eta^2+1}$ and $\Im(\sigma'_v)=0$ for \[k^i_\pm(\eta)\coloneqq\pm\frac{\sqrt{4\eta^3+4\eta-s}}{\sqrt{12\eta}},\]
we obtain from $k^r_\pm= k^i_\pm$ the real solution
\[\eta^\textnormal{dr}_v\coloneqq\frac{-4\sqrt[3]{9}+\sqrt[3]{3}\kappa^2}{12\kappa}<0,\quad \textnormal{where}\quad \kappa\coloneqq \sqrt[3]{-9s+\sqrt{192+81s^2}}>0.\]
Note that the real solution to $k^r_\pm=k^i_\pm$ is unique due to symmetry and we further define
\[k^r_\pm(\eta^\textnormal{dr}_v)=k^i_\pm(\eta^\textnormal{dr}_v)=k^\textnormal{dr}_\pm\coloneqq\pm\sqrt{\frac{1}{2}+\frac{\sqrt[3]{3}}{\kappa^2}+\frac{\kappa^2}{16\sqrt[3]{3}}}.\]
Therefore, the spectral curve $\sigma_v$ possesses two complex conjugated cusps for $\eta=\eta^\textnormal{dr}_v$ which are located at $\sigma_v(k^\textnormal{dr}_\pm;\eta^\textnormal{dr}_v)=\lambda_v^\textnormal{dr}$. Using again~\cite[Proposition 2.3.1]{fiedler03}, we obtain that the double roots $\lambda^\textnormal{dr}_v$ are the rightmost points in $\Sabs(\calL^+_v)$ and thus they are pinched. Moreover, one readily verifies that these double roots are also simple.

\medskip
\textbf{(iii)} First, we note that $\eta^\textnormal{tr}_v<\eta^\textnormal{dr}_v<0$ for all $s>0$ and independent of $\mu$. Next, we obtain again by~\cite[Proposition 2.3.1]{fiedler03} and~\cite[Lemma 4.2]{rademacher2006geometric}, where the latter one describes the Morse index at intersections of two curves in $\Sess^\eta$, that further elements of $\Sabs(\calL^+_v)$ apart from the interval $(-\infty,\lambda^\textnormal{tr}_v]$ (see step (i)) occur only at self-intersections of $\sigma_v$ for $\eta^\textnormal{tr}_v\le\eta\le\eta^\textnormal{dr}_v$. These intersections are located at points for which for any wavenumbers $k_1$ and $k_2$ with $k_1\neq k_2$ it holds that $\sigma_v(k_1;\tilde{\eta})=\sigma_v(k_2;\tilde{\eta})$ with $\eta^\textnormal{tr}_v\le\tilde{\eta}\le\eta^\textnormal{dr}_v$. Equating the real and imaginary parts, we obtain that under the condition $\eta^\textnormal{tr}_v\le\tilde{\eta}\le\eta^\textnormal{dr}_v$ the intersections are located at
\[\check{\sigma}_v(\tilde{\eta})\coloneqq 24\tilde{\eta}^4+8\tilde{\eta}^2+\frac{7}{2}s\tilde{\eta}+\mu+\frac{8s\tilde{\eta}+s^2}{16\tilde{\eta}^2}\pm\mbi\left(-4\tilde{\eta}^3-\frac{s}{2}+\frac{1}{2}\sqrt{-(8\tilde{\eta}^3+s)(32\tilde{\eta}^3+8\tilde{\eta}+s)}\right)\tilde{\kappa},\]
where
\[\tilde{\kappa}\coloneqq \sqrt{12\tilde{\eta}^2+4-\frac{1}{\tilde{\eta}}\sqrt{-(8\tilde{\eta}^3+s)(32\tilde{\eta}^3+8\tilde{\eta}+s)}}.\]
One further verifies that $\check{\sigma}_v(\eta^\textnormal{tr}_v)=\lambda^\textnormal{tr}_v$ and $\check{\sigma}_v(\eta^\textnormal{dr}_v)=\lambda^\textnormal{dr}_v$, as expected.

The last step is to show that $\Re(\check{\sigma}'_v)>0$ for all $\eta^\textnormal{tr}_v<\tilde{\eta}<\eta^\textnormal{dr}_v$. To this end, we first obtain that $\sqrt{-(8\tilde{\eta}^3+s)(32\tilde{\eta}^3+8\tilde{\eta}+s)}$ is positive for all $\eta^\textnormal{tr}_v<\tilde{\eta}<\eta^\textnormal{dr}_v$, since $8\tilde{\eta}^3+s$ is positive and $32\tilde{\eta}^3+8\tilde{\eta}+s$ is negative for those weights (recall that $\eta^\textnormal{tr}_v$ is the real root of $20\eta^3+4\eta+s$), and where we note that $\eta^\textnormal{dr}_v$ solves precisely $32\eta^3+8\eta+s=0$. Thus, the real part of $\check{\sigma}'_v$ reads
\[\Re(\check{\sigma}'_v(\tilde{\eta}))=96\tilde{\eta}^3+16\tilde{\eta}+\frac{7}{2}s+\frac{s}{2\tilde{\eta}^2}-\frac{8s \tilde{\eta}+s^2}{8\tilde{\eta}^3}\]
for which we obtain that it increases monotonically for any $\tilde{\eta}<0$ and that $\Re(\check{\sigma}'_v(\eta^\textnormal{tr}_v))>0$ as well as $\Re(\check{\sigma}'_v(\eta^\textnormal{dr}_v))>0$.

Summing up, the absolute spectrum of $\calL^+_v$ is given by the real interval $(-\infty,\lambda^\textnormal{tr}_v]$, where $\lambda^\textnormal{tr}_v<-1$, and the two complex conjugated branches given by $\check{\sigma}_v(\tilde{\eta})$ for $\eta^\textnormal{tr}_v\le\tilde{\eta}\le\eta^\textnormal{dr}_v$, connecting $\lambda^\textnormal{tr}_v$ to the left and $\lambda^\textnormal{dr}_v$ to the right, and where the real part of these branches increases strictly from left to right. This concludes the proof.
\end{Proof}

We emphasize again that the branch points $\lambda^\textnormal{dr}_v$, i.e.\ the simple and pinched double roots, are the most unstable elements of $\Sabs(\calL^+_v)$ for any $s>0$ and $\mu$. Based on the monotonicity of $\check{\sigma}(\tilde{\eta})$ with respect to $\tilde{\eta}$ it follows that remnant, absolute, and pointwise instability coincide in the Swift-Hohenberg, which is also the case for the KPP equation.

\end{appendix}

\bibliographystyle{abbrv}
\bibliography{RemnantBib}

\def\cprime{$'$}
\begin{thebibliography}{10}

\bibitem{AW78}
D.~G. Aronson and H.~F. Weinberger.
\newblock Multidimensional nonlinear diffusion arising in population genetics.
\newblock {\em Advances in Mathematics}, 30(1):33--76, 1978.

\bibitem{beck14}
M.~Beck, T.~T. Nguyen, B.~Sandstede, and K.~Zumbrun.
\newblock Nonlinear stability of source defects in the complex
  {G}inzburg-{L}andau equation.
\newblock {\em Nonlinearity}, 27(4):739--786, 2014.

\bibitem{bers84}
A.~{Bers}.
\newblock {Space-time evolution of plasma instabilities-absolute and
  convective}.
\newblock In {A.~A.~Galeev \& R.~N.~Sudan}, editor, {\em Basic Plasma Physics:
  Selected Chapters, Handbook of Plasma Physics, Volume 1}, pages 451--517,
  1984.

\bibitem{brevdo96}
L.~Brevdo and T.~J. Bridges.
\newblock Absolute and convective instabilities of spatially periodic flows.
\newblock {\em Philos. Trans. Roy. Soc. London Ser. A}, 354(1710):1027--1064,
  1996.

\bibitem{bricmont92}
J.~Bricmont and A.~Kupiainen.
\newblock Renormalization group and the {G}inzburg-{L}andau equation.
\newblock {\em Comm. Math. Phys.}, 150(1):193--208, 1992.

\bibitem{briggs}
R.~J. Briggs.
\newblock {\em Electron-Stream Interaction with Plasmas}.
\newblock MIT Press, Cambridge, 1964.

\bibitem{eckmann94}
J.-P. Eckmann and C.~E. Wayne.
\newblock The nonlinear stability of front solutions for parabolic partial
  differential equations.
\newblock {\em Comm. Math. Phys.}, 161(2):323--334, 1994.

\bibitem{faye19}
G.~Faye and M.~Holzer.
\newblock Asymptotic stability of the critical {F}isher-{KPP} front using
  pointwise estimates.
\newblock {\em Z. Angew. Math. Phys.}, 70(1):Art. 13, 21, 2019.

\bibitem{faye20}
G.~Faye and M.~Holzer.
\newblock Asymptotic stability of the critical pulled front in a
  {L}otka-{V}olterra competition model.
\newblock {\em Journal of Differential Equations}, 269(9):6559 -- 6601, 2020.

\bibitem{faye17}
G.~Faye, M.~Holzer, and A.~Scheel.
\newblock Linear spreading speeds from nonlinear resonant interaction.
\newblock {\em Nonlinearity}, 30(6):2403--2442, may 2017.

\bibitem{fiedler03}
B.~Fiedler and A.~Scheel.
\newblock Spatio-temporal dynamics of reaction-diffusion patterns.
\newblock In {\em Trends in nonlinear analysis}, pages 23--152. Springer,
  Berlin, 2003.

\bibitem{gallay94}
T.~Gallay.
\newblock Local stability of critical fronts in nonlinear parabolic partial
  differential equations.
\newblock {\em Nonlinearity}, 7(3):741--764, 1994.

\bibitem{gallay04}
T.~Gallay, G.~Schneider, and H.~Uecker.
\newblock Stable transport of information near essentially unstable localized
  structures.
\newblock {\em Discrete and Continuous Dynamical Systems-Series B},
  (2):349--390, 2004.

\bibitem{gardner98}
R.~A. Gardner and K.~Zumbrun.
\newblock The gap lemma and geometric criteria for instability of viscous shock
  profiles.
\newblock {\em Comm. Pure Appl. Math.}, 51(7):797--855, 1998.

\bibitem{ghazaryan07}
A.~Ghazaryan and B.~Sandstede.
\newblock Nonlinear convective instability of {T}uring - unstable fronts near
  onset: A case study.
\newblock {\em SIAM Journal on Applied Dynamical Systems}, 6(2):319--347, 2007.

\bibitem{goh11}
R.~N. Goh, S.~Mesuro, and A.~Scheel.
\newblock Spatial wavenumber selection in recurrent precipitation.
\newblock {\em SIAM Journal on Applied Dynamical Systems}, 10(1):360--402,
  2011.

\bibitem{holzer14}
M.~Holzer.
\newblock Anomalous spreading in a system of coupled {F}isher-{KPP} equations.
\newblock {\em Phys. D}, 270:1--10, 2014.

\bibitem{holzer16}
M.~Holzer.
\newblock A proof of anomalous invasion speeds in a system of coupled
  {F}isher-{KPP} equations.
\newblock {\em Discrete Contin. Dyn. Syst.}, 36(4):2069--2084, 2016.

\bibitem{holzerscheelLV}
M.~Holzer and A.~Scheel.
\newblock A slow pushed front in a lotka--volterra competition model.
\newblock {\em Nonlinearity}, 25(7):2151, 2012.

\bibitem{holzerscheel14}
M.~Holzer and A.~Scheel.
\newblock Criteria for pointwise growth and their role in invasion processes.
\newblock {\em J. Nonlinear Sci.}, 24(4):661--709, 2014.

\bibitem{howard02}
P.~Howard.
\newblock Pointwise estimates and stability for degenerate viscous shock waves.
\newblock {\em J. Reine Angew. Math.}, 545:19--65, 2002.

\bibitem{howard12}
P.~Howard and B.~Kwon.
\newblock Asymptotic stability analysis for transition front solutions in
  {C}ahn-{H}illiard systems.
\newblock {\em Phys. D}, 241(14):1193--1222, 2012.

\bibitem{huerre90}
P.~{Huerre} and P.~A. {Monkewitz}.
\newblock {Local and global instabilities in spatially developing flows}.
\newblock {\em Annual Review of Fluid Mechanics}, 22:473--537, 1990.

\bibitem{kapitula98}
T.~Kapitula and B.~Sandstede.
\newblock Stability of bright solitary-wave solutions to perturbed nonlinear
  {S}chr\"odinger equations.
\newblock {\em Phys. D}, 124(1-3):58--103, 1998.

\bibitem{kirchgassner92}
K.~Kirchgässner.
\newblock On the nonlinear dynamics of travelling fronts.
\newblock {\em Journal of Differential Equations}, 96(2):256 -- 278, 1992.

\bibitem{lunardi}
A.~Lunardi.
\newblock {\em Analytic semigroups and optimal regularity in parabolic
  problems}.
\newblock Modern Birkh\"{a}user Classics. Birkh\"{a}user/Springer Basel AG,
  Basel, 1995.
\newblock [2013 reprint of the 1995 original] [MR1329547].

\bibitem{rademacher2006geometric}
J.~D.~M. Rademacher.
\newblock Geometric relations of absolute and essential spectra of wave trains.
\newblock {\em SIAM J. Appl. Dyn. Syst.}, 5(4):634--649, 2006.

\bibitem{rademacher07}
J.~D.~M. Rademacher, B.~Sandstede, and A.~Scheel.
\newblock Computing absolute and essential spectra using continuation.
\newblock {\em Phys. D}, 229(2):166--183, 2007.

\bibitem{sandstede00}
B.~Sandstede and A.~Scheel.
\newblock Absolute and convective instabilities of waves on unbounded and large
  bounded domains.
\newblock {\em Phys. D}, 145(3-4):233--277, 2000.

\bibitem{sattinger}
D.~H. Sattinger.
\newblock On the stability of waves of nonlinear parabolic systems.
\newblock {\em Advances in Math.}, 22(3):312--355, 1976.

\bibitem{vansaarloos03}
W.~van Saarloos.
\newblock Front propagation into unstable states.
\newblock {\em Physics Reports}, 386(2-6):29 -- 222, 2003.

\bibitem{zumbrun11}
K.~Zumbrun.
\newblock Instantaneous shock location and one-dimensional nonlinear stability
  of viscous shock waves.
\newblock {\em Quart. Appl. Math.}, 69(1):177--202, 2011.

\bibitem{zumbrun98}
K.~Zumbrun and P.~Howard.
\newblock Pointwise semigroup methods and stability of viscous shock waves.
\newblock {\em Indiana Univ. Math. J.}, 47(3):741--871, 1998.

\end{thebibliography}

\end{document}